\documentclass[12pt]{article}
\usepackage{amsfonts}
\usepackage{amssymb}
\usepackage{polski}

\input{tcilatex}
\begin{document}

\begin{center}
\textbf{The Go\l \k{a}b-Schinzel and Goldie functional equations in Banach
algebras}

\textbf{by \\[0pt]
N. H. Bingham and A. J. Ostaszewski}\\[0pt]

\bigskip
\end{center}

\noindent \textbf{Abstract.} We are concerned below with the
characterization in a unital commutative real Banach algebra $\mathbb{A}$ of
continuous solutions of the Go\l \k{a}b-Schinzel functional equation
(below), the general Popa groups they generate and the associated Goldie
functional equation. This yields general structure theorems involving both
linear and exponential homogeneity in $\mathbb{A}$ for both these functional
equations and also explict forms, in terms of the recently developed theory
of multi-Popa groups [BinO3,4], both for the ring $C[0,1]$ and for the case
of $\mathbb{R}^{d}$ with componentwise product, clarifying the context of
recent developments in [RooSW]. The case $\mathbb{A=C}$ provides a new
viewpoint on continuous complex-valued solutions of the primary equation by
distinguishing analytic from real-analytic ones.

\bigskip

\noindent \textbf{Classification}: 26A03, 26A12, 33B99, 39B22.

\section{Equations and groups}

\subsection{Functional equations linked to that of Cauchy}

\textit{The Cauchy, Go\l \k{a}b-Schinzel, Goldie and Levi-Civita functional
equations. }General regular variation [BinO1] has recently emerged as
embracing three kinds of univariate regular variation (RV) due in turn to:
Karamata (classical), Bojani\'{c}-Karamata-de Haan, and Beurling (for the
Beurling Tauberian Theorem), for which see [BinGT]. Underlying this
unification is the Go\l \k{a}b-Schinzel functional equation and an
associated group structure (below). The equation reads%
\begin{equation}
S(x+S(x)y)=S(x)S(y)  \tag{$GS$}
\end{equation}%
with $x,y$ ranging over a half-line in $\mathbb{R}$ (`S for survival
probability') and has positive (continuous) solutions which necessarily take
the `canonical' form%
\begin{equation}
S(x)=S_{\rho }(x):=1+\rho x\qquad (\rho \geq 0),  \tag{$Can$}
\end{equation}%
so with a maximal connected domain $\mathbb{G}_{\rho }:=\{x:1+\rho
x>0\}=(-\rho ^{-1},\infty ):$ cf. Cor. 5.1. (For $\rho =0$, interpret $-\rho
^{-1}$ as $-\infty .)$ It is significant that the associated group operation 
$\circ $ on $\mathbb{G}_{\rho }$ (the Popa group, below) allows a
re-statement of $(GS)$ as a homomorphism equation%
\[
S(x\circ y)=S(x)S(y), 
\]%
so is fundamentally the multiplicative variant of the \textit{Cauchy
functional equation}. We study only \textit{continuous solutions,} leaving
aside the issues of automatic continuity of homomorphisms, for which see
e.g. [Ros], [Dal], cf. [Ost1, \S 1].

A further key tool in generalized regular variation is the \textit{Goldie
functional equation}:%
\begin{equation}
K(x+y)=K(x)+g(x)K(y)  \tag{$GFE$}
\end{equation}%
in the pair $(K,g)$ of real-valued functions; this is very closely related
to $(GS)$ (see [BinO1]). (Equation $(GFE)$ is a special case of a
Levi-Civita equation [Lev] -- see [Stet, Ch. 5]; cf. [AczD, Ch.14].) Here $g$
is necessarily multiplicative:%
\[
g(x+y)=g(x)g(y), 
\]%
and so is termed the \textit{multiplicative auxiliary. }It is either trivial
($g\equiv 1$) or exponential on $\mathbb{R}$. Correspondingly, $K$ is
additive and so linear, or in the non-trivial case (characterized by $g=1$
only for $x=0$), it is monotone (cf. [BinO4, Lemma 1]). Then, as $+$ is
commutative and $K,g$ are real-valued,%
\[
K(y)+g(y)K(x)=K(x)+g(x)K(y):\quad \lbrack g(y)-1]K(x)=[g(x)-1]K(y): 
\]
\[
K(x)=\text{const}\times (1-g(x))=c\frac{1-e^{-\gamma x}}{1-e^{-\gamma }}, %
\tag{$Exp$} 
\]%
with $c$ a constant. It is thematic here and below that $(Exp)$ reduces to $%
K(x)=cx$ for $\gamma =0$ under the \textit{L'Hospital convention}: we shall
see curvilinear variants of this linear-versus-exponential dichotomy in a
more general setting.

We note that any decreasing solution $g$ contributes a notable solution to $%
(GFE),$ namely%
\[
K(x)=1-e^{-\gamma x}\qquad \text{with }\gamma >0, 
\]%
one that is an exponential probability distribution on $[0,\infty );$ for
the background here, see [BinGT, Ch. 3] and [BinO1].

As pointed out in [Ost5] (cf. [Ost3]), assuming the solution $S$ of $(GS)$
to be injective, as will be the case for $S=S_{\rho }$ with $\rho >0$,
replacement of $S(x)$ by $e^{x}$ and of $S^{-1}(e^{x})$ by $K(x)$ yields $%
(GFE)$ with $g(x)=e^{x}.$ That is, $(GFE)$ and $(GS)$ are then formally
equivalent.

\bigskip

\textit{Infinite-dimensional settings. }Both $(GS)$ and $(GFE)$ are capable
of a natural interpretation when $x,y$ range over a topological vector space 
$X$ and the functions $S,K,g$ are real-valued. Here, by the Brillou\"{e}%
t-Dhombres-Brzd\k{e}k theorem [BriD, Prop. 3], [Brz1, Th. 4], the continuous
solutions of $(GS)$ take the form, for $x\in X$%
\begin{equation}
S(x)=S_{\rho }(x):=1+\rho (x)\qquad (\rho \in X^{\ast }),  \tag{$Can_{X}$}
\end{equation}%
with $X^{\ast }$ the dual space of continuous linear functionals on $X.$
Here too there are: an analogous maximal connected open domain 
\begin{equation}
\mathbb{G}_{\rho }(X):=\{x\in X:1+\rho (x)>0\} 
\tag{$\QTR{Bbb}{G}_{\rho
}(X)$}
\end{equation}%
and an associated abelian group structure $\circ _{\rho }$ on this domain.
In this context, Goldie's equation generalizes to 
\begin{equation}
K(x\circ _{\rho }y)=K(x)+g(x)K(y),  \tag{$GGE_{\rho g}$}
\end{equation}%
$g$ is again necessarily multiplicative on $\mathbb{G}_{\rho }$:%
\begin{equation}
g(x\circ _{\rho }y)=g(x)g(y),  \tag{$G$}
\end{equation}%
and, again by the commutativity of $\circ ,$ for some constant $c$%
\[
K(x)=c(1-g(x)). 
\]%
It further emerges [BinO3,4] that with $K$ injective, as here, $(GFE)$ may
be equivalently rewritten as%
\[
K(x+y)=K(x)\circ _{\sigma }K(y)\qquad (x,y\in X)\qquad \text{for }\sigma
(z):=g(K^{-1}(z)), 
\]%
with $\sigma =\sigma _{g}\in X^{\ast }.$ So again a Cauchy equation. Its
(continuous) solutions are explicitly characterized in [BinO3,4] in the more
general context%
\[
K(x\circ _{\rho }y)=K(x)\circ _{\sigma }K(y)\qquad (x,y\in \mathbb{G}_{\rho
}(X)) 
\]%
by reference to the fundamental homomorphisms: linear maps, exponentials,
and the two $GS$-type functions: $(1+\rho (x))$ and $(1+\sigma (x))$, these
being the basic building blocks.

\bigskip

\textit{Results. }Below we pursue both equations in the more general setting
of finite- and infinite-dimensional Banach algebras. There are six main
results in this general setting: the \textit{Decomposition Theorem}, Th. 2.1
(into linear and non-linear parts); \textit{Phantom Linear Characterization}%
, Th.2.2; \textit{Pencil Theorem}, Th. 4.1, showing that neighbourhoods of
the origin are spanned by a pencil through $0$ of otherwise disjoint
subgroups isomorphic to canonical (abelian) Popa groups; a \textit{First
Characterization Theorem}, Th. 5.1, from which it emerges that $(GS)$ and $%
(GGE_{SS}),$ another generalized form of $(GFE),$ are as inseparable as
husband and wife; a \textit{Second Characterization Theorem}, Th. 5.2,
showing that the solution of $(GFE)$ in a Banach algebra exhibits the
exponential homogeneity familiar in the real line setting $(Exp),$ albeit in
curvilinear form; finally, a \textit{Third Characterization Theorem}, Th.
5.3, giving under technical assumptions a differential characterization to a
region where $S$ takes the linear-plus-one form (`$1_{\mathbb{A}}+\mathrm{%
Linear}$').

As a preliminary, the analysis will begin with the finite-dimensional Banach
algebras provided by Euclidean space. Here the main result is the \textit{%
Structure Theorem}, Th. 3.3, where the general $d$-dimensional Popa group is
broken down into irreducible building blocks (as with the decomposition of
finite groups into finite simple groups). There is a consequent analogue for
the Banach algebra $C[0,1]$, based on the Stone-Weierstrass Theorem. The
signature of the Euclidean decomposition is the partition $\mathcal{P}$ of $%
\mathbb{N}_{d}=\{1,...,d\}$. This breaks down the $d$ coordinates of the
group elements into parts, whose coordinates are interchangeable with each
other but not with those from other parts. This `signature decomposition'
reveals structure at two levels for each factor, $i\in I$ and $I\in $ $%
\mathcal{P}$, and so three levels altogether (cf., say, `continents,
countries and counties'). Such structure may not be previously visible.

\bigskip

\textit{Connections with statistics. }The emergence of such `unsuspected
structure' can be very important, even in two dimensions. To give a classic
instance: in Silverman's book on density estimation [Sil, \S 4.2.3, Fig.
4.7] he gives an account of a study of a certain disease. A two-dimensional
contour plot of an estimated density revealed (in the manner of an Ordnance
Survey map) two `peaks'. Medical investigation showed that the disease under
study occurred in two forms, corresponding to these peaks. With this
difference identified, it emerged that the two forms were best treated in
different ways.

\bigskip

\textit{Connections with probability. }In probability theory, the theory of
independent sums is of central importance (infinitely-divisible laws, L\'{e}%
vy processes, L\'{e}vy-It\^{o} decomposition, L\'{e}vy-Khintchine formula,
etc.) A central role here is played by the stable laws -- those obtainable
from a single sequence of independent copies of a random variable, rather
than a doubly-indexed array of distributions (infinitely divisible laws) or
a singly-indexed one (self-decomposability). The role played by functional
equations in the study of stability directly (rather than by specialisation
from infinite-divisibility) has been considered by Pitman and Pitman [PitP]
and the second author [Ost4]. The functional equations relevant here are
those of Goldie, Cauchy and Levi-Civita.

In extreme-value theory (EVT), the role of sums above is played instead by
maxima (cf. [BinGT, \S 8.15]). A survey of the regular-variation aspects of
EVT was recently given by the authors ([BinO5]; cf. [BinO4]). It emerges
that the key functional equation there is the Goldie equation. It is
striking that the Goldie equation plays a central role in the theory of both
sums and maxima, two important areas whose similarities are striking but
whose differences are even more so.

This is in one dimension; in multidimensional situations in probability and
statistics, it is interesting to see the effect of dimensionality on how the
relevant limits are parametrised. In EVT in one dimension, the limits are
parametric (Fisher-Tippett theorem: three families classically, one if one
uses generalised extreme-value laws, GEV). But as soon as the dimension $d$
is at least two, limits become non-parametric (more precisely,
semi-parametric: one scalar radial parameter, one spectral measure on the
unit sphere). If one specialises to vines (cf. [BinO5, \S 2 Dependence
structure]), matters decompose into bivariate copulas (non-parametric),
linked by nested trees ($O(d)$ parameters). In Popa groups, there is no such
abrupt discontinuity as the dimension increases through $d=1,1<d<\infty $
and $d=\infty $, and it is striking that the link with Popa groups is lost
as soon as $d>1$. Here the multi-dimensional feature permits alternative
intepretations of $(GS)\ $according to the side-conditions (e.g.
collinearity or other co-dependencies) imposed on its two free variables:
see the comments preceding Prop. 1.1 below.

\subsection{Associated Groups and Banach algebras}

\textit{Popa groups. }Following Popa's analysis [Pop] of (Lebesgue)
measurable solutions of $(GS)$, we equip $\mathbb{G}_{\rho }\subseteq 
\mathbb{R}$ above with the operation%
\[
x\circ y=x\circ _{\rho }y:=x+(1+\rho x)y, 
\]%
turning $\mathbb{G}_{\rho }$ into a group, which under $S$ is isomorphic to $%
(\mathbb{R}_{+},\times )$. One may also follow Javor [Jav] by applying $%
\circ _{\rho }$ to $\mathbb{G}_{\rho }^{\ast }:=\{x:1+\rho x\neq 0\}=\mathbb{%
R}\backslash \{-\rho ^{-1}\}.$ We term $-\rho ^{-1}$ the \textit{Popa centre}%
.

This \textit{univariate} Popa group-structure provides the group theory,
previously lacking, with which to export transparently Karamata theory
(whose underlying group structure was explicitly recognized by Baj\v{s}%
anski-Karamata [BajK] and Balkema [Bal, Ch. 9]) to the other RV\ theories.

The equation $(GS),$ and likewise the corresponding group structure, refers
to the \textit{ring} structure of $\mathbb{R}$ and so extends to the context
of a unital commutative \textit{real} Banach algebra $\mathbb{A}$. Three
examples here are: the complex numbers $\mathbb{C}$, the setting for the
study of complex regularly varying functions, for which see [BinGT, A1.2]
(in Corollary 5.6 below we characterize the continuous solutions of $(GS_{%
\mathbb{C}})$); the Euclidean algebra $\mathbb{R}^{d}$ equipped with
componentwise (Hadamard) product, corresponding to the statistics of
sea-level measurement at $d$ locations -- see e.g. [RooSW], [KirRSW]; the
ring of continuous functions $C[0,1],$ corresponding to measurements of
sea-levels along a coastline parametrized by $[0,1]$. The latter two cases,
characterized in \S 3, provide a setting for the \textit{location and scale}
standardization of their statistics: the vector-space structure allows for
the translation (re-location) of each component random variable according to
its mean, together with \textit{uniform} scaling (i.e. a scale common to all
components); furthermore, the componentwise product structure permits 
\textit{individual} scaling of each component of a random variable or
stochastic process by its variance (equivalently its precision).

For a unital Banach algebra $\mathbb{A}$, we denote by $\mathbb{A}^{-1}$ the
open subset of invertible elements of $\mathbb{A}$, viewed as a
multiplicative group, and by $\mathbb{A}_{1}$ the connected component of the
identity (the \textit{principal component} [Ric, Def. 1.4.9]), a
multiplicative subgroup of $\mathbb{A}^{-1}$, coinciding with $\exp (\mathbb{%
A}),$ the exponential elements [Rud, 10.34]. For a solution $S:\mathbb{A}%
\rightarrow \mathbb{A}$ of $(GS),$ we equip the sets%
\[
\mathbb{G}_{S}^{\ast }(\mathbb{A}):=\{x\in \mathbb{A}:S(x)\in \mathbb{A}%
^{-1}\},\text{ and }\mathbb{G}_{S}(\mathbb{A}):=\{x\in \mathbb{A}:S(x)\in 
\mathbb{A}_{1}\}
\]%
with the operation $\circ _{S},$ 
\begin{equation}
x\circ _{S}y:=x+S(x)y,  \tag{$\circ _{S}$}
\end{equation}%
generating the \textit{Popa group} corresponding to $S$. The case $S(x):=1-x$
yields the \textit{circle operation} of the well established group of
`quasi-regular' elements of $\mathbb{A}$ [Ric, Ch. 1 \S 4]; cf. [Ost2].
Since $\mathbb{R}_{1}=(0,\infty ),$ $\mathbb{G}_{\rho }^{\ast }(\mathbb{R})=%
\mathbb{G}_{S}^{\ast }(\mathbb{R})$ and $\mathbb{G}_{\rho }(\mathbb{R})=%
\mathbb{G}_{S}(\mathbb{R})$ correspond to $S=S_{\rho }$ as in $(Can).$ For
further examples see \S 7 (Appendix).

We show in \S 4 that $\mathbb{G}_{S}(\mathbb{A})$ usually contains, as
abelian subgroups, copies of $\mathbb{G}_{\rho }(\mathbb{A})$ for $\rho \in 
\mathbb{A}^{-1}$; we also derive in Theorem 5.1 a characterization for $S$
by reference to an auxiliary function that solves an equation of Goldie
type, $(GGE_{SS}),$ generalizing $(GFE),$ by assuming the differentiability
of its multiplicative auxiliary $S$. This complements the Wo\l od\'{z}ko
approach [Wol] as encapsulated in a theorem in [Jav] and cited later in
[BriD]; that earlier approach readily extends to the present context
provided the invertible elements of the commutative field there are
interpreted as referring to $\mathbb{A}^{-1}$; see \S 6.1 with further
details in \S 7 (Appendix).

Other re-interpretations of $(GS)$ are possible. Of particular interest
below are the multivariate Popa groups of [BinO3,4] over a topological
vector space $X$ (note the larger category involved here), briefly the 
\textit{multi-Popa groups} (to maintain a clear distinction), with an
operation defined via a continuous linear functional $\rho \in X^{\ast },$ by%
\[
x\circ y:=x+(1+\rho (x))y 
\]%
for $x,y$ ranging over the half-space $\mathbb{G}_{\rho }(X):=\{x\in
X:1+\rho (x)>0\}.$ So here $S(x)=1+\rho (x)$ generalizes the canonical form $%
(Can)$ but with $S:X\rightarrow \mathbb{R}$, i.e. mapping to $\mathbb{R}$
rather than back to its domain; nevertheless, here $x$ acts affinely on $y,$
thus allowing location and uniform scaling. The groups $\mathbb{G}_{\rho
}(X) $ helpfully contribute a structure theorem describing the more general
Popa groups $\mathbb{G}_{S}(\mathbb{R}^{d})$ and $\mathbb{G}_{S}(C[0,1]).$

We do not pursue yet another interpretation, studied in [BriD], in which $%
S:X\rightarrow GL(X),$ for $X\ $a Banach space, so that $x$ acts on $y$ via $%
S(x)y$, but we do note below occasional similarities. Such similarities are
inevitable since $S(a)\in GL(\mathbb{A})$ for $a\in \mathbb{G}_{S}^{\ast }(%
\mathbb{A})$ (the map $x\mapsto S(a)x$ being a continuous automorphism of $%
\mathbb{A}$ since $||S(a)x||\leq ||S(a)||.||x||).$ It is all the more
unsurprising given that the componentwise product $x\cdot y$ of two vectors
in $\mathbb{R}^{d}$ may be presented as a matrix product $P(x)y$ in which $%
P(x):=$ \textrm{diag}$(x)$ is the diagonal matrix generated by $x$ (mapping
the product to composition: $P(x\cdot y)=P(x)P(y)$). Corresponding results
from [BriD, Ths 7, 8] will be seen from the present context as `degenerate'
variants of our results (while sometimes our results are special cases of
theirs): see the Remark after Th.4.2.

\bigskip

\textit{Banach algebras. }We are concerned below with the characterization,
in a unital commutative real Banach algebra $\mathbb{A}$, of those solutions 
$S:\mathbb{A\rightarrow A}$ of the Go\l \k{a}b-Schinzel equation%
\begin{equation}
S(x+yS(x))=S(x)S(y)\qquad (x,y\in \mathbb{A})  \tag{$GS_{\QTR{Bbb}{A}}$}
\end{equation}%
that, when restricted to $\mathbb{G}_{S}^{\ast }=\mathbb{G}_{S}^{\ast }(%
\mathbb{A}):=\{x:S(x)\in \mathbb{A}^{-1}\},$ are \textit{Fr\'{e}chet
differentiable} at the points of $\mathbb{G}_{S}(\mathbb{A})$ relative to
its range (in contrast to $\mathbb{A}$\textit{-differentiability}: see \S %
5). In this case, with $1_{\mathbb{A}}$ denoting the \textit{identity element%
} of $\mathbb{A}$ (under multiplication), a significant role is played by
the \textit{adjustor}, defined in Th. 5.1 as the map%
\[
N(x):=S(x)-1_{\mathbb{A}}-(S(1_{\mathbb{A}})-1_{\mathbb{A}})x\qquad (x\in 
\mathbb{G}_{S}^{\ast }). 
\]%
We assume here and below that 
\begin{equation}
1_{\mathbb{A}}\in \mathbb{G}_{S}^{\ast },  \tag{$1_{\QTR{Bbb}{A}}$}
\end{equation}%
i.e. that $S(1_{\mathbb{A}})\in \mathbb{A}^{-1}.$ Thus $N\ $measures
divergence from the canonical affine form. This is the central theme of \S 5.

Here and below $\mathbb{A}$ is always a unital commutative \textit{real}
Banach algebra (and so below `linear' means `$\mathbb{R}$-linear', unless
otherwise indicated), and the quantifier over $x,y$ in $(GS_{\mathbb{A}})$
is restricted to $\mathbb{G}_{S}^{\ast }(\mathbb{A}).$ It is noteworthy
that, unlike in the case of $\mathbb{A=R}$ (cf. [BinO1]), in a
multi-dimensional context \textit{quantifier weakening} (which we do not
pursue here) will broaden the nature of a solution function $S,$ as was
pointed out by Marshall and Olkin in [MarO1] for the similar context of the
multivariate Cauchy functional equation; see also [MarO2, MarO3] and \S 6.2.

Henceforth we view $S$ as a homomorphism. Our starting point is to establish
the group structure it generates on its domain, the significance of its
kernel $\mathcal{N}$ (`N for null') for its image, and `invariance of
openness' under $\circ _{S}$-shifts. Later in \S 5 we will be concerned with
the adjustor function $N$ above, which takes values in $\mathcal{N}$.

\bigskip

\noindent \textbf{Proposition 1.1. }\textit{Suppose }$S:\mathbb{A\rightarrow
A}$ \textit{satisfies }$(GS_{\mathbb{A}})$\textit{.}

(i)\textit{\ }$(\mathbb{G}_{S}^{\ast },\circ _{S})$ \textit{is a group and }$%
\mathbb{G}_{S}$\textit{\ a subgroup of} $\mathbb{G}_{S}^{\ast }$\textit{. }

(ii) \textit{Furthermore,}%
\[
\mathcal{N}=\mathcal{N}_{S}:=\{a\in \mathbb{A}:S(a)=1_{\mathbb{A}}\}\mathbb{%
\subseteq G}_{S}^{\ast } 
\]%
\textit{is a subgroup of }$\mathbb{G}_{S}$\textit{\ on which }$+$ \textit{%
and }$\circ _{S}$\textit{\ agree, so an additive subspace of }$\mathbb{A}$%
\textit{.}

(iii) \textit{If }$0$ \textit{is an interior point of }$\mathbb{G}_{S}(%
\mathbb{A})$\textit{, then so is }$w$ \textit{for any }$w\in \mathbb{G}_{S}(%
\mathbb{A}),$\textit{\ and, conversely, if }$w$\textit{\ is an interior
point of }$\mathbb{G}_{S}(\mathbb{A}),$ \textit{then so also is }$0.$ 
\textit{The same holds relativized to a line }$\langle u\rangle .$

(iv) \textit{For }$S$\textit{\ continuous, both }$(\mathbb{G}_{S},\circ
_{S}) $ \textit{and }$(\mathbb{G}_{S}^{\ast },\circ _{S})$ \textit{are
topological groups (in the subspace topology induced by }$\mathbb{A}$) 
\textit{and }$\mathcal{N}$\textit{\ is closed.}

\bigskip

\noindent \textbf{Proof. }For parts (i) and (ii), which are routine, see \S %
7 (Appendix).

(iii) We know from (i) that for $w\in \mathbb{G}_{S}(\mathbb{A})$ the map $%
x\mapsto w\circ _{S}x$ takes $\mathbb{G}_{S}(\mathbb{A})$ into $\mathbb{G}%
_{S}(\mathbb{A}).$ So if $B_{\delta }(0)\subseteq \mathbb{G}_{S}(\mathbb{A}%
), $ then%
\[
w+S(w)B_{\delta }(0)=w\circ _{S}B_{\delta }(0)\subseteq \mathbb{G}_{S}(%
\mathbb{A}). 
\]%
Here $S(w)B_{\delta }(0)$ is an open neighbourhood of $0$, since
multiplication by an invertible element of $\mathbb{A}\ $is a homeomorphism
of $\mathbb{A}$. So $w$ is interior to $\mathbb{G}_{S}(\mathbb{A}).$

It now also follows that, for $w\in \langle u\rangle $, if $x\in \langle
u\rangle \cap B_{\delta }(0)\subseteq \langle u\rangle \cap \mathbb{G}_{S}(%
\mathbb{A}),$ then $w+(S(w)B_{\delta }(0)\cap \langle u\rangle )\subseteq 
\mathbb{G}_{S}(\mathbb{A})\cap \langle u\rangle $ and, as $S(w)B_{\delta
}(0) $ is open, its intersection with $\langle u\rangle $ is relatively open.

Conversely, for $w\in \mathbb{G}_{S}(\mathbb{A}),$ by (i) the inverse $%
w_{S}^{-1}\in \mathbb{G}_{S}(\mathbb{A}),$ so if $(w+B_{\delta
}(0))\subseteq \mathbb{G}_{S}(\mathbb{A}),$ then again by (i),%
\begin{eqnarray*}
w_{S}^{-1}+S(w_{S}^{-1})(w+B_{\delta }(0)) &=&w_{S}^{-1}\circ (w+B_{\delta
}(0))\subseteq \mathbb{G}_{S}(\mathbb{A}): \\
S(w)^{-1}B_{\delta }(0) &=&w_{S}^{-1}+S(w)^{-1}(w+B_{\delta }(0))\subseteq 
\mathbb{G}_{S}(\mathbb{A}),
\end{eqnarray*}%
as $S(w)^{-1}=S(w_{S}^{-1}),$ $S$ being a homomorphism. Here $%
S(w)^{-1}B_{\delta }(0)$ is an open neighbourhood of $0$, again as
multiplication by an invertible element of $\mathbb{A}\ $is a homeomorphism
of $\mathbb{A}$.

A similar argument also holds under relativization to $\langle u\rangle .$

(iv) Now take $S$ continuous. Here $a\mapsto -aS(a)^{-1}$ is continuous,
since inversion is continuous on $\mathbb{A}^{-1}$ [Rud, p. 268, Th. 10.34],
[Con, VII Th.2.2]; clearly $(a,b)\mapsto a+S(a)b$ is continuous: $\mathbb{G}%
_{S}$ is a topological group.\hfill $\square $

\bigskip

\noindent \textbf{Corollary 1.1. }\textit{The multiplicative group} $(S(%
\mathbb{G}_{S}(\mathbb{A})),.)$\textit{\ is isomorphic to }$\mathbb{G}%
_{S}^{\ast }/\mathcal{N}$.

\bigskip

\noindent \textbf{Proof. }This is immediate, since $S(a\circ
_{S}b)=S(a)S(b); $ so by Prop. 1.1 $S$ is a homomorphism from $\mathbb{G}%
_{S}^{\ast }$ to $(S(\mathbb{G}_{S}(\mathbb{A})),.)$ with kernel $\mathcal{N}
$.\hfill $\square $

\bigskip

\noindent \textbf{Remarks. }1. As $\mathbb{A}^{-1}$ and $\mathbb{A}_{1}$ are
open [Con, Ch. 7 Th. 2.2], so too is $\mathbb{G}_{S}^{\ast }(\mathbb{A}%
)=S^{-1}(\mathbb{A}^{-1})$ and $\mathbb{G}_{S}(\mathbb{A})=S^{-1}(\mathbb{A}%
_{1}),$ for $S$ continuous. Following [DalF] say that $\mathbb{A}$ has 
\textit{dense invertibles} if $\mathbb{A}^{-1}$ is dense in $\mathbb{A},$ a
convenient property whenever invertibility is needed. In the present
circumstances this condition holds iff the \textit{topological stable rank}
of $\mathbb{A}$ is 1; this in turn is equivalent to the existence of a dense
set of points whose spectra have empty interior: see [CorS, Cor. 1.10].
Spectra emerge in \S 5.

\noindent 2.\textbf{\ }$S[\mathbb{G}_{S}^{\ast }(\mathbb{A})]$ is an
(abelian) multiplicative subgroup of $\mathbb{A}^{-1},$ as $S\ $is a
homomorphism (by $(GS)).$

\noindent 3. $\mathbb{A}_{1},$ the connected component of the identity,
coincides with the subgroup generated by the set of elements which have a
logarithm [Ric, Th.1.4.10], [Rud,Th. 10.34(c)], with connection from $1_{%
\mathbb{A}}$ to $g=e^{h}$ provided by $e^{th}.$

\noindent 4. If the operation $\circ _{S}$ is \textit{commutative}, then for 
$a,b\in \mathbb{A}^{-1}\cap \mathbb{G}_{S}^{\ast }(\mathbb{A})$ 
\[
a+S(a)b=b+S(b)a:\quad (1_{\mathbb{A}}-S(b))b^{-1}=(1_{\mathbb{A}%
}-S(a))a^{-1}=\text{constant }=-\rho ,\text{ say} 
\]%
and so 
\[
S(a)=1_{\mathbb{A}}+\rho a\quad (a\in \mathbb{A}^{-1}\cap \mathbb{G}_{S}(%
\mathbb{A})), 
\]%
hence for all $a\in \mathbb{G}_{S}$ if $\mathbb{A}$ has dense invertibles%
\textit{. }For $\mathbb{R},$ this is implicit in [GolS, Lemma 5] and
explicit for Popa [Pop, Prop. 3], where it is key. If $1_{\mathbb{A}}\in 
\mathbb{G}_{S}^{\ast }(\mathbb{A}),$ this Remark is non-vacuous. Evidently,
the operation $x\circ _{\rho }y$ applied to $x,y$ in any commutative ring is 
\textit{commutative}.

\bigskip

The Popa groups $\mathbb{G}_{\rho }(\mathbb{A)}$ below emerge in Prop 4.4 as
subgroups of $\mathbb{G}_{S}(\mathbb{A}).$

\bigskip

\noindent \textbf{Proposition 1.2.} \textit{The multiplicative group} $(S(%
\mathbb{G}_{S}(\mathbb{A})),.)$\textit{\ is isomorphic to the Popa group} $%
\mathbb{G}_{\rho }(\mathbb{A)}:=(\rho ^{-1}(\mathrm{ran}S-1_{\mathbb{A}%
}),\circ _{\rho })$, \textit{for each }$\rho \in \mathbb{A}^{-1}.$ \textit{%
The latter contains }$\{\rho ^{-1}(e^{t\rho }-1_{\mathbb{A}}):$\textit{\ }$%
t\in \mathbb{R}\}$\textit{\ as a one-parameter connected subgroup of} $%
\mathbb{G}_{\rho }(\mathbb{A)}$.

\bigskip

\noindent \textbf{Proof.} Put $c=\rho ^{-1}$; then $y=\eta _{\rho }(x):=1_{%
\mathbb{A}}+\rho x$ iff $\eta _{\rho }^{-1}(y)=c(y-1_{\mathbb{A}}).$ So%
\[
\eta _{\rho }^{-1}(\mathrm{ran}S)=c(\mathrm{ran}S-1_{\mathbb{A}});\quad \eta
_{\rho }(g)\eta _{\rho }(h)=\eta _{\rho }(g\circ _{\rho }h):\qquad g\circ
_{\rho }h=\eta _{\rho }^{-1}(\eta _{\rho }(g)\eta _{\rho }(h)). 
\]%
The final assertion is clear, since each $e^{t\rho }\in \mathbb{A}^{-1}.$%
\hfill $\square $

\section{Beyond Brillou\"{e}t-Dhombres-Brzd\k{e}k and linear-plus-1}

The theorems of this section are motivated by the observation that $(GS)$
may be solved with $S$ a Fr\'{e}chet differentiable function in the form%
\[
S(x)=1_{\mathbb{A}}+\gamma _{S}(x), 
\]%
for $\gamma =\gamma _{S}$ linear and continuous, provided $\gamma $ has the
following property which we may term $\mathbb{A}$\textit{-homogeneity }over $%
\mathbb{G}_{S}$: 
\[
\gamma (u\gamma (v))=\gamma (u)\gamma (v)\qquad (u,v\in \mathbb{G}_{S}). 
\]%
A weakened version of the property, considered at the end of the section, is
also relevant in describing solutions to the \textit{tilting equation} of
Section 5. To illustrate the property consider the following examples:

\noindent (i) $\gamma (x)=\gamma \cdot x$ for some $\gamma \in \mathbb{A}$,

\noindent (ii) $\gamma (x)=\gamma (x)\cdot 1_{\mathbb{A}}$ for some
continuous linear $\gamma :\mathbb{A}\rightarrow \mathbb{R}$: which includes
the case:

\noindent (iii) $\gamma (x)=x(\theta )$ for $x\in \mathbb{A}=C[0,1]$ and
some fixed $\theta $ with $0\leq \theta \leq 1.$

We will see later in Cor. 5.4\ that in $\mathbb{R}^{d}$ the property (ii)
may holds `partwise'. (That is, there is a (fixed) partition of $\{1,...,d\}$
such that (ii) holds for pairs of vectors restricted to the subspace
generated by the natural base vectors corresponding to any one part.) In
such cases $\gamma $ is $\mathbb{R}^{d}$-homogeneous.

\bigskip

Our main result in this section gives two decompositions of a Fr\'{e}chet
differentiable solution of $(GS)$ into a linear part and a part that is
orthogonal relative to both of the symmetric bilinear forms generated by $%
\gamma :$%
\[
\langle a,b\rangle =\gamma (ab)\text{ and }\langle a,b\rangle _{\gamma
}=\gamma (a\gamma (b)), 
\]%
though we have yet to prove (see below) the symmetry of the latter form.
When working relative to the latter, we speak of $\gamma $-orthogonality.
One decomposition yields a \textit{linear} and a non-linear part, the other
an $\mathbb{A}$\textit{-differentiable} part (in the sense of \S 5) --
ultimate source here of the distinction between analytic and real-analytic
solutions. We need a preliminary calculation.

\bigskip

\noindent \textbf{Proposition 2.1}. \textit{For }$S\ $\textit{Fr\'{e}chet
differentiable, satisfying }$(GS)$ 
\[
S^{\prime }(c)=S(c)S^{\prime }(0)S(c)^{-1}\qquad (c\in \mathbb{G}_{S}^{\ast
}). 
\]%
\noindent \textbf{Proof. }For fixed $a\in \mathbb{G}_{S}^{\ast },$ the
`affine' map $b\mapsto a\circ _{S}b=a+S(a)b$ is Fr\'{e}chet differentiable
on $\mathbb{A}$ and is onto (as $S(a)$ is invertible and multiplication is
continuous) with derivative $S(a)$. For $b\in \mathbb{G}_{S}^{\ast },$ as $%
a+S(a)b=a\circ _{s}b\in \mathbb{G}_{S}^{\ast }$ and $\mathbb{G}_{S}^{\ast }$
is open, with $S\ $Fr\'{e}chet differentiable, the Chain Rule applies [Ber,
Th. 2.1.15]: differentiating $(GS)$ with respect to $b$ and setting $%
a=b_{S}^{-1}$ yields 
\begin{eqnarray*}
S^{\prime }(a+bS(a))S(a) &=&S(a)S^{\prime }(b): \\
S^{\prime }(b) &=&S(b)S^{\prime }(0)S(b)^{-1},
\end{eqnarray*}%
since $S(a)^{-1}=S(a_{S}^{-1})=S(b)$. \hfill $\square $

The next result for $S$ real-valued and $\gamma $ injective can apply to $%
\mathbb{R}^{d}$ only for $d=1$, as then $\gamma $ has rank $1$ and nullity $%
0 $. This points up the essential difference between the linear form in
Theorem 2.1 and that of the Brillou\"{e}t-Dhombres-Brzd\k{e}k
characterization appropriate to Euclidean spaces (of $1+\rho (x)$ with $\rho 
$ in (the dual of) $\mathbb{R}^{d}$).

\bigskip

\noindent \textbf{Theorem 2.1 (Decomposition Theorem).}\textit{\ For }$S$%
\textit{\ Fr\'{e}chet differentiable satisfying }$(GS)$\textit{\ with }$%
\gamma =S^{\prime }(0)$, \textit{there exist two functions }$m,n$\textit{\
both with range} $\mathcal{N}_{\gamma }=\{x:\gamma (x)=0\}$ \textit{with}%
\newline
(i)\textit{\ }$n(x)$\textit{\ orthogonal in either of the above senses to }$%
\gamma (x)$ \textit{and with }$n(x)=o(x)$ \textit{as }$x\rightarrow 0,$%
\newline
(ii) $m(x)$\textit{\ orthogonal to }$\gamma (x)$ \textit{and }$\gamma $%
\textit{-orthogonal to both }$\gamma (x)$\textit{\ and }$\gamma (1_{\mathbb{A%
}})x,$\textit{\ and}%
\[
S(x)=1_{\mathbb{A}}+\gamma (x)+n(x)=1_{\mathbb{A}}+\gamma (1_{\mathbb{A}%
})x+m(x). 
\]%
\textit{In particular,}%
\[
S(\mathcal{N})\subseteq \{y:\gamma (y-1_{\mathbb{A}})=0\}=1_{\mathbb{A}}+%
\mathcal{N}_{\gamma }, 
\]%
\textit{and if }$\gamma $\textit{\ is injective}%
\[
S(x)=1_{\mathbb{A}}+\gamma (1_{\mathbb{A}})x. 
\]%
This will be a corollary of the following result in which we prove an $%
\mathbb{A}$-homogeneity property weakened by placing an additional $\gamma ,$
`like a mask', over the desired relation, as in $(\mathrm{pH})$ below.

\bigskip

\noindent \textbf{Remark.} For the special case of $S\ $real-valued on $%
\mathbb{A=C}$, the map $\gamma $ is homogeneous, being real-valued, so $%
\gamma (1)\zeta $ is complex, implying compensation by a necessarily complex
component $m(\zeta )$. So for this case, the decomposition is uninformative.
See Cor. 5.5.

\bigskip

\noindent \textbf{Theorem 2.2 (Phantom characterization of linearity-plus-1).%
} \textit{For }$S$\textit{\ Fr\'{e}chet differentiable with }$\gamma
=S^{\prime }(0)$\textit{: }$S$\textit{\ satisfies}%
\begin{equation}
S(a+bS(a))S(a)=S(a)S(b)\qquad (a,b\in \mathbb{G}_{S}\mathbb{)}  \tag{$GS$}
\end{equation}%
\textit{iff both }(i)%
\begin{equation}
\gamma (S(c)h)=\gamma ((1_{\mathbb{A}}+\gamma (c))h)\qquad (h\in \mathbb{A}%
,c\in \mathbb{G}_{S}\mathbb{)},  \tag{$\ast $}
\end{equation}%
\textit{and }(ii)%
\begin{equation}
\gamma (\gamma (k)h)=\gamma (k\gamma (h))\qquad (k,h\in \mathbb{A}), 
\tag{$\ast \ast $}
\end{equation}%
\textit{subject to the similarity relations }$S^{\prime }(c)=S(c)\gamma
S(c)^{-1},$ \textit{in which case for all} $a\in \mathrm{ran}(S),$%
\begin{equation}
\gamma ((a^{-1}\gamma a)(k)h)=\gamma (k\gamma (h))\qquad (k,h\in \mathbb{A}).
\tag{$\ast \ast \ast $}
\end{equation}%
\textit{In particular, }($\ast $) \textit{implies the phantom (or masked)
linearity-plus-1}%
\begin{equation}
\gamma (S(c))=\gamma (1_{\mathbb{A}}+\gamma (c))\qquad (c\in \mathbb{G}_{S}%
\mathbb{)}  \tag{\QTR{rm}{pL}}
\end{equation}%
\textit{\ (`p for phantom, L for linear'), while }($\ast \ast $) \textit{%
implies symmetry }%
\[
\langle a,b\rangle _{\gamma }=\langle b,a\rangle _{\gamma }, 
\]%
\textit{and }($\ast \ast \ast $) \textit{implies}%
\[
\gamma (\gamma (h)-\gamma (1_{\mathbb{A}})h)=0\qquad (h\in \mathbb{A}). 
\]%
\textit{This last implies the following phantom homogeneity:}%
\begin{equation}
\gamma (\gamma (a\gamma (b)))=\gamma (\gamma (a)\gamma (b))=\gamma (\gamma
(b\gamma (a))).  \tag{$\QTR{rm}{pH}$}
\end{equation}

\noindent \textbf{Proof.} Relegating the routine details to \S 7 (Appendix),
increment both $a$ and $b$ in $(GS)\ $by $h.$ Expansion of $S(a+h)$ and $%
S(b+h)$ to order $o(h)$ and use of Prop 2.1 give%
\begin{equation}
\gamma (S(c)h)=\gamma (h)+\gamma (c\gamma h),  \tag{$\ast ^{\prime }$}
\end{equation}%
with $c$ for the inverse of $b$ under $\circ ,$ all steps being reversible.
Differentiating with respect to $c$ in direction $k$ leads to%
\[
\gamma (k\gamma h)=\gamma (S^{\prime }(c)(k)h)=\gamma (S(c)\gamma
S(c)^{-1}(k)h), 
\]%
which for $a=S(c)$ yields the claim $(\ast \ast \ast )$ and for $c=0$ the
claim $(\ast \ast ).$ Furthermore, writing $S(c)k$ for $k$ gives via $(\ast
^{\prime })$ and $(\ast \ast )\ $that%
\[
\gamma (S(c)h)=\gamma (h)+\gamma (c\gamma h)=\gamma (h)+\gamma (h\gamma
(c))=\gamma (h(1+\gamma (c)). 
\]%
That is, $(\ast )\ $holds. Evidently $(\ast )$ and $(\ast \ast )$ yield $%
(\ast ^{\prime }),$ and so the conjuction of $(\ast )$ and $(\ast \ast )$
yields $(GS).$ The remaining conclusions are routine. \hfill $\square $

\bigskip

\noindent \textbf{Proof of Theorem 2.1. }Since $\gamma (\gamma (1_{\mathbb{A}%
})x)=\gamma (1_{\mathbb{A}}\gamma (x)),$%
\[
\gamma (S(x)-1_{\mathbb{A}}-\gamma (1_{\mathbb{A}})x)=\gamma (S(x)-1_{%
\mathbb{A}}-\gamma (x))=0. 
\]%
Now take $m(x):=S(x)-1_{\mathbb{A}}-\gamma (1_{\mathbb{A}})x\in \mathcal{N}%
_{\gamma }=\{x:\gamma (x)=0\}.$ By ($\ast \ast $) 
\[
\langle \gamma (x),m(x)\rangle =\gamma (\gamma (x)m(x))=\gamma (x\gamma
(m(x))=\gamma (0)=0, 
\]%
so $m(x)$ is orthogonal to $\gamma (x).$ Again by ($\ast \ast $) 
\begin{eqnarray*}
\langle \gamma (1_{\mathbb{A}})x,m(x)\rangle _{\gamma } &=&\gamma (\gamma
(1_{\mathbb{A}})x\gamma (m(x)))=0=\gamma (\gamma (\gamma (1_{\mathbb{A}%
})x)m(x)) \\
&=&\gamma (\gamma (\gamma (x))m(x))=\gamma (\gamma (x)\gamma (m(x)))=\langle
\gamma (x),m(x)\rangle _{\gamma }
\end{eqnarray*}%
so $m(x)$ is $\gamma $-orthogonal to both $\gamma (1_{\mathbb{A}})x$ and to $%
\gamma (x).$

Put $n(x)=S(x)-1_{\mathbb{A}}-\gamma (x);$ then $n(x)=o(x),$ and by ($\ast $)%
\begin{eqnarray*}
\langle \gamma (x),n(x)\rangle &=&\gamma (\gamma (x)(S(x)-1_{\mathbb{A}%
}-\gamma (x)))=0, \\
\langle \gamma (x),n(x)\rangle _{\gamma } &=&\gamma (\gamma (\gamma
(x))(S(x)-1_{\mathbb{A}}-\gamma (x)))=0.
\end{eqnarray*}%
If $S(x)=1_{\mathbb{A}},$ then clearly $\gamma (S(x)-1_{\mathbb{A}})=0,$
proving the final claim. \hfill $\square $

\bigskip

\noindent \textbf{Corollary 2.1}\textit{. With the assumptions of Theorem
2.1, for} $c\in \mathcal{N},$ \textit{both} $\gamma (c)\in \mathcal{N}%
_{\gamma }$ \textit{and} $\gamma (1_{\mathbb{A}})c\in \mathcal{N}_{\gamma }$%
. \textit{If }$\gamma $\textit{\ is injective, then}%
\[
S(c)=1_{\mathbb{A}}+\gamma (c)=1_{\mathbb{A}}+\gamma (1_{\mathbb{A}})c\qquad
(c\in \mathbb{G}_{S}\mathbb{)}, 
\]%
\textit{\ so }$\gamma $\textit{\ is }$\mathbb{A}$\textit{-homogeneous and }$%
\mathcal{N=N}_{\gamma }=\{0\}.$

\bigskip

\noindent \textbf{Proof. } Since $S(c)=1_{\mathbb{A}}$ for $c\in \mathcal{N}$%
, by Theorem 2.1%
\[
\gamma (1_{\mathbb{A}})=\gamma (S(c))=\gamma (1_{\mathbb{A}})+\gamma (\gamma
(c)), 
\]%
so $0=\gamma (\gamma (c))=\gamma (\gamma (1_{\mathbb{A}})c),$ yielding both $%
\gamma (c)\in \mathcal{N}_{\gamma }$ and $\gamma (1_{\mathbb{A}})c\in 
\mathcal{N}_{\gamma }.$

For $\gamma $ injective, $\mathcal{N}_{\gamma }=\{0\}$ and $S(c)=1_{\mathbb{A%
}}+\gamma (c),$ for $c\in \mathbb{G}_{S}.$ Likewise $\gamma \gamma
(c)=\gamma (\gamma (1_{\mathbb{A}})c)$ implies $\gamma (c)=\gamma (1_{%
\mathbb{A}})c,$ so $\gamma $ is $\mathbb{A}$-homogeneous on $\mathbb{G}_{S}$%
. Now if $c\in \mathcal{N},$ then $\gamma (c)=0,$ as $S(c)=1_{\mathbb{A}},$
so $\{0\}\mathcal{\subseteq N\subseteq N}_{\gamma }=\{0\},$ giving $\mathcal{%
N=N}_{\gamma }.$\hfill $\square $

\bigskip

\noindent \textbf{Remarks.} 1. When $S$ is real-valued and $\gamma \neq 0,$
the Cor. 2.1 captures the traditional and well established fact that $%
S(x)=1+\gamma (x)$ with $\gamma $ linear. The relevant statistical
literature includes Oakes and Dasu [OakD].\newline
2. Later Theorem 5.3 identifies $\mathcal{N}_{\gamma }$ as the maximal
vector subspace $\mathcal{H}$ of $\mathcal{N}$, whence above $\gamma (c)\in 
\mathcal{H}$ and $\gamma (1_{\mathbb{A}})c\in \mathcal{H}$ for $c\in 
\mathcal{N}$.

\bigskip

\noindent \textbf{Definition. }Recalling, from e.g.[Kec], the set-theoretic
notation\textbf{\ }$\omega =\{0,1,2,..\},$ say that $\gamma $ is $\omega $-%
\textit{homogeneous} if the homogeneity property $\gamma (v\gamma
(u))=\gamma (v)\gamma (u)$ holds for any $u$ and $v\in \{u\gamma (u)^{k}:$ $%
k\in \omega \}$. This is equivalent to a \textit{power-raising} (-shifting)
multiplicative effect of $u$ under $\gamma :$

\bigskip

\noindent \textbf{Proposition 2.2.} $\gamma $ \textit{is} $\omega $-\textit{%
homogeneous iff}%
\begin{equation}
\gamma (u\gamma (u)^{k})=\gamma (u)^{k+1}\text{ for all }u\text{ and all }%
k=0,1,...  \tag{$\times $}
\end{equation}%
\noindent \textbf{Proof}. A routine induction establishes this. See \S 7
(Appendix).\hfill $\square $

Below $f$ denotes both a real-analytic function and its natural extension to 
$\mathbb{A}.$ See \S 5 for applications (Prop. 5.1) and an extension of the
domain of validity (Theorem S).

\bigskip

\noindent \textbf{Corollary 2.2.} \textit{For }$f$\textit{\ with all Taylor
coefficients non-zero and radius }$R>0:$ \textit{the linear continuous }$%
\gamma $ \textit{satisfies }$(\times )$ \textit{for }$\gamma (u)\in \mathbb{A%
}^{-1}$ \textit{iff }%
\[
f(t\gamma (u))=\gamma (f(t\gamma (u))\cdot u/\gamma (u))\qquad (0\leq
t||\gamma (u)||<R). 
\]%
\textit{Thus this equivalence holds both for }$f(x):=e^{x}-1_{\mathbb{A}}$ 
\textit{and its inverse:} 
\begin{eqnarray*}
e^{t\gamma (u)}-1_{\mathbb{A}} &=&\gamma ((e^{t\gamma (u)}-1_{\mathbb{A}%
})\cdot u/\gamma (u))\qquad (t\geq 0), \\
\log (1_{\mathbb{A}}+t\gamma (v)) &=&\gamma (\log (1_{\mathbb{A}}+t\gamma
(v))\cdot v/\gamma (v))\qquad (0\leq t||\gamma (v)||<1).
\end{eqnarray*}

\noindent \textbf{Proof}. Assuming $(\times )$, apply $\gamma $ term by term
to the series expansion of $f$:%
\[
\gamma (uf(t\gamma (u))/\gamma (u))=\tsum\nolimits_{n=0}^{\infty
}a_{n}t^{n}\gamma (u\gamma (u)^{n})/\gamma (u)=f(t\gamma (u))\qquad (0\leq
t||\gamma (v)||<R). 
\]%
Conversely, compare coefficients at $t^{k}$ to obtain $(\times ).$\hfill $%
\square $

\bigskip

\noindent \textbf{Proposition 2.3.} \textit{Suppose that }$S$\textit{\
continuous satisfies }$(GS)$\textit{\ and, for some continuous linear }$%
\gamma ,$\textit{\ takes the form}%
\[
S(x)=1_{\mathbb{A}}+\gamma (x)+e(x). 
\]%
\textit{Then}%
\[
\lim_{x\rightarrow 0}e(x)/||x||^{2}=0\text{ implies the power-raising
property }(\times )\text{ above.} 
\]%
\textit{The converse holds for }$e(x\circ x)/2e(x)$\textit{\ bounded away
from }$1_{\mathbb{A}}$\textit{\ (as }$x\rightarrow 0$\textit{).}

\bigskip

\noindent \textbf{Proof}. For the sake of continuity, we defer this to \S 7
(Appendix).\hfill $\square $

\bigskip

\noindent \textbf{Remark.} The proof is somewhat reminiscent of the
Hyers-Ulam stability theorem with its near additivity: see e.g. [CabC]. We
conjecture that non-additivity of $e(.)$ implies boundedness away from
unity. Example 7.3 in the Appendix is illuminating here.

\section{Finite dimensions and $C[0,1]$}

Here we denote by $\cdot $ (rather than by $\odot $) the componentwise 
\textit{Hadamard-Schur} product applied to vectors in $\mathbb{R}^{d},$
which turns the $d$-vectors into a Banach algebra under the Euclidean norm.
For $S$ a continuous solution of $(GS)$ on this Banach space, $\mathbb{G}%
_{S}^{\ast }(\mathbb{R}^{d})$ is a topological group under the Euclidean
norm topology by Prop. 1.1, hence by the Montgomery-Zippin theorem this is a
Lie group [MonZ] or [Tao, Th. 1.1.13]. Then $S$ is $C^{\infty }$ in the
real-variable sense, cf. [BriD]. We use this fact in Corollary 5.6 below to
give a new treatment of $(GS_{\mathbb{C}})$ based on the differentiability
of the adjustor $N$ of \S 1, viewing $\mathbb{C}$ as a two-dimensional real
Banach algebra. Our first result below characterizes which `linear'
functions solve $(GS)$ in $\mathbb{G}_{S}^{\ast }(\mathbb{R}^{d})$. That
these are indeed the only \textit{non-degenerate} continuous solutions (i.e.
truly $d$-variate, below) is asserted thereafter in Th. 3.2. It is
convenient here to use the language of \textit{partitions} $\mathcal{P}$ of $%
\{1,...,d\}$ into (disjoint) subsets, termed \textit{parts} $I.$ When
needed, $|I|$ denotes the \textit{cardinality} of the part $I.$ Later, in
the context of $C[0,1],$ we use partitions of $[0,1]$ into compact parts $K.$

\bigskip

\noindent \textbf{Theorem 3.1 (Euclidean Characterization Theorem).}\newline
\noindent (i)\ \textit{The continuous solutions }$S:\mathbb{G}_{S}^{\ast }(%
\mathbb{R}^{d})\rightarrow \mathbb{R}^{d}$\textit{\ of }%
\begin{equation}
S(x+S(x)\cdot y)=S(x)\cdot S(y)  \tag{$GS$}
\end{equation}%
\textit{for }$S(x)=(...,1+\sigma _{i}(x),...)$ \textit{taking the form}%
\[
S(x):=\mathbf{1}+\Sigma x\text{ with }\Sigma =(\sigma _{ij}), 
\]%
\textit{where }$\mathbf{1}:=(1,1,...,1)^{\prime }$, \textit{have matrices }$%
\Sigma =(\sigma _{ij})$ \textit{satisfying, for} 
\[
\sigma _{i}(x):=\Sigma _{j}\sigma _{ij}x_{j}, 
\]%
\textit{\ }%
\[
\sigma _{ij}=0\text{ or }\sigma _{i}(x)\equiv \sigma _{j}(x)\qquad (1\leq
i,j\leq d). 
\]%
\noindent (ii) \textit{Hence there are: a linear map} $\sigma :\mathbb{R}%
^{d}\rightarrow \mathbb{R}^{d}$ \textit{with}%
\[
S(x):=\mathbf{1}+\sigma (x)\text{,} 
\]%
\textit{a `generator' functional} $\rho :\mathbb{R}^{d}\rightarrow \mathbb{R}
$ 
\[
\rho (x):=\Sigma _{i}\rho _{i}x_{i}, 
\]%
\textit{and a partition} $\mathcal{P}$ \textit{of }$\{1,...,d\}$\textit{\
with parts }$I$ \textit{so that, with }$e_{I}$\textit{\ the projection onto
the span} $\langle \{e_{i}:i\in I\}\rangle $ \textit{of the corresponding
natural base vectors }$e_{i}=(\delta _{ij}),$%
\begin{equation}
\sigma (x):=\sum\nolimits_{I\in \mathcal{P}}\sum\nolimits_{i\in I}\rho
(e_{I}x)e_{i}.  \tag{$\dag $}
\end{equation}%
\noindent (iii) \textit{For instance, partitioning of }$\{1,...,d\}$\textit{%
\ into two parts }$I,J$\textit{\ generates the solutions}%
\[
\sigma _{i}(x)=\sigma _{I}(x)=\Sigma _{k\in I}\rho _{k}x_{k}\text{ }(i\in
I),\qquad \text{ }\sigma _{j}(x)=\sigma _{J}(x):=\Sigma _{k\in J}\rho
_{k}x_{k}\text{ }(j\in J). 
\]%
\noindent (iv) \textit{In particular, in }$\mathbb{R}^{3}$ \textit{the
solutions }$S=(S_{1},S_{2},S_{3})^{T}$\textit{\ take one of the following
three forms:}%
\[
S_{i}(x)=1+\rho _{i}x_{i}\text{ for }i=1,2,3; 
\]%
\textit{or with }$(i,j,k)$\textit{\ a permutation of }$(1,2,3):$\textit{\ }%
\[
S_{i}(x)=S_{j}(x)=1+\rho _{i}x_{i}+\rho _{j}x_{j}\text{ and }S_{k}(x)=1+\rho
_{k}x_{k}; 
\]%
\textit{or}%
\[
S_{1}=S_{2}=S_{3}=1+\rho _{1}x_{1}+\rho _{2}x_{2}+\rho _{3}x_{3}. 
\]%
\textit{Thus the set }$\mathcal{N}_{S}=\{x\in \mathbb{R}^{3}:S(x)=\mathbf{1}%
\}$\textit{\ is a vector subspace of corresponding dimension }$0,1,2.$

\bigskip

\noindent \textbf{Proof. }(i) We compute the two sides of $(GS):$\textbf{\ } 
\begin{eqnarray*}
S(x)\cdot S(y) &=&(1+\Sigma x)\cdot (1+\Sigma y)=1+\Sigma x+\Sigma y+\Sigma
x\cdot \Sigma y, \\
S(x+S(x)y) &=&1+\Sigma (x+y+(\Sigma x)\cdot y)=1+\Sigma x+\Sigma y+\Sigma
((\Sigma x)\cdot y).
\end{eqnarray*}%
On comparing, $(GS)$ reduces to%
\[
\Sigma ((\Sigma x)\cdot y)=\Sigma x\cdot \Sigma y. 
\]%
We compute the $i$-th component on each side:%
\begin{eqnarray*}
RHS_{i} &=&(\Sigma _{j}\sigma _{ij}x_{j})(\Sigma _{k}\sigma
_{ik}y_{k})=\Sigma _{jk}\sigma _{ij}\sigma _{ik}x_{j}y_{k}, \\
LHS_{i} &=&\Sigma _{k}\sigma _{ik}\Sigma _{j}\sigma _{kj}x_{j}y_{k}=\Sigma
_{jk}\sigma _{ik}\sigma _{kj}x_{j}y_{k}.
\end{eqnarray*}%
Comparison of the coefficient of $x_{j}y_{k}$ on each side yields that, for
all $ijk,$%
\[
\sigma _{ij}\sigma _{ik}=\sigma _{ik}\sigma _{kj}:\quad \sigma _{ik}=0\text{
or }\sigma _{ij}=\sigma _{kj}, 
\]%
as asserted.

(ii) The partition statements are corollaries of (i) as follows.

For $\rho \in \mathbb{R}^{d}$ and $K\subseteq \{1,...,d\}$ put%
\[
\sigma _{K}(x):=\Sigma _{k\in K}\text{ }\rho _{k}x_{k}. 
\]%
Thus $(\sigma _{K})_{k}=0$ for $k\notin K.$ Partition $\{1,...,d\}$ into $%
I,J $ and take%
\[
\sigma _{i}=\sigma _{I}\text{ for }i\in I,\qquad \sigma _{j}=\sigma _{J}%
\text{ for }j\in J. 
\]%
Then $\sigma _{ij}=0$ for $i\in I,j\in J$; indeed, $\sigma _{ij}=0$ for $%
i\in I,$ as $j\notin I.$ Similarly for $i\in J\ $and $j\in I.$ So%
\begin{eqnarray*}
\sigma _{i} &=&\sigma _{k}\text{ for }i,k\in I\text{ and }\sigma _{j}=\sigma
_{k}\text{ for }j,k\in J, \\
\sigma _{ij} &=&0\text{ for }(i,j)\in (I\times J)\cup (J\times I).
\end{eqnarray*}%
Thus $\sigma _{ik}=0$ or $\sigma _{i}=\sigma _{k}.$

(iii) and (iv) are now immediate corollaries of (ii). \hfill $\square $

\bigskip

\noindent \textbf{Remark. }For $\mathbb{A}=\mathbb{R}^{2},$ $S$ either takes
the \textit{independent form}:%
\[
S(x)=S(x_{1},x_{2})=1+\rho x=(1+\rho _{1}x_{1},1+\rho _{2}x_{2}),\text{ with 
}\rho _{1}\rho _{2}\neq 0, 
\]%
and then $\mathcal{N}(\rho ):=\{(x_{1},x_{2}):(\rho _{1}x_{1},\rho
_{2}x_{2})=(0,0)\}=\{(0,0)\}$ has co-dimension 2; or it takes the \textit{%
co-dependent form} with its two components given by $S_{1}(x)=S_{2}(x)=1+%
\rho _{1}x_{1}+\rho _{2}x_{2},$ where the corresponding $\mathcal{N}(\rho )$
has co-dimension 1.

In $\mathbb{R}^{3}$ the co-dimensions can be $3,2,$ or $1:$%
\[
\mathcal{N}(1,1,1)=\{0,0,0\},\mathcal{N}(0,1,1)=\mathbb{R\times \{(}0,0)\},%
\mathcal{N}(0,0,1)=\mathbb{R}^{2}\mathbb{\times \{}0\}. 
\]%
Here the corresponding ranges $\mathcal{N}(\rho )$ are of dimensions $0,1,2.$

\bigskip

We complement Th. 3.1 by proving that it is \textit{exhaustive}, i.e.
includes all the continuous solutions, as qualified below. The proof is in
principle straightforward: it reduces the analysis of $(GS)$ to a series of
interconnected scalar equations. This is unavoidably lengthy and messy. For
the sake of simplicity we confine ourselves to the case $d=2$. For the sake
of continuity, we defer the proof to \S 7 (Appendix).

\bigskip

\noindent \textbf{Theorem 3.2 (Exhaustivity})\textbf{. }\textit{The
non-trivial continuous solutions for }$s:\mathbb{R}^{2}\rightarrow \mathbb{R}%
^{2}$\textit{\ to the Go\l \k{a}b-Schinzel equation below in the algebra of }%
$\mathbb{R}^{2}$%
\begin{equation}
s(a+s(a)b)=s(a)s(b)  \tag{$GS$}
\end{equation}%
\textit{take for some }$\rho \in \mathbb{R}^{2}$ \textit{either the
`co-dependent form': }%
\begin{eqnarray*}
s(x_{1},x_{2}) &=&(1+\rho _{1}x_{1}+\rho _{2}x_{2},1+\rho _{1}x_{1}+\rho
_{2}x_{2})=(\eta _{\rho }(x),\eta _{\rho }(x)) \\
&=&\mathbf{1}+x(\rho ^{T},\rho ^{T})\qquad \text{(matrix form),}
\end{eqnarray*}%
\textit{or the `independent form':}%
\[
s(x_{1},x_{2})=(1+\rho _{1}x_{1},1+\rho _{2}x_{2}), 
\]%
\textit{or the `degenerate' (univariate) form : }%
\[
s(x_{1},x_{2})=(\eta _{\rho }(x_{1}),\tau (x_{1}))\text{ or }(\tau
(x_{2}),\eta _{\rho }(x_{2})), 
\]%
\textit{with }$\rho \in \mathbb{R}$, \textit{where }$\tau $\textit{\ solves
the homomorphism equation}%
\begin{equation}
\tau (x_{1}\circ _{\rho }y_{1})=\tau (x_{1})\tau (y_{1}).  \tag{$Hom$}
\end{equation}%
\textit{\ For instance, corresponding to }$\rho =0,\rho \in (0,\infty )$%
\textit{\ and} $\rho =\infty ,$ \textit{for some }$\gamma \in \mathbb{R}$%
\[
s(x_{1},x_{2})=(1,e^{\gamma x_{1}})\text{ or }(1+\rho x_{1},(1+\rho
x_{1})^{\gamma })\text{ or }(x_{1},x_{1}^{\gamma }),\text{ respectively.} 
\]

\bigskip

\noindent \textbf{Remark.} \textit{Degenerate solutions} of $(GS)$ in the
context $S:\mathbb{R}^{d}\rightarrow \mathbb{R}^{d}$ occur when the $d$
components of $x$ may be partitioned, as $x=(u,v)$ say, and $S(x)=S(u,0)$ so
that $S$ is \textit{not truly }$d$\textit{-variate}; then for the
corresponding partition $S(x)=(s(u),t(u)),$ the $(GS)$ equation reduces to
the two equations%
\[
s(a+bs(a))=s(a)s(b),\qquad t(a+bs(a))=t(a)t(b). 
\]%
A degenerate solution thus departs from the `1-plus-linear' form and couples
a lower-dimensional solution $s$ of $(GS)$ (of 1-plus-linear form') with a
sequence of components $s(u)^{\gamma _{i}}$ (\textquotedblleft $s$%
-homomorphism\textquotedblright ), or in the maximally degenerate case when $%
s\equiv 1$ with a sequence of exponential components $e^{\langle \gamma
_{i},u\rangle }$ as in [BriD, Th. 8].

\bigskip

We turn now to the Structure Theorem, Th. 3.3 below, and its proof. We use
`tilde on tilde' $\approx $ for `is isomorphic to' and $d_{I}:=|I|$ (the
cardinality of a part $I$ of the partition $\mathcal{P}$ of $\{1,...,d\}).$

\bigskip

\noindent \textbf{Theorem 3.3 (Structure Theorem).} (i) \textit{The
Banach-algebra Popa group }$\mathbb{G}_{S}^{\ast }(\mathbb{R}^{d})$ \textit{%
generated by a continuous solution }$S$\textit{\ (with invertible values) of
the Go\l \k{a}b-Schinzel equation }$(GS)$ \textit{on }$\mathbb{R}^{d}$%
\textit{\ equipped with Hadamard product is isomorphic to a direct product
of multi-Popa groups:}%
\[
\mathbb{G}_{S}^{\ast }(\mathbb{R}^{d})\approx \otimes _{I\in \mathcal{P}}%
\mathbb{G}_{\sigma _{I}}\mathbb{(R}^{d_{I}}) 
\]%
\textit{for a partition }$\mathcal{P}$ \textit{of }$\{1,...,d\}$ \textit{and
some linear maps} $\sigma _{I}:\mathbb{R}^{d_{I}}\rightarrow \mathbb{R}$%
\textit{\ with }$I\in \mathcal{P}$\textit{.}

(ii) \textit{Likewise, the Popa group }$\mathbb{G}_{S}^{\ast }(C[0,1])$ 
\textit{is isomorphic to a direct product of multi-Popa groups:}%
\[
\mathbb{G}_{S}^{\ast }(C[0,1])\approx \otimes _{K\in \mathcal{P}}\mathbb{G}%
_{S_{K}}(C[K]) 
\]%
\textit{for some} \textit{partition }$\mathcal{P}$ \textit{of }$[0,1]$ 
\textit{into compact subsets and some 1-plus-linear maps} $%
S_{K}:C(K)\rightarrow \mathbb{R}$ \textit{for} $K\in \mathcal{P}.$

\bigskip

\noindent \textbf{Proof. }(i)\textbf{\ }For $\mathbb{G}_{S}^{\ast }(\mathbb{R%
}^{d})$: by Theorem 3.1 for each $I$ the relevant affine map is%
\[
S_{I}:=S|\mathbb{R}^{I}:\mathbb{R}^{d_{I}}\rightarrow \mathbb{R}. 
\]%
Furthermore, projecting $x+y\cdot S_{I}(y)$ to $\mathbb{R}^{d_{I}}$ yields
for $x_{I}:=\langle x_{i}:i\in I\rangle $ and $y_{I}:=\langle y_{i}:i\in
I\rangle $ 
\begin{eqnarray*}
\langle x_{i} &:&i\in I\rangle +\langle y_{i}S_{I}(x_{I}):i\in I\rangle
=\langle x_{i}:i\in I\rangle +S_{I}(x_{I})\langle y_{i}:i\in I\rangle \\
&=&x_{I}+y_{I}S_{I}(x_{I}),
\end{eqnarray*}%
which is the same binary operation as that of the multivariate Popa group $%
\mathbb{G}_{S_{I}}\mathbb{(R}^{d_{I}}).$ The functions $\{S_{I}:I\in 
\mathcal{P}\}$ may then be `merged', as in Th. 3.1, to yield a single linear
generator function $\rho :\mathbb{R}^{d}\rightarrow \mathbb{R}$ satisfying $%
(\dag )$, say with%
\[
\rho (x):=\Sigma _{i=1}^{d}\sigma _{i}x_{i}. 
\]

(ii) This follows from (i) by Stone-Weierstrass approximation; see \S 7
(Appendix).\hfill $\square $

\bigskip

\textbf{Remark.} In Th. 3.3 for $x\in \mathbb{G}_{S}^{\ast }(C[0,1])$ domain
restriction $x\mapsto x_{K}$ represents the isomorphism map via%
\[
x\mapsto (x_{K})_{K\in \mathcal{P}},\quad S(x)\mapsto (S(x)_{K})_{K\in 
\mathcal{P}}. 
\]%
Thus, if $\{k\}=K\in \mathcal{P}$ with $k\in \lbrack 0,1],$ then $%
x_{K}=x(k), $ giving $C(K)=\mathbb{R}$ and $\mathbb{G}_{\sigma _{K}}(C[K])=%
\mathbb{G}_{\sigma _{k}}(\mathbb{R}),$ where $\sigma _{K}(t)=1+\sigma _{k}t$%
, say; so%
\begin{eqnarray*}
\sigma _{K}(x) &=&1+\sigma _{k}x(k): \\
S(x)(k) &\mapsto &S(x)_{K}=1+\sigma _{k}x(k).
\end{eqnarray*}%
In particular, for $\mathcal{P}=\{\{k\}:k\in \lbrack 0,1]\}$, take $\rho
(k):=\sigma _{k}$ for $k\in \lbrack 0,1];$ then%
\[
S(x)=1+\rho x 
\]%
with $\rho \in C[0,1]$ (since $\rho =S(1)-1\in C[0,1]$). Here $\rho \neq 0:$
for any $s\neq t,$ since $K_{t}=\{t\},$ there is $x_{s}$ with $\rho
(t)x_{s}(t)\neq \rho (s)x_{s}(s).$

\bigskip

\noindent \textbf{Remark. }Our results subsume and extend the similar
results in [BriD], concerned with $S(x)$ in diagonal form, specifically
their Th. 7(i) for $\mathbb{R}^{2}$; likewise their Th. 8 refers to those
functions $S(x_{1},...,x_{d})$ in diagonal form which, like the univariate
types above, are from our perspective degenerate solutions to $(GS)$ through
not being properly $d$-variate.

\section{Spanning Pencil Theorem}

Our main theorem here, Th. 4.1 below, shows that neighbourhoods of the
origin are spanned by a pencil of disjoint (modulo $0$) subgroups isomorphic
to canonical Popa groups. As a preliminary, we study general properties of
the set $\mathcal{N}=\mathcal{N}_{S}=S^{-1}(1_{\mathbb{A}})$ for $S$ a 
\textit{continuous} solution of $(GS_{\mathbb{A}})$ in the context of a
unital commutative real Banach algebra $\mathbb{A}$. We have seen in Prop
1.1 that $\mathcal{N}$ is additive. A significant issue, adressed both in
this and in the next section when discussing $(GFE)$ in the Banach algebra
setting, is whether $\mathcal{N}$ is a vector subspace; see e.g. Lemma 4.2
and Lemma 5.4. This is not an easy matter to verify for $\mathbb{A}$ of
dimension higher than 2 and seems to require additional hypotheses either on
the range of $S$ or on appropriate solubility of the tilting equation ($T$)
in \S 5 (see Prop. 5.3). Even for an image isomorphic to a (commutative)
field such as $\mathbb{R}$ or $\mathbb{C}$, quite some effort may be
required: see the verification in [BriD, Prop. 3] and within the proof of
[Brz, Th. 3]. (There the function is not assumed continuous: instead of
being open, the associated Popa group $\mathbb{G}_{S}^{\ast }(X)$ is assumed
to have the algebraic interior point property [Lyu, \S 2.2].) That said,
note that passing to $\mathbb{C}$-valued functions on $\mathbb{C}$, $(GFE)$
is satisfied by $(K,g)$ with multiplicative auxiliary $g(\zeta )=e^{\zeta }$
whose level set $g=1$ is discrete. (Taking $K(\zeta ):=(1-e^{\zeta })/2,$
the level set $K=1$ is also discrete, being the solution set of $e^{\zeta
}=-1.)$

When $S$ is Fr\'{e}chet differentiable at $0$ and $\mathcal{N}$ is a vector
subspace, Theorem 5.3 gives a \textit{differential characterization} of $%
\mathcal{N}$ as 
\[
\mathcal{N}=\{u:DS(0)u=0\}, 
\]%
with the adjustor $N$ (of \S 1) linear on $\mathcal{N}.$ As a result this
allows decomposition of any continuous solution of $(GEE)$ into the sum of a
linear function and one that is exponential in a curvilinear sense (as a
corollary of Th. 5.2).

Below, for $\Sigma \subseteq \mathbb{A},$ $\langle \Sigma \rangle $ denotes
the vector subspace of $\mathbb{A}$ generated by $\Sigma .$ Some of the
initial results here, like Lemma 4.1, are known for $\mathbb{R}_{+}$: see
e.g. [Mur, Lemma 1] (the proof needs only invertibility) which generalizes
to the invertible elements $\mathbb{A}^{-1}$.

\bigskip

\noindent \textbf{Lemma 4.1 (}[GolS, Lemma 1], [Mur, Lemma1])\textbf{.} 
\textit{For }$a,b\in \mathbb{G}_{S}^{\ast }$\textit{, if }$S(a)=S(b),$%
\textit{\ then:}

\noindent (i)\textit{\ }$S(c)=S(c+a-b)$\textit{\ for any }$c\in \mathbb{G}%
_{S}^{\ast }$, \textit{so in particular, }$S(a-b)=S(0)=1_{\mathbb{A}}$%
\textit{;}

\noindent (ii)\textit{\ }$S(c)=S(c+S(z)(a-b))$ \textit{for any} $z\in 
\mathbb{G}_{S}^{\ast }$, \textit{so in particular }$S(z)a\in \mathcal{N}$%
\textit{\ for }$a\in \mathcal{N}$ \textit{and so }%
\[
S(\mathbb{G}_{S}^{\ast })\mathcal{N}=\mathcal{N}. 
\]

\noindent \textbf{Proof: }See \S 7 (Appendix).\hfill $\square $

\bigskip

\noindent \textbf{Theorem B }([Bou, VII \S 2, Prop. 3]). \textit{Every
non-discrete closed additive subgroup of }$\mathbb{R}^{n}$\textit{\ contains
a one-dimensional vector subspace.}

\bigskip

The proof rests on local compactness. Notice the inherent limitation: the
subgroup $\mathbb{R\times Z}$ is closed and non-discrete and contains a
vector subspace but is not itself a vector space; we shall exclude this for $%
\mathcal{N}$ below in the case $\mathbb{A=R}^{2}.$

\bigskip

\noindent \textbf{Proposition 4.1.}\textit{\ For }$S$\textit{\ a continuous
solution of }$(GS),$ $\mathcal{N}$\textit{\ is closed and either the
singleton }$\{0\}$\textit{\ or dense-in-itself. For}%
\[
\mathcal{N}_{0}:=\mathcal{N\cap }\mathbb{G}_{S}, 
\]%
$\mathcal{N}_{0}$\textit{\ and so }$\mathcal{N}$ \textit{is closed and
either the singleton }$\{0\}$\textit{\ or dense-in-itself.}

\bigskip

\noindent \textbf{Proof.} For the details see \S 7 (Appendix). In brief: for 
$b$ close enough to $0$ and any $0\neq $ $a\in \mathcal{N}$, $S(b)a\in 
\mathcal{N}$ by Lemma 4.1(ii). Then $||a-aS(b)||\leq ||a||.||1_{\mathbb{A}%
}-S(b)||,$ so $a$ is an accumulation point of $\mathcal{N}$, unless $S=1_{%
\mathbb{A}}$ near $0$. But then $B_{\varepsilon }(0)\subseteq \mathcal{N}$
for some $\varepsilon >0$ and so $\mathcal{N}$, being an additive subgroup,
contains $\bigcup\nolimits_{n\in \mathbb{N}}B_{n\varepsilon }(0)=\mathbb{A}$%
, and so is dense-in-itself (and also $S\equiv 1$).\hfill $\square $

\bigskip

\noindent \textbf{Corollary 4.1. }\textit{For }$0\neq a\in \mathcal{N},$%
\textit{\ either} $\langle a\rangle \subseteq \mathcal{N}$ \textit{or there
is }$b\in \mathbb{G}_{S}^{\ast }$, \textit{with }$S(b)\notin \langle
1_{A}\rangle .$

\bigskip

\noindent \textbf{Proof. }W.l.o.g. $\mathcal{N\neq }\{0\}$. Following [Dal],
put $B_{\varepsilon }^{\bullet }(0):=B_{\varepsilon }(0)\backslash \{0\}$.
Note first that $(\langle a\rangle \backslash \{a\})\cap S(B_{\varepsilon
}^{\bullet }(0))a\neq \emptyset $, provided $S(b)a=ta$ or $S(b)=t1_{\mathbb{A%
}},$ for some $t\neq 1$.

By Prop. 4.1, for each $a\in \mathcal{N},$ one of two cases may arise:

\noindent (i) $(\langle a\rangle \backslash \{a\})\cap S(B_{\varepsilon
}(0))a\neq \emptyset $ for each $\varepsilon >0,$ i.e. $S(B_{\varepsilon
}(0))\cap \langle 1_{\mathbb{A}}\rangle $ contains points other than $1_{%
\mathbb{A}}$;

\noindent (ii) for some $\varepsilon >0,$ $S(B_{\varepsilon }^{\bullet
}(0))a\cap \langle a\rangle =\emptyset .$

If case (ii) never arises, then $\langle a\rangle \cap \mathcal{N}$ is
closed and dense and so $\langle a\rangle \subseteq \mathcal{N}$ for each $%
a\in \mathcal{N}$ (and then $\mathcal{N}$ is a vector subspace, being an
additive subgroup).

Otherwise, there is some $\varepsilon >0$ and $a_{0}\in \mathcal{N}$ such
that $S(b)a_{0}\notin \langle a_{0}\rangle $ for all $b\in B_{\varepsilon
}^{\bullet }(0).$ So for each $t\in \mathbb{R},$ $S(b)a_{0}\neq ta_{0},$
implying $S(b)\neq t1_{\mathbb{A}}$ (otherwise $S(b)=t1_{\mathbb{A}}$
implies $S(b)a=ta),$ and then $S(b)\notin \langle 1_{\mathbb{A}}\rangle .$
\hfill $\square $

\bigskip

\noindent \textbf{Lemma 4.2. }\textit{If} $S(\mathbb{G}_{S}^{\ast }(\mathbb{%
A)})\supseteq \mathbb{R}_{+}1_{\mathbb{A}},$ \textit{then} $\mathcal{N}$%
\textit{\ is a vector space.}

\bigskip

\noindent \textbf{Proof.} For $0\neq a\in \mathcal{N}$ and each $t>0$ there
is $b_{t}$ with $S(b_{t})=t$, so $ta=t1_{\mathbb{A}}a=S(b_{t})a\in S(\mathbb{%
G}_{S}(\mathbb{A)})\mathcal{N}\subseteq \mathcal{N}$. By Theorem B $\langle
a\rangle \cap \mathcal{N}$, being up to the isomorphism $t\mapsto ta$ a
dense subgroup of $\mathbb{R}$, is all of $\mathbb{R}$, so $\langle a\rangle
\subseteq S(\mathbb{G}_{S}^{\ast }(\mathbb{A)})\mathcal{N}\subseteq \mathcal{%
N}$ and so $\mathcal{N}$ is a vector subspace. \hfill $\square $

\bigskip

The condition in Lemma 4.2 is not fulfilled in Example 7.3 in the Appendix.
We return to this matter below.

\bigskip

A simple corollary is that for a two-dimensional $\mathbb{A}$ the subgroup $%
\mathcal{N}$ is a vector subspace. For this we need the following

\bigskip

\noindent \textbf{Lemma G }(cf. [Geb])\textbf{. }\textit{For }$a\in \mathbb{A%
}^{-1}$\textit{\ with} $\langle a\rangle \subseteq \mathcal{N}$ \textit{and }%
$b\in \mathbb{A}$, \textit{if }$S(b)\notin \langle 1_{A}\rangle ,$\textit{\
then }$\langle a,S(b)a\rangle $\textit{\ is a two-dimensional vector
subspace of }$\mathcal{N}$. \textit{Likewise,} \textit{for }$a\in \mathbb{A}$%
\textit{\ with} $\langle a\rangle \subseteq \mathcal{N}$ \textit{and }$b\in 
\mathbb{A}$, \textit{if }$S(b)a\notin \langle a\rangle ,$\textit{\ then }$%
\langle a,S(b)a\rangle $\textit{\ is a two-dimensional vector subspace of }$%
\mathcal{N}$.

\bigskip

\noindent \textbf{Proof. }Two-dimensionality of $\langle a,S(b)a\rangle $ is
clear, as otherwise $S(b)a=ta$ for some $t\in \mathbb{R}$, and then $%
S(b)=t1_{A},$ a contradiction. Evidently, $tS(b)a=S(b)ta\in \mathcal{N}$ as
by Lemma 4.1(ii) $S(b)\langle a\rangle \subseteq \mathcal{N}$, i.e $\langle
S(b)a\rangle \subseteq \mathcal{N}$, and so $\langle a,S(b)a\rangle
\subseteq \mathcal{N}$.\hfill $\square $

\bigskip

\noindent \textbf{Corollary 4.2. }\textit{If }$\mathbb{A}$ \textit{has
dimension }$\leq 2$\textit{, then }$\mathcal{N}$\textit{\ is a vector
subspace.}

\bigskip

\noindent \textbf{Proof. }W.l.o.g.\textbf{\ }$\mathcal{N\neq }\{0\}.$ By
Prop. 4.1 and Theorem B there is $a\in \mathcal{N}$ with $\langle a\rangle
\subseteq \mathcal{N}$. So either $\mathcal{N=}\langle a\rangle ,$ a vector
subspace, or otherwise by Lemma G $\mathcal{N}$ contains a two-dimensional
subspace, and so $\mathcal{N}=\mathbb{A}$, again a vector subspace.\hfill $%
\square $

\bigskip

Lemma G has a natural extension.

\bigskip

\noindent \textbf{Lemma 4.3. }\textit{If} $\mathcal{N}$ \textit{contains an }%
$m$\textit{-dimensional subspace generated by }$\{S(b_{1})a,...,S(b_{m})a\}$%
\textit{\ and }\textrm{ran}$S(B_{\varepsilon }(0))$\textit{\ is
topologically at least }$(m+1)$\textit{-dimensional, then} $\mathcal{N}$ 
\textit{contains an }$(m+1)$\textit{-dimensional subspace }$\langle
a,S(b_{1})a,...,S(b_{m+1})a\rangle $.

\bigskip

\noindent \textbf{Proof. }W.l.o.g. $\langle 1_{A}\rangle \cap $ \textrm{ran}$%
S(B_{\varepsilon }(0))=\emptyset .$ Take $\Sigma =\{S(b_{1}),...,S(b_{m})\};$
then $\langle \Sigma \rangle $ is $m$-dimensional, so there is $b$ with 
\[
S(b)a\notin \langle \Sigma a\rangle , 
\]%
and so $\langle S(b)a,\Sigma \rangle $ is $(m+1)$-dimensional. Then as $S(G)%
\mathbb{R}a\subseteq S(G)\mathcal{N}\subseteq \mathcal{N}$, for each such $%
b, $ $\langle S(b)a,\Sigma \rangle \subseteq \mathcal{N}$, i.e. $\mathcal{N}$
contains an $(m+1)$-dimensional subspace.\hfill $\square $

\bigskip

\noindent \textbf{Corollary 4.3. }\textit{If }\textrm{ran}$S$\textit{\ and }$%
\langle \mathcal{N}\rangle $ \textit{are both }$n$\textit{-dimensional, then 
}$\mathcal{N}$\textit{\ is a vector subspace.}

\bigskip

\noindent \textbf{Proof. }As before by Prop. 4.1 and Theorem B, $\mathcal{N}$
contains a $1$-dimensional subspace. Now apply the preceeding Lemma $n-1$
times.\hfill $\square $

\bigskip

We now study the condition $\mathrm{ran}S\supseteq \mathbb{R}_{+}1_{\mathbb{A%
}}.$ Our main tool is the functional equation for $g:\mathbb{R}%
_{+}\rightarrow \mathbb{A}$%
\begin{equation}
g(s)+sg(t)=g(st)\qquad (s,t\in \mathbb{R}_{+}),  \tag{$g_{\QTR{Bbb}{R}}$}
\end{equation}%
which, as we shall see below, has solution%
\[
g(t)=c_{g}(t-1)\qquad (t\in \mathbb{R}_{+}) 
\]%
for some $c_{g}\in \mathbb{A}$. We return to a more general variant of $(g_{%
\mathbb{R}})$ in Prop. 4.4 and again in \S 5 below where we see the more
general appearance of the Popa group $\mathbb{G}_{\rho }(\mathbb{A}),$ here
in the form $\rho ^{-1}(\mathbb{R}_{+}\mathbb{-}1_{\mathbb{A}})=(-1,\infty
)\rho ^{-1},$ isomorphic to $\mathbb{G}_{1}(\mathbb{R})$ and so to $(\mathbb{%
R}_{+},\mathbb{\times )}$.

\bigskip

\noindent \textbf{Proposition 4.2.} \textit{For }$S$\textit{\ a solution of }%
$(GS),$ \textit{if} $\mathrm{ran}S\supseteq \mathbb{R}_{+}1_{\mathbb{A}}$%
\textit{, then for some }$c\in \mathbb{A}$ 
\[
S(c(t-1))=t\qquad (t\in \mathbb{R}_{+}). 
\]%
\textit{So if }$c\in \mathbb{A}^{-1}$\textit{, then for }$\rho =c^{-1}$%
\[
S(w)=1_{\mathbb{A}}+\rho w\qquad (w\in \rho ^{-1}(\mathbb{R}_{+}\mathbb{-}1_{%
\mathbb{A}})=(-1,\infty )\rho ^{-1}); 
\]%
\textit{in particular }$\mathbb{G}_{\rho }(\mathbb{A})=\rho ^{-1}(\mathbb{R}%
_{+}\mathbb{-}1_{\mathbb{A}})$ \textit{is a subgroup of }$\mathbb{G}_{S}(%
\mathbb{A})$\textit{.}

\textit{Conversely, if} $S(w)=1_{\mathbb{A}}+\rho w$ \textit{for some} $\rho
\in \mathbb{A}^{-1}$ \textit{and all }$w\in (-1,\infty )\rho ^{-1},$\textit{%
\ then}%
\[
\mathrm{ran}S\supseteq \mathbb{R}_{+}1_{\mathbb{A}}. 
\]

\bigskip

\noindent \textbf{Proof. }Suppose that $\mathrm{ran}S\supseteq \mathbb{R}%
_{+}1_{\mathbb{A}}.$ For $t>0,$ select $w(t)$ with $S(w(t))=t1_{\mathbb{A}}.$
Then%
\begin{eqnarray*}
S(w(s)+w(t)S(w(s)) &=&S(w(s)\circ _{S}w(t))=S(w(s))S(w(t)) \\
&=&st1_{\mathbb{A}}=S(w(st)): \\
w(s)+sw(t) &=&w(st)\text{ mod }\mathcal{N},
\end{eqnarray*}%
by Lemma. 4.1(i). Thus $w$ satisfies $(g_{\mathbb{R}})$ $\func{mod}$ $%
\mathcal{N}.$ Put $c:=w(2).$ Then for $t>0$, with $s=2$, 
\begin{eqnarray*}
w(2)+2w(t) &=&w(2t)=w(t2)=w(t)+tw(2)\text{ mod }\mathcal{N}, \\
w(t) &=&w(2)(t-1)\text{ mod }\mathcal{N}\text{ }=c(t-1)\text{ mod }\mathcal{N%
}, \\
w(t) &=&c(t-1)\text{ }+n(t),\text{ say, with }n(t)\in \mathcal{N}.
\end{eqnarray*}%
So, as $S(n(t))=1_{\mathbb{A}},$ 
\begin{eqnarray*}
S(w(t)) &=&S(n(t)+c(t-1)\text{ })=S(n(t)+c(t-1)S(n(t))) \\
&=&S(n(t))S(c(t-1))=S(c(t-1))=t1_{\mathbb{A}}.
\end{eqnarray*}

For $c\in \mathbb{A}^{-1},$ as $w=c(t-1)$ iff $1_{\mathbb{A}}+c^{-1}w=t1_{%
\mathbb{A}},$ by Prop. 1.2 the remaining assertions are clear.\hfill $%
\square $

\bigskip

\noindent \textbf{Proposition 4.3 }(cf. [BriD, Prop. 3])\textbf{.} \textit{%
For }$S$\textit{\ a solution of }$(GS),$ \textit{if }$1_{\mathbb{A}}\in 
\mathbb{G}_{S}$, \textit{and }$S(\mathbb{G}_{S})$ \textit{contains an
interval on }$\langle 1_{\mathbb{A}}\rangle $ \textit{contiguous with }$1_{%
\mathbb{A}}$\textit{, then} $\mathcal{N}$\textit{\ is a vector subspace;
this is so when }$S(\mathbb{G}_{S})\cap \langle 1_{\mathbb{A}}\rangle $%
\textit{\ is non-meagre on }$\langle 1_{\mathbb{A}}\rangle $\textit{. If }$S(%
\mathbb{G}_{S})\supseteq \mathbb{A}_{1},$ \textit{it is also an ideal; this
is so when }$S(\mathbb{G}_{S})$\textit{\ is non-meagre.}

\bigskip

\noindent \textbf{Proof.} The two assertions follow from $S(\mathbb{G}_{S})%
\mathcal{N}=\mathcal{N}$. For the case of \textit{\ }$\langle 1_{\mathbb{A}%
}\rangle ,$ as $\mathcal{N}$ is additive and $\mathbb{R}$ is the union of
all the iterated vector sums of any non-empty open interval, it follows that 
$a\mathbb{R}1\mathbb{_{\mathbb{A}}\subseteq }\mathcal{N}$, for each $a\in 
\mathcal{N}.$ As $\mathcal{N}$ is additive, this in turn implies that $%
\mathcal{N}$ is a vector subspace. The second assertion is similarly proved,
since for any $a\in \mathcal{N}$ the linear span of $B_{\delta }(1_{\mathbb{A%
}})$ is $\mathbb{A}$, and so $\mathbb{A}\mathcal{N}\subseteq \mathcal{N}$,
i.e. $\mathcal{N}$ is a closed ideal.

The two particular cases asserted follow from the Interior-point Theorem for
category (Steinhaus-Piccard-Pettis theorem [Oxt, Th. 4.8] or [BinO2]):
indeed, $\mathbb{G}_{S}$, as an open subset of a Banach space, is analytic,
so $S(\mathbb{G}_{S}),$ being the continuous image of an analytic set, is
analytic, so has the Baire property, by Nikodym's theorem [Rog]. So $1_{%
\mathbb{A}}\in \mathrm{int}(S(\mathbb{G}_{S})S(\mathbb{G}_{S})^{-1})\mathbb{%
\subseteq }S(\mathbb{G}_{S}),$ the latter inclusion holding because $S(%
\mathbb{G}_{S})$ is a group.\hfill $\square $

\bigskip

We recall that $\mathbb{G}_{S}:=S^{-1}(\mathbb{A}_{1})$ with $\mathbb{A}_{1}$
the connected component of the identity $1_{\mathbb{A}}$. Below our
assumptions imply that $1_{\mathbb{A}}$ is not isolated in $\mathrm{ran}S=S(%
\mathbb{G}_{S}(\mathbb{A})).$ The latter is a natural property, as otherwise 
$S^{-1}(\{1_{\mathbb{A}}\})$ is an open neighbourhood of $0$ on which $%
S\equiv 1_{\mathbb{A}}.$ The analysis below is modelled after that of Prop.
4.2, but with some significant differences which uncover a \textit{spanning
pencil} of isomorphic abelian Popa subgroups $\mathbb{G}_{\rho }(\mathbb{A})$
continuously covering an open set contiguous to $0$. An illustrative example
follows.

\bigskip

\noindent \textbf{Theorem 4.1 (Spanning Pencil Theorem).} \textit{For }$S$%
\textit{\ a solution of }$(GS),$ \textit{if }$1_{\mathbb{A}}$\textit{\ is an
accumulation point of }$\mathrm{ran}S\cap (1_{\mathbb{A}}-\mathrm{ran}S)$%
\textit{, then for some }$c\in \mathbb{A}$%
\begin{equation}
S(c(g-1_{\mathbb{A}}))=g\qquad (g\in S(\mathbb{G}_{S}(\mathbb{A}))). 
\tag{Scale$(c)$}
\end{equation}%
\textit{So if }$c\in \mathbb{A}^{-1}$\textit{, then for }$\rho :=c^{-1}$%
\textit{\ }%
\[
S(w)=1_{\mathbb{A}}+\rho w\qquad (w\in \mathbb{G}_{\rho }(\mathbb{A})=\rho
^{-1}(\mathrm{ran}S-1_{\mathbb{A}})). 
\]%
\textit{In particular,} $\mathbb{G}_{\rho }(\mathbb{A})$ \textit{under }$%
\circ _{\rho }$\textit{\ is an abelian subgroup of }$\mathbb{G}_{S}(\mathbb{A%
}).$

\textit{Furthermore, if also }$\mathrm{Scale}(d)$ \textit{(that is, }$(%
\mathrm{Scale}(c))$\textit{\ with }$c$\textit{\ replaced by }$d$\textit{)
holds for some }$d\in \mathbb{A}$\textit{, then} $c-d\in \mathcal{N}$\textit{%
. So, the sets} $\{\mathbb{G}_{\rho }(\mathbb{A})\cap \mathbb{A}^{-1}:\rho
\in \mathbb{A}^{-1}\}$\textit{\ are mutually disjoint and induce a
continuous partition of }$\mathbb{G}_{S}(\mathbb{A})\cap \mathbb{A}^{-1}$ 
\textit{in some neighbourhood of }$0.$

\textit{Conversely, if} $S(w)=1_{\mathbb{A}}+\rho w$ \textit{for some} $\rho
\in \mathbb{A}^{-1}$ \textit{and all }$w\in \mathbb{G}_{\rho }(\mathbb{A}),$%
\textit{\ then }$1_{\mathbb{A}}$\textit{\ is an accumulation point of }$%
\mathrm{ran}S\cap (1_{\mathbb{A}}-\mathrm{ran}S)$.

\bigskip

\noindent \textbf{Proof. }By assumption we may choose $k\in S(\mathbb{G})$
with both $||k||<1$ and $1_{\mathbb{A}}-k\in S(\mathbb{G});$ then $[1_{%
\mathbb{A}}-k]^{-1}\in S(\mathbb{G}),$ since $\mathrm{ran}S$ is a
multiplicative subgroup of $\mathbb{A}$. By Lemma 4.2(ii),%
\[
\lbrack 1_{\mathbb{A}}-k]^{-1}\mathcal{N\subseteq }S(\mathbb{G})\mathcal{N}=%
\mathcal{N}. 
\]%
For $g\in S(\mathbb{G})\subseteq \mathbb{A}^{-1},$ select $W(g)\in \mathbb{G}
$ with $S(W(g))=g$ and put $c:=W(k)(1_{\mathbb{A}}-k)^{-1}.$ For $g,h\in 
\mathrm{ran}S$%
\begin{eqnarray*}
S(W(g)+W(h)S(W(g)) &=&S(W(g)\circ _{S}W(h))=S(W(g))S(W(h)) \\
&=&gh=S(W(gh)): \\
W(g)+gW(h) &=&W(gh)\text{ mod }\mathcal{N},
\end{eqnarray*}%
by Lemma 4.1(i). Thus $W$ satisfies $(g_{\mathbb{G}})$ $\func{mod}$ $%
\mathcal{N}.$ By commutativity of $\mathbb{A}$ 
\[
S(S(W(g)+W(h)S(W(g)))=gh=hg=S(W(h)+W(g)S(W(h)). 
\]%
So again by Lemma 4.1(i), writing $=_{\mathcal{N}}$ for equality mod $%
\mathcal{N}$, with $k$ for $h$%
\begin{eqnarray*}
W(h)+kW(g) &=&_{\mathcal{N}}\text{ }W(g)+gW(k), \\
W(g)[k-1_{\mathbb{A}}] &=&W(k)[g-1_{\mathbb{A}}]+n(g),\text{ with }n(g)\in 
\mathcal{N}\text{ say:} \\
W(g) &=&[k-1_{\mathbb{A}}]^{-1}W(k)[g-1_{\mathbb{A}}]+[k-1_{\mathbb{A}%
}]^{-1}n(g)=_{\mathcal{N}}c(g-1_{\mathbb{A}}).\text{ }
\end{eqnarray*}%
Write%
\[
W(g)=c(g-1_{\mathbb{A}})+n_{k}(g),\text{ with }n_{k}(g)\in \mathcal{N}; 
\]%
then, as $S(n_{k}(g))=1_{\mathbb{A}},$ 
\begin{eqnarray*}
g &=&S(W(g))=S(n_{k}(g)+c(g-1_{\mathbb{A}})\text{ }%
)=S(n_{k}(g)+S(n_{k}(g))c(g-1_{\mathbb{A}})) \\
&=&S(n_{k}(g))S(c(g-1_{\mathbb{A}}))=S(c(g-1_{\mathbb{A}})).
\end{eqnarray*}%
So $($\textrm{Scale}$\mathcal{(}c\mathcal{))}$ holds.

For $c\in \mathbb{A}^{-1},$ as $w=c(g-1_{\mathbb{A}})$ iff $1_{\mathbb{A}%
}+c^{-1}w=g,$%
\[
S(w)=1_{\mathbb{A}}+c^{-1}w. 
\]%
Hence with $\rho =c^{-1},$ $\mathbb{G}_{\rho }(\mathbb{A})$ is an abelian
subgroup of\textit{\ }$\mathbb{G}_{S}(\mathbb{A}).$

Furthermore, if also $($\textrm{Scale}$\mathcal{(}d\mathcal{))}$ holds for
some $d,$ then 
\[
S(c(g-1_{\mathbb{A}}))=g=S(d(g-1_{\mathbb{A}}))\qquad (g\in S(\mathbb{G}_{S}(%
\mathbb{A}))). 
\]%
As $1_{\mathbb{A}}$ is not isolated in $\mathrm{ran}S\cap (1_{\mathbb{A}}-%
\mathrm{ran}S),$we may take $g$ with $||g||<1$ and $1_{\mathbb{A}}-g\in S(%
\mathbb{G}_{S}(\mathbb{A})).$ As $g-1_{\mathbb{A}}\in \mathbb{A}^{-1}$, by
Lemma 4.1(i) there is $n\in \mathcal{N}$ with%
\[
c(g-1_{\mathbb{A}})=d(g-1_{\mathbb{A}})+n:\qquad c=d+n(g-1_{\mathbb{A}%
})^{-1}=d+n_{g}, 
\]%
with $n_{g}\in \mathcal{N}$ say, as $S(\mathbb{G}_{S}(\mathbb{A}))\mathcal{N}%
=\mathcal{N}$ by Lemma 4.1(ii) (and also as $(g-1_{\mathbb{A}})^{-1}\in S(%
\mathbb{G}_{S}(\mathbb{A}))$, which is a multiplicative group). So $c-d\in 
\mathcal{N}.$

If $w\in c(\mathrm{ran}S-1_{\mathbb{A}})\cap d(\mathrm{ran}S-1_{\mathbb{A}%
})\cap \mathbb{A}^{-1}$ for $c,d\in \mathbb{A}^{-1},$ then%
\[
1_{\mathbb{A}}+c^{-1}w=S(w)=1_{\mathbb{A}}+d^{-1}w:\qquad c=d. 
\]%
First take $V$ an open neighbourhood of $0$ with $V\subseteq S^{-1}(B_{1}(1_{%
\mathbb{A}})).$ Now take $w\in \mathbb{G}_{S}\cap \mathbb{A}^{-1}\cap V;$
then $y=S(w)\in B_{1}(1_{\mathbb{A}}).$ Put $c:=w(y-1_{\mathbb{A}})^{-1}\in 
\mathbb{A}^{-1}.$ Then $y-1_{\mathbb{A}}=c^{-1}w=\rho w,$ say. Then $%
w=c(y-1_{\mathbb{A}})\in c(\mathrm{ran}S-1_{\mathbb{A}})=\mathbb{G}_{\rho }(%
\mathbb{A}),$ and so $w\in \mathbb{G}_{\rho }(\mathbb{A})\cap \mathbb{A}%
^{-1}.$ The map $w\mapsto c=c(w):=$ $w(S(w)-1_{\mathbb{A}})^{-1}$ is
continuous on $\mathbb{G}_{S}\cap \mathbb{A}^{-1}\cap V.$

For the converse with $c=\rho ^{-1}\in \mathbb{A}^{-1}$, note that $\mathbb{G%
}_{S}:=S^{-1}(\mathbb{A}_{1})$ is a non-empty open neighbourhood of $0$.
Choose in $\mathbb{G}_{S}$ non-zero $h_{n}\rightarrow 0$ with $h_{n}\in c(%
\mathrm{ran}S-1_{\mathbb{A}});$ this is possible since $1_{\mathbb{A}}$ is
not isolated and for some $\delta >0,$ $S(B_{\delta }(0))$ is connected and
contains $1_{\mathbb{A}}$ (because $B_{\delta }(0)\subseteq S^{-1}(\mathbb{A}%
_{1})$ for some $\delta >0).$ Put 
\[
g_{n}:=S(h_{n})=1_{\mathbb{A}}+\rho h_{n}\neq 1_{\mathbb{A}};\qquad
k_{n}:=-\rho h_{n}\rightarrow 0. 
\]%
Then, for $w_{n}:=-\rho ^{-1}-h_{n}=-\rho ^{-1}+\rho
^{-1}k_{n}=c(k_{n}-1)\in c(\mathrm{ran}S-1_{\mathbb{A}}),$%
\[
S(w_{n})=1+\rho (-\rho ^{-1}-h_{n})=k_{n}=1-g_{n}:\qquad g_{n}=1-k_{n}. 
\]%
So $g_{n},k_{n}\in \mathrm{ran}S$ and $g_{n}\rightarrow 1_{\mathbb{A}}.$ So $%
1_{\mathbb{A}}$ is a limit point of $\mathrm{ran}S\cap (1_{\mathbb{A}}-%
\mathrm{ran}S).$\hfill $\square $

\bigskip

\noindent \textbf{Corollary 4.4 }(\textbf{Illustrative example). }\textit{Up
to isomorphism,} $\mathbb{G}_{S}(\mathbb{C})$ \textit{is either} $\mathbb{G}%
_{\rho }(\mathbb{C})$ \textit{or a Popa product} $\mathbb{G}_{\alpha }^{\ast
}(\mathbb{R})\times _{\sigma }\mathbb{G}_{\beta }^{\ast }(\mathbb{R}).$

\bigskip

\noindent \textbf{Proof. }See \S 7 (Appendix).\hfill $\square $

\bigskip

In $\mathbb{G}_{\rho }(\mathbb{R)}$, recall that $\mathcal{N=}\{0\},$ so
that $S(a)=1$ for $a\in \mathcal{N}$. Our final result generalizes this to
the Banach-algebra setting. A similar observation in a different context,
and with an altogether different proof, is made in [BriD]. In [BriD, Lemma
6], the lower bound $||S(a)||\geq ||1_{\mathbb{A}}||$ holds for all $a,$
because the context of $GL(\mathbb{A})$ forces $S(a)$ to be an automorphism
(invertible).

\bigskip

\noindent \textbf{Proposition 4.4 (A dichotomy) }(cf. [BriD, Lemma 6]). 
\textit{If }$S$\textit{\ is defined on }$\mathbb{A}$ \textit{and satisfies }$%
(GS)$\textit{\ on }$\mathbb{A}$\textit{, then for any} $a\in \mathbb{A}$, 
\textit{if }$1-S(a)$ \textit{is invertible, then }$S(a(1_{\mathbb{A}%
}-S(a))^{-1})=0$\textit{. In particular,} \textit{either }$||S(a)||\geq ||1_{%
\mathbb{A}}||$ \textit{or }$S(a(1_{\mathbb{A}}-S(a))^{-1})=0.$

\bigskip

\noindent \textbf{Proof.} Assume w.l.o.g. $||1_{\mathbb{A}}||=1,$ and take $%
a\in \mathbb{A}$. Suppose $1_{\mathbb{A}}-S(a)$ is invertible in $\mathbb{A}$%
, so that $S(a)\neq 1_{\mathbb{A}}.$ Take $b:=a(1_{\mathbb{A}}-S(a))^{-1},$
then 
\begin{eqnarray*}
b &=&a+bS(a)=a\circ _{S}b:\quad S(b)=S(a)S(b): \\
0 &=&S(b)(1_{\mathbb{A}}-S(a)):\quad S(b)=0.
\end{eqnarray*}
So $b\notin \mathbb{G}_{S}\mathbb{(A)}$, i.e. $S(b)$ is not invertible, as
claimed.

Now suppose $||S(a))||<1.$ Then [Rud, 10.7] $1_{\mathbb{A}}-S(a)$ is
invertible in $\mathbb{A}$, so $S(b)=0.$\hfill $\square $

\bigskip

\noindent \textbf{Remarks. }1. The final argument above fails if $(GS)$
holds only on $\mathbb{G}_{S}\mathbb{(A)}$, as $b\notin \mathbb{G}_{S}%
\mathbb{(A)}$.

\noindent 2. It is instructive to consider the case of $\mathbb{G}_{\rho }%
\mathbb{(R)}$ with $\rho >0.$ If $|S(a)|=|1+\rho a|<1,$ then automatically $%
a\neq 1_{\rho }=0,$ and $-2\rho ^{-1}<a<0$. Here $1-S(a)=1-(1+\rho a)=-\rho
a $ and so $b:=a(-\rho a)^{-1}=-\rho ^{-1},$ the Popa centre, the only real
that is not a member of $\mathbb{G}_{\rho }^{\ast }\mathbb{(R)}$. So $%
S(b)=0, $ as $0$ is the only non-invertible here. So also $S$ satisfies $%
(GS) $ on all of $\mathbb{R}$. More generally, prompted by the case of $%
\mathbb{R} $- and $\mathbb{C}$-valued functions considered in [Brz1]:

\bigskip

\noindent \textbf{Corollary 4.5} ([Brz1, Cor. 2]). \textit{If }$S$\textit{\
is as in Prop. 4.4 with }$S$ \textit{taking values only in} $\mathbb{A}%
^{-1}, $\textit{\ then }$1_{\mathbb{A}}-S$\textit{\ is never invertible}.%
\textit{\ In particular, if }$S:\langle u\rangle \mathbb{\ \rightarrow R}1_{%
\mathbb{A}} $\textit{\ with }$S$\textit{\ taking only non-zero values, then }%
$S|\langle u\rangle \equiv 1_{\mathbb{A}}.$

\bigskip

\noindent \textbf{Proof.} If $1_{\mathbb{A}}-S(a)$ were invertible for some $%
a$, then, for $b:=a(1-S(a))^{-1},$ $S(b)\notin \mathbb{A}^{-1}$, a
contradiction. So $1_{\mathbb{A}}-S$ is never invertible. In particular, for 
$S:\langle u\rangle \mathbb{\ \rightarrow R}1_{\mathbb{A}},$ since $%
uS(\langle u\rangle )\subseteq \langle u\rangle ,$ take $\mathbb{A}$ to be$%
\mathbb{\ }\langle u\rangle $ and $S$ to be $uS|\langle u\rangle $ to
conclude that as $\mathrm{ran}(S)\subseteq 1_{\mathbb{A}}\mathbb{R}$, $1_{%
\mathbb{A}}-S\equiv 0.$ \hfill $\square $

\bigskip

The next result is distilled from [Brz1, Th. 3] and included here, as it
pursues the linearity theme of $\mathcal{N}$. Key here is a density
argument. Exceptionally, we do not assume that $S$ is continuous; instead,
following [Brz1], we assume that, for each $w,$ if $\langle w\rangle \cap 
\mathbb{G}_{S}(\mathbb{A})$ is non-empty, then it has an interior point --
the \textit{algebraic interior point property} [Lyu, \S 2.2], weaker than $0$
being in the interior of $G_{S}(\mathbb{A})$. By Prop. 1.1 (ii) w.l.o.g. we
may assume below that $0$ is the relevant algebraic interior point.

\bigskip

\noindent \textbf{Corollary 4.6.} \textit{Suppose }$S:\mathbb{A\rightarrow R}%
1_{\mathbb{A}}$ \textit{satisfies }$(GS)$\textit{\ and that, for each} $w,$ 
\textit{a non-empty intersection} $\langle w\rangle \cap \mathbb{G}_{S}(%
\mathbb{A})$\textit{\ has an interior point. Then} $\mathcal{N}$ \textit{is
a vector space.}

\bigskip

\noindent \textbf{Proof.} Suppose otherwise, then $S$ is not identically $1_{%
\mathbb{A}}$ on some line $\langle u\rangle $ with $u\in \mathcal{N}$
(otherwise for $u,v\in \mathcal{N},$ $\langle u\rangle ,\langle v\rangle
\subseteq \mathcal{N},$ and so $su+tv\in \mathcal{N}$ for $s,t\in \mathbb{R}$%
, by additivity of $\mathcal{N}$). So by Cor. 4.5 $S(tu)$ vanishes for some $%
t\in \mathbb{R}.$

By the algebraic interior point property, $\langle u\rangle \cap B_{\delta
}\subseteq $ $\mathbb{G}_{S}(\mathbb{A})$ for some $\delta >0$. By Lemma
4.1(ii), $\mathbb{R}\mathcal{N}=\mathcal{N}$ and so, as $\mathcal{N}$ $\neq
\{0\}$, $\mathcal{N}$ is dense-in-itself (as in Prop. 4.1). So there is $%
su\in \mathcal{N\cap }\langle u\rangle $ with $|(-tu)-su|<\delta $ and so $%
(tu+su)\in \mathbb{G}_{S}(\mathbb{A}).$ But $%
S(su+tu)=S(su+S(su)tu)=S(su)S(tu)=0,$ contradicting that $S(su+tu)$ is
invertible. So $\mathcal{N}$ is a vector space.\hfill $\square $

\section{Banach algebra characterisations}

Our first main result, Theorem 5.1, is an analogue in the Banach-algebra
context of Th. 3.1, which characterises in the Euclidean context the
continuous solutions $S$ of $(GS)$. Here this characterizes continuous
solutions of $(GS)$ over the Popa group $\mathbb{G}_{S}^{\ast }$ of a Banach
algebra as $1+\rho x+N(x),$ where the adjustor $N(x)$ satisfies a Goldie
equation: this incorporation of the Goldie equation bestows on the results
of Section 2 a more satisfactory presentation of the non-linear contribution 
$n(x)$ of $S(x)$ as $N(x)$. Our second main result, Theorem 5.2, identifies
a dichotomy in the behaviour of $N(x):$ it is either linear or exhibits a
curvilinear exponential homogeneity. We term the latter \textit{exponential
tilting}. (It is futher studied in Proposition 5.1.) Finally, in Theorem 5.3
with the hypothesis that $\mathcal{N}\ $is a vector space and $N$ is Fr\'{e}%
chet differentiable, we give a differential characterization of $\mathcal{N}$
and show that $N$ is linear on $\mathcal{N}$ but not beyond. However, it is
more appealing not to assume $\mathcal{N}$ is a vector space: the focus then
shifts to the maximal vector subspace of $\mathcal{N}$, denoted $\mathcal{H}$
below on account of its defining homogeneity property, on which $N$ is
linear. In any case, Prop. 5.3 below connects circumstances of solubility of
the tilting equation ($T$) (below Th. 5.3) to whether $\mathcal{N}$ is a
vector subspace.

We will need Lemma 4.1\ and the following two lemmas.

\bigskip

\noindent \textbf{Lemma 5.1. }\textit{For normed vector spaces }$X,Y$\textit{%
, if }$F:X\rightarrow Y$\ \textit{satisfies:}\newline
\noindent (i)\textit{\ }$F$\textit{\ is Fr\'{e}chet differentiable at every }%
$x\in X$\textit{\ with derivative }$F^{\prime }(x),$\newline
\noindent (ii)\textit{\ }$F(0)=0$\textit{, and}\newline
\noindent (iii)\textit{\ for some continuous linear }$L:X\rightarrow Y$ 
\textit{and} 
\[
F^{\prime }(x)(h)=L(h)\qquad (x,h\in X) 
\]%
-- \textit{then }$F$ \textit{is linear and }$F=L$\textit{.}

\bigskip

\noindent \textbf{Proof. }Fix $u\in X.$ For $t\in \mathbb{R}$ take $%
f(t):=F(tu).$ Then, with $D_{u}$ the directional derivative, 
\[
f^{\prime }(t)=\lim_{s\rightarrow 0}(F((t+s)u)-F(tu))/s=D_{u}F(tu)=L(u). 
\]%
Integrating from $0$ to $t$,%
\[
f(t)-f(0)=L(u)t:\qquad F(tu)=L(u)t. 
\]%
Now take $t=1.$\hfill $\square $

\bigskip

The result below shows that in general the solution of $(GS)$ involves not
only the `canonical example' of a `GS function' from the context of $\mathbb{%
R}$, namely $1+\rho x,$ but also an adjustor $N(x)$ whose characterising
equation is a particular case of the \textit{generalized Goldie equation }$%
(GGE_{SS})$ below (with the subscript indicating that the auxiliaries on the
inside and outside of $N$ are $S$, cf. [Jab1,2]). The adjustor is at best 
\textit{Fr\'{e}chet differentiable} and need not be linear -- see Example
7.3 in the Appendix. (Given \textit{associativity} of the circle operation, $%
(GGE_{SS})$ implies for non-trivial $N$ that $S,$ like $g$ in \S 1, is
multiplicative: it satisfies $(GS).$)

\bigskip

\noindent \textbf{Definition.} Following the notion of $\mathbb{C}$%
-differentiability, say that $f:\mathbb{A}\rightarrow \mathbb{A}$ is $%
\mathbb{A}$\textit{-differentiable} at $a\in \mathbb{A}$ if for some $m\in 
\mathbb{A}$%
\[
\lim_{h\rightarrow 0}h^{-1}[f(a+h)-f(a)]=m,\qquad (h\in \mathbb{A}^{-1}) 
\]%
in which case we will write $m=f^{\prime }(a).$

This is a far stronger property than Fr\'{e}chet differentiability: by Th.
5.1 below $\mathbb{G}_{S}^{\ast }(\mathbb{A})$ has to be abelian if $S$ is $%
\mathbb{A}$-differentiable (cf. Example 7.1-3 in the Appendix); the
preceding definition involves not just the vector but also the ring
structure, as Lemma 5.2 clarifies.

\bigskip

\noindent \textbf{Lemma 5.2}. \textit{For }$\mathbb{A}$ \textit{with dense
invertibles, if} $f:\mathbb{A}\rightarrow \mathbb{A}$ \textit{is continuous
near }$a$ \textit{and} $\mathbb{A}$\textit{-differentiable} \textit{at }$a$%
\textit{\ with }$f^{\prime }(a)=m,$\textit{\ then }$f$\textit{\ is Fr\'{e}%
chet differentiable at }$a$ \textit{with derivative }$f^{\prime }(a)h=mh.$

\bigskip

\noindent \textbf{Proof. }The map $L:h\mapsto mh$ is linear and bounded (as $%
||mh||\leq ||m||.||h||$). For $\varepsilon >0$ and $||h||$ small enough (and
in $\mathbb{A}^{-1}$)%
\[
||h^{-1}[f(a+h)-f(a)-mh]||\leq \varepsilon 
\]%
holds, which implies for $h\in \mathbb{A}^{-1}$ that 
\begin{eqnarray*}
||hh^{-1}[f(a+h)-f(a)-mh]|| &\leq &||h||.||h^{-1}[f(a+h)-f(a)-mh]||: \\
||f(a+h)-f(a)-mh|| &\leq &\varepsilon ||h||,
\end{eqnarray*}%
the latter extending to all small enough $h,$ by density of $\mathbb{A}^{-1}$
and continuity of $f$ near $a$. That is, $f$ is Fr\'{e}chet differentiable
at $a$ with derivative $L.$ \hfill $\square $

\bigskip

In the following theorem the adjustor $N$ is typically linear on $\mathcal{N}%
:$ see Th. 5.3. Its derivative behaves as does that of $S$ in Prop. 2.1.

\bigskip

\noindent \textbf{Theorem 5.1 (First Banach Algebra Characterization Theorem)%
}. \textit{If }$S:\mathbb{G}_{S}^{\ast }$ $\mathbb{\rightarrow A}$\textit{\
satisfies }$(GS_{\mathbb{A}}),$\textit{\ and }$1_{\mathbb{A}}\in \mathbb{G}%
_{S}^{\ast }(\mathbb{A}),$ \textit{then with }$\rho :=S(1_{\mathbb{A}})-1_{%
\mathbb{A}}$\textit{\ there is} $N:\mathbb{G}_{S}^{\ast }\rightarrow 
\mathcal{N}$ \textit{such that}%
\begin{equation}
S(x)=1_{\mathbb{A}}+\rho x+N(x),  \tag{$\ddag $}
\end{equation}%
\textit{where }$N$\textit{\ satisfies the adjustor equation }%
\begin{equation}
N(x+S(x)y)=N(x)+S(x)N(y)\qquad (x,y\in \mathbb{G}_{S}^{\ast }). 
\tag{$GGE_{SS}$}
\end{equation}%
\textit{In particular, }%
\[
N(0)=N(1_{\mathbb{A}})=0. 
\]%
\textit{\ Moreover, if }$N$\textit{\ is Fr\'{e}chet differentiable as a map
into }$\mathcal{N}$\textit{, then its derivative satisfies the similarity
relation}%
\[
N^{\prime }(x)=S(x)N^{\prime }(0)S(x)^{-1}=S(x)N^{\prime
}(0)S(x_{S}^{-1})\qquad (x\in \mathbb{G}_{S}^{\ast }). 
\]%
\textit{\ If }$N$\textit{\ is linear over }$\mathbb{G}_{S}^{\ast }$\textit{,
then for some projection} $\pi :\mathbb{G}_{S}^{\ast }\rightarrow \mathcal{N}
$ \textit{and some linear map }$L,$ \textit{with }$L\pi $ \textit{linear and
injective (into }$\mathcal{N}$\textit{),}%
\[
S(x)=1_{\mathbb{A}}+\rho x+L(\pi (x)). 
\]%
\textit{Furthermore, for }$\mathbb{A}$ \textit{with dense invertibles, if }$%
S $\textit{\ is} $\mathbb{A}$\textit{-differentiable at }$0$\textit{\
(equivalently everywhere), then }$N(x)\equiv 0$ \textit{and }$\mathbb{G}%
_{S}^{\ast }(\mathbb{A})$ \textit{is abelian, since}%
\[
S(x)=1_{\mathbb{A}}+\rho x. 
\]

\noindent \textbf{Proof. }The argument here is an extension of that in Prop.
4.2 and 4.4. As $S(x\circ 1_{\mathbb{A}})=S(1_{\mathbb{A}})S(x)=S(1_{\mathbb{%
A}}\circ x),$ by Lemma 4.1(i), for each $x\in \mathbb{G}_{S}^{\ast }$ there
is $N(x)\in \mathcal{N}$ with%
\begin{eqnarray*}
x+S(x) &=&1_{\mathbb{A}}+xS(1_{\mathbb{A}})+N(x): \\
S(x) &=&1_{\mathbb{A}}+x[S(1_{\mathbb{A}})-1_{\mathbb{A}}]+N(x)=1_{\mathbb{A}%
}+\rho x+N(x).
\end{eqnarray*}%
It now follows that $N(0)=N(1_{\mathbb{A}})=0.$ Substituting for $S$ into $%
(GS):$%
\[
1_{\mathbb{A}}+\rho (x\circ _{S}y)+N(x\circ _{S}y)=S(x\circ _{S}y)=[1_{%
\mathbb{A}}+\rho x+N(x)][1_{\mathbb{A}}+\rho y+N(y)]. 
\]%
Now $\rho (x\circ _{S}y)=\rho (x+yS(x))=\rho (x+y+\rho xy+N(x)y),$ so
multiplying out%
\[
1_{\mathbb{A}}+\rho x+\rho y+\rho ^{2}xy+N(x)\rho y+N(x\circ _{S}y)=[1_{%
\mathbb{A}}+\rho x+N(x)][1_{\mathbb{A}}+\rho y+N(y)]. 
\]%
Re-arrangement gives the asserted equation$.$

Now proceed, as in Prop. 2.1, differentiating $(GGE_{SS})$ w.r.t. $y,$ then
setting $y=0.$ Recalling that $S$ takes invertible values, this gives for $%
x,y,h\in \mathbb{G}_{S}^{\ast }:$%
\begin{eqnarray*}
N^{\prime }(x+S(x)y)(S(x)h) &=&S(x)N^{\prime }(y)h:\quad S(x)^{-1}N^{\prime
}(x+S(x)y)(S(x)h)=N^{\prime }(y)h, \\
S(x)^{-1}N^{\prime }(x)(S(x)h) &=&N^{\prime }(0)h.
\end{eqnarray*}%
Equivalently,%
\[
N^{\prime }(x)=S(x)N^{\prime }(0)S(x)^{-1}=S(x)N^{\prime
}(0)S(x_{S}^{-1})\qquad (x\in \mathbb{G}_{S}^{\ast }), 
\]%
the latter since $1_{\mathbb{A}}=S(x\circ _{S}x_{S}^{-1})=S(x)S(x_{S}^{-1}).$

If $N$ is linear, then $\mathrm{ran}N$ is a vector subspace of $\mathcal{N}$%
. Take $V_{1}:=\ker (N)$ and $V_{0}$ a complementary subspace to $V_{1}.$
Take $\pi $ to be projection onto $V_{0}$ parallel to $V_{1}$ and define $%
L:V_{0}\rightarrow V_{0}$ to be $N|_{V_{0}.}$ For $x\in \mathbb{A}$,%
\begin{eqnarray*}
N(x) &=&N(\pi x+(x-\pi x))=N(\pi x)=L(\pi (x)), \\
S(x) &=&1_{\mathbb{A}}+\rho x+N(x)=1_{\mathbb{A}}+\rho x+L(\pi (x)),
\end{eqnarray*}%
with $L\pi :V_{0}\rightarrow \mathcal{N}$ linear and injective.

Finally, suppose that $\mathbb{A}^{-1}$ is dense in $\mathbb{A}$ and $S$ is $%
\mathbb{A}$-differentiable at $0.$ For $x\in \mathbb{G}_{S}^{\ast }(\mathbb{A%
})$ and $h\in \mathbb{A}^{-1}$ take $k=hS(x)^{-1}\in \mathbb{A}^{-1};$ take $%
h\rightarrow 0,$ then as $||k||\leq ||h||.||S(x)^{-1}||,$ $k\rightarrow 0$
and so $(GS)\ $gives 
\begin{eqnarray*}
h^{-1}[S(x+h)-S(x)] &=&h^{-1}S(x)[S(hS(x)^{-1})-1_{\mathbb{A}}] \\
&=&k^{-1}[S(k)-S(0)]\rightarrow S^{\prime }(0).
\end{eqnarray*}%
Thus $S^{\prime }(x)=S^{\prime }(0).$ So Lemma 5.1 applies at each $x\in 
\mathbb{G}_{S}^{\ast }(\mathbb{A})$ to $F(x):=S(x)-S(0),$ since $S$ is Fr%
\'{e}chet differentiable (by continuity of $S$ and Lemma 5.2) with $%
S^{\prime }(x)h=S^{\prime }(0)h$, giving%
\[
S(x)=S(0)+S^{\prime }(0)x. 
\]%
Taking $x=1_{\mathbb{A}}$ gives $S^{\prime }(0)=\rho $ and so $S(x)=1+\rho
x, $ which is equivalent to commutativity of $\circ _{S},$ as noted in the
Remark 4 after Cor. 1.1. \hfill $\square $

\bigskip

\noindent \textbf{Corollary 5.2. }(i)\textbf{\ }\textit{If} $\mathcal{N=\{}%
0\},$ \textit{then }$N(x)\equiv 0$ \textit{and so} 
\[
S(x)=1_{\mathbb{A}}+\rho x\text{ \qquad for }\rho :=S(1_{\mathbb{A}})-1_{%
\mathbb{A}} 
\]%
\textit{Thus, for} $\mathbb{A}=\mathbb{R}$ \textit{positive solutions of} $%
(GS_{\mathbb{R}})$\textit{\ take this form with }$\rho \geq 0.$

\noindent (ii) \textit{The adjustor }$N$\textit{\ takes the form }$kx$%
\textit{\ with }$k\in \mathbb{A}$\textit{\ iff }$\circ _{S}$ \textit{is
commutative.}

\bigskip

\noindent \textbf{Proof. }(i) The first assertion is clear. For $\mathbb{A}=%
\mathbb{R}$, since $S$ takes values in $\mathbb{R}^{-1}=\mathbb{R}\backslash 
\mathbb{\{}0\mathbb{\}}$, by Cor 4.2 either $\mathcal{N}=\mathbb{R}$ and
then $S(x)\equiv 1_{\mathbb{A}}$ or $\mathcal{N=\{}0\}.$

(ii) If $\circ _{S}$ is commutative, then by Remark 4 after Cor. 1.1. above, 
$S(x)=1_{\mathbb{A}}+\sigma x$ for some $\sigma ;$ then $1_{\mathbb{A}}+\rho
x+N(x)=1_{\mathbb{A}}+\sigma x,$ and so $N(x)=(\sigma -\rho )x$ is linear.
(This also follows directly from $(GGE_{SS}).)$ Conversely, if $N(x)=kx$ for
some $k,$ then $S(x)=1_{\mathbb{A}}+\sigma x$ for $\sigma :=(\rho +k),$ so $%
\circ _{S}$ is commutative.\hfill $\square $

\bigskip

Our next result, Th. 5.2, draws on a result in [BinO4] where the context
refers to real-valued functions on a topological vector space $X$; however,
the calculations there may be reinterpreted upon replacing $X$ there by $%
\mathbb{A}$ here, thereby introducing also $\mathbb{A}$-valued functions and
using symbolic calculus (Riesz-Dunford functional calculus). This requires
relevant elements of $\mathbb{A}$, such as $e^{\gamma (u)}-1_{\mathbb{A}}$
below for $\gamma (u):=S_{u}^{\prime }(0)$, the directional derivative, to
be invertible; see e.g. [Rud, Ch. 10]. Below $e^{\gamma (u)}-1_{\mathbb{A}}$
will be invertible iff $e^{\lambda }\neq 1$ for all $\lambda $ in the
spectrum of $\gamma (u)$ [Rud, Th. 10.28]. See Example 7.3 in the Appendix.
However, the analysis here departs from [BinO4] in requiring the stronger
assumption of Fr\'{e}chet differentiability. To use symbolic calculus in
Theorem 5.2 below, we need the following convergence result.

\bigskip

\noindent \textbf{Lemma 5.3.} (i) \textit{For }$\mathbb{C}^{\dag }$ \textit{%
an unbounded, connected open subset of the complex plane containing the
origin and for a sequence of positive reals }$t(n)\rightarrow t>0:$ \textit{%
\ }%
\[
\frac{(1+z/n)^{nt(n)}-1}{(1+z/n)^{n}-1}\rightarrow \frac{e^{tz}-1}{e^{z}-1},%
\text{ or }t\text{ if }e^{z}=1, 
\]%
\textit{convergence being uniform on compact subsets of} $\mathbb{C}^{\dag
}. $

(ii) \textit{The map }$z\mapsto \mu (z):=(e^{z}-1)/z,$\textit{\ }or\textit{\ 
}$1$\textit{\ }if\textit{\ }$z=0,$\textit{\ is holomorphic and invertible
near }$(w,z)=(1,0)$\textit{.}

\bigskip

\noindent \textbf{Proof.} (i) The assumptions allow the use of a logarithm
on $\mathbb{C}^{\dag }.$ So w.l.o.g assume $\mathbb{C}^{\dag }:=\mathbb{C}%
\backslash (-\infty ,0]$. For $t>0,$ the function $g_{t}$, defined for $%
\zeta \in $ $\mathbb{C}^{\dag }$ by%
\[
g_{t}(\zeta )=\frac{e^{t\log \zeta }-1}{\zeta -1}\text{ with }g_{t}(1)=t, 
\]%
is differentiable, because it is differentiable at $\zeta =1:$ using
L'Hospital twice,%
\[
\lim_{w\rightarrow 0}\frac{g_{t}(1+w)-t}{w}=\frac{t(t-1)}{2}. 
\]%
Hence%
\[
f(\zeta ):=g_{t}(e^{\zeta })=\frac{e^{t\zeta }-1}{e^{\zeta }-1},\text{ or }t%
\text{ if }e^{\zeta }=1, 
\]%
is holomorphic, as is for each $n=1,2,...$

\[
f_{n}(\zeta ):=g_{t(n)}((1+\zeta /n)^{n})=\frac{(1+\zeta /n)^{nt(n)}-1}{%
(1+\zeta /n)^{n}-1},\text{ or }t(n)\text{ if }(1+\zeta /n)^{n}=1. 
\]%
Furthermore, since $(1+\zeta /n)^{n}\rightarrow e^{\zeta }$ and $%
t(n)\rightarrow t>0,$ then for each $\zeta \in \mathbb{C}^{\dag }$ 
\[
f_{n}(\zeta )\rightarrow f(\zeta )=\frac{e^{t\zeta }-1}{e^{\zeta }-1},\text{
or }t\text{ if }e^{\zeta }=1, 
\]%
convergence being uniform on compact subsets of $\mathbb{C}^{\dag }.$

(ii) The first part is clear from the preceding paragraph. As for
invertibility of $w=w(\zeta )$ near $w=1$, take $F(\omega ,\zeta ):=\omega
-\mu (\zeta ),$ which is analytic (in the two complex variables sense) with $%
F(1,0)=0$, and note that $F_{\zeta }=-(1+e^{\zeta }(\zeta -1))/\zeta ^{2}$
or $-1/2$ if $\zeta =0.$ By the Implicit Function Theorem [Gam1, Ch. 3]
there is a solution $z=z(\omega )$ near $\omega =1$ with $z(1)=0.$ \hfill $%
\square $

\bigskip

\noindent \textbf{Corollary 5.3.} \textit{With }$\mathbb{C}^{\dag }$ \textit{%
as in Lemma 5.3, in any Banach algebra with dense invertibles, for }$a$ 
\textit{an element with spectrum satisfying} $\mathrm{spec}(a)\subseteq 
\mathbb{C}^{\dag }\cup \{0\}:$%
\[
\frac{(1+a/n)^{nt(n)}-1}{(1+a/n)^{n}-1}\rightarrow \frac{e^{ta}-1}{e^{a}-1},%
\text{ or }t\text{ if }1\in \exp (\mathrm{spec}(a)). 
\]

\bigskip

\noindent \textbf{Proof.} If $0\notin \mathrm{spec}(a),$ the result follows
from Lemma 5.3 by [Rud, Th. 10.27]. If $0\in \mathrm{spec}(a),$ then $a$ is
not invertible. Assuming dense invertibles, choose invertible elements with $%
a_{k}\rightarrow a.$ Again by Lemma 5.3 and [Rud, Th. 10.27] for each $k$%
\[
\frac{(1+a_{k}/n)^{nt(n)}-1}{(1+a_{k}/n)^{n}-1}\rightarrow \frac{e^{ta_{k}}-1%
}{e^{a_{k}}-1},\text{ or }t\text{ if }1\in \exp (\mathrm{spec}(a_{k})). 
\]%
Now choose $k(n)$ so that $(1+a_{k(n)}/n)^{n}\rightarrow a.$ So if $\mathrm{%
spec}(a)\subseteq \mathbb{C}^{\dag },$ then%
\[
\frac{(1+a_{k(n)}/n)^{nt(n)}-1}{(1+a_{k(n)}/n)^{n}-1}\rightarrow \frac{%
e^{ta}-1}{e^{a}-1},\text{ or }t\text{ if }1\in \exp (\mathrm{spec}(a)). 
\]%
and so%
\[
\frac{(1+a/n)^{nt(n)}-1}{(1+a/n)^{n}-1}\rightarrow \frac{e^{ta}-1}{e^{a}-1},%
\text{ or }t\text{ if }1\in \exp (\mathrm{spec}(a)). 
\]%
Indeed, if $1\notin \exp (\mathrm{spec}(a)),$ then $e^{a}-1$ is invertible
and so also for large $n$ is $(1+a/n)^{n}-1$, and so also is its approximand 
$(1+a_{k(n)}/n)^{n}-1.$\hfill $\square $

\bigskip

Recall that the spectrum of an element is compact. So its complement is
open: in order to have a logarithm available for Lemma 5.3 to hold, the
connected component of $0$ must be unbounded; one may term this a \textit{no
encirclement} condition.

\bigskip

\noindent \textbf{Theorem 5.2 (Second Characterization Theorem: Curvilinear
exponential homogeneity)}. \textit{Suppose that the solution }$S$\textit{\
to }$(GS)$ \textit{is Fr\'{e}chet differentiable and that }$N$\textit{\
solves the Goldie equation in} $\mathbb{A}:$%
\[
N(x+S(x)y)=N(x)+S(x)N(y). 
\]%
\textit{For any }$u$\ \textit{with }$\gamma (u):=S_{u}^{\prime }(0)$\textit{%
\ having a spectrum not separating }$0$\textit{\ from }$\infty ,$ \textit{%
the radiality formulae below hold for }$t\geq 0$\textit{: } 
\begin{equation}
N(u(e^{t\gamma (u)}-1)/\gamma (u))=\lambda _{u}(t)N(u(e^{\gamma
(u)}-1)/\gamma (u)),  \tag{\QTR{rm}{Rad}}
\end{equation}%
\textit{\ for }%
\[
\lambda _{u}(t):=(e^{t\gamma (u)}-1)/(e^{\gamma (u)}-1)\text{,} 
\]%
\textit{with the L'Hospital convention that, when }$1\in \exp (\mathrm{spec}%
(\gamma (u))),$%
\[
N(tu)=t1_{\mathbb{A}}N(u)=tN(u)\qquad (t\geq 0). 
\]%
\textit{Furthermore, the exponential tilting map (with }$\mu $\textit{\ as
above) }%
\[
T(u):=u\mu (\gamma (u))=u(e^{\gamma (u)}-1)/\gamma (u):\qquad T(tu)=\lambda
_{u}(t)T(u)\qquad (t\geq 0) 
\]%
\textit{\ has invertible multiplier }$\mu (\gamma (u))$\textit{\ for }$u\in 
\mathbb{G}_{S}(\mathbb{A}),$\textit{\ and exhibits the `curvilinear
exponential homogeneity' under }$N$\textit{:}%
\[
N(T(tu))=N(\lambda _{u}(t)T(u))=\lambda _{u}(t)N(T(u))\qquad (t\geq 0). 
\]

\noindent \textbf{Proof.} Referring to the polynomials $\wp
_{n}(x)=1+x+...+x^{n-1}$ and rational polynomials $[\wp _{m}/\wp _{n}](x),$
[BinO4] Lemma 5 gives, taking $g=h=S,$ that 
\begin{equation}
N(\wp _{m}(S(u/n))u/n)=[\wp _{m}/\wp _{n}](S(u/n))N(\wp _{n}(S(u/n))u/n), 
\tag{\QTR{rm}{Pre}$Rad$}
\end{equation}%
for any $u.$ By Fr\'{e}chet differentiability, write 
\[
S(u/n)-S(0)=S^{\prime }(0)(u/n)+\varepsilon _{n}(u), 
\]%
with $n\varepsilon _{n}(u)\rightarrow 0.$ Put $\xi _{n}=S(u/n)$ and $\gamma
(u):=S^{\prime }(0)u.$ Note that $\gamma (u)/n\rightarrow 0,$ as $||\gamma
(u)||\leq ||S^{\prime }(0)||.||u||.$ For all $n$ large enough $1_{\mathbb{A}%
}+\gamma (u)/n+\varepsilon _{n}(u)\in \mathbb{A}_{1},$ the connected
component of unity; then $\eta _{n}(u):=n\log [1_{\mathbb{A}}+\gamma
(u)/n+\varepsilon _{n}(u)]$ is well defined ([Rud,10.43c], [Ric, 1.4.12])
and $\eta _{n}(u)\rightarrow 0.$ Fix $t>0$ and choose $m=m(n)$ so that $%
m/n=m(n)/n=t(n)\rightarrow t.$ By Lemma 5.3,%
\begin{eqnarray*}
\lbrack \wp _{m}/\wp _{n}](\xi _{n}) &=&\frac{\exp \{t(n)[n\log [1_{\mathbb{A%
}}+S_{u}^{\prime }(0)(1/n)+\varepsilon _{n}(u)]\}-1}{\exp \{[n\log [1_{%
\mathbb{A}}+S_{u}^{\prime }(0)(1/n)+\varepsilon _{n}(u)]\}-1} \\
&=&\frac{\exp (t(n)\eta _{n}(u))-1}{\exp (\eta (n))-1}\rightarrow \frac{%
e^{\gamma (u)t}-1}{e^{\gamma (u)}-1},\text{ or }t\text{ if }1\in \exp (%
\mathrm{spec}(\gamma (u))).
\end{eqnarray*}%
Likewise, provided $\xi _{n}\neq 1$ 
\[
\wp _{m}(\xi _{n})/n=\frac{\exp \{t(n)[n\log [1_{\mathbb{A}}+S_{u}^{\prime
}(0)(1/n)+\varepsilon _{n}(u)]\}-1}{n(\xi _{n}-1)} 
\]

\[
=\frac{\exp (t(n)\eta _{n}(u))-1}{S^{\prime }(0)u+n\varepsilon _{n}(u)}%
\rightarrow \mu (\gamma (u))=\frac{e^{\gamma (u)t}-1}{\gamma (u)},\text{ or }%
t\text{ if }1\in \exp (\mathrm{spec}(\gamma (u))). 
\]%
The equations now follow from $($\textrm{Pre}$Rad)$ and Corollary 1, the
case $t=0$ being trivial.

As $\gamma (u)$ is (linear and so) homogeneous,%
\begin{eqnarray*}
T(tu) &=&u(e^{t\gamma (u)}-1)/\gamma (u)=u\mu (\gamma (u))\lambda
_{u}(t)=\lambda _{u}(t)T(u), \\
N(u(e^{t\gamma (u)}-1)/\gamma (u)) &=&N(u\lambda _{u}(t)(e^{\gamma
(u)}-1)/\gamma (u))=\lambda _{u}(t)N(u(e^{\gamma (u)}-1)/\gamma (u)).
\end{eqnarray*}%
By Lemma 5.4(ii) $\mu (\gamma (u))$ is invertible in $\mathbb{A}.$\hfill $%
\square $

\bigskip

\noindent \textbf{Remarks: }1. \textit{Exponential tilting}\textbf{. }One
may interpret the adjustor $N$ as comprising a linear action on $\mathcal{N}$
and, on complementary directions, a homogeneous action after the \textit{%
tilting} $T$ with (exponential)\textit{\ scaling }$\lambda _{u}$ -- so a
kind of shearing. (The term exponential tilting here is borrowed from
probability theory, where it is used as a synonym for the Esscher transform
of collective risk theory [GerS].)

\noindent 2. Interestingly, $\lambda _{u}$ satisfies a pair of Goldie
equations with parameter $\gamma (u)$:%
\begin{eqnarray*}
\lambda _{u}(s+t) &=&\lambda _{u}(s)+e^{\gamma (u)s}\lambda _{u}(t)\qquad
(s,t\in \mathbb{R}); \\
\lambda _{u+v}(t) &=&\lambda _{u}(t)+e^{\gamma (u)t}\lambda _{v}(t)\qquad
(t\in \mathbb{R}).
\end{eqnarray*}%
We now study the interplay between $\mathcal{N}$ and the set of $N$%
-invariant directions, which we denote by $\mathcal{H}$ for `homogeneous':%
\[
\mathcal{H}:=\{x:(\forall t\in \mathbb{R})N(tx)=tN(x)\}. 
\]%
This is closed by continuity of $N$. If $x\in \mathcal{H}$, then, for $%
s,t\in \mathbb{R}$, $N(t(sx))=tsN(x)=tN(sx),$ so $sx\in \mathcal{H}$, i.e. $%
\langle x\rangle \subseteq \mathcal{H}$: so $\mathcal{H}$ is homogeneous.

\bigskip

\noindent \textbf{Lemma 5.4.} (i) \textit{For }$S$ \textit{the solution to }$%
(GS)$ \textit{and }$N$\textit{\ its adjustor, by }$(\ddag )$%
\[
\mathcal{N}=\{x:N(x)=-\rho x\}, 
\]%
\textit{so }$N$\textit{\ is additive on }$\mathcal{N}$\textit{, and so when }%
$1\in \exp (\mathrm{spec}(\gamma (u)))$\textit{, then}%
\[
N(tu)=tN(u)\qquad (u\in \mathcal{N},t\in \mathbb{R}). 
\]%
(ii) \textit{Furthermore,}%
\[
\mathcal{N}\subseteq \mathcal{H}\text{ iff }\mathcal{N}\text{ is a vector
subspace.} 
\]%
\textit{That is, }$\mathcal{N}$\textit{\ is a vector subspace iff }$N$%
\textit{\ acts homogeneously on }$\mathcal{N}$\textit{\ iff }$N$\textit{\ is
linear on }$\mathcal{N}$\textit{.}

\noindent (iii) \textit{If }$\langle u\rangle \subseteq \mathcal{N}$\textit{%
, then }$\langle u\rangle \subseteq \mathcal{H}$\textit{.}

\noindent (iv) $\mathcal{H}\cap \mathcal{N}$\textit{\ is a vector subspace
(so if }$u\in \mathcal{H}\cap \mathcal{N},$\textit{\ then }$\langle u\rangle
\subseteq \mathcal{N}$\textit{).}

\bigskip

\noindent \textbf{Proof.} (i) Indeed $S(x)=1_{\mathbb{A}}+\rho x+N(x)$
implies the equivalence of $S(x)=1$ and $N(x)=-\rho x.$ For $u,v\in \mathcal{%
N}$, as $u+v\in \mathcal{N}$,%
\[
N(u+v)=\rho (u+v)=\rho u+\rho v=N(u)+N(v). 
\]%
In particular, for $u\in \mathcal{N}$, $0=N(u-u)=N(u)+N(-u)$, i.e. $%
N(-u)=-N(u),$ and so since $1\in \exp (\mathrm{spec}(\gamma (u)))$ implies $%
1\in \exp (\mathrm{spec}(\gamma (-u))),$ the claim of linearity follows from
Th. 5.2.

(ii) Suppose that $\mathcal{N}\subseteq \mathcal{H}$. For $u,v\in \mathcal{N}
$ and positive integers $p,q,r,$ additivity gives $pu+qv\in \mathcal{N}$,
from which it follows that $(pu+qv)/r\in \mathcal{H}$, by homogeneity of $%
\mathcal{H},$ i.e. that $su+tv\in \mathcal{H}$ for $s,t$ in $\mathbb{Q}_{+}$%
, and so, as $\mathcal{H}$ is closed, that $su+tv\in \mathcal{H}$ for $s,t$
in $\mathbb{R}_{+}$. Further, with $u,v,p,q,r$ as above, as $pu+qv\in 
\mathcal{H}$, by additivity of $N$ on $\mathcal{N}$,%
\begin{eqnarray*}
N((pu+qv)/r) &=&N(pu+qv)/r=(pN(u)+qN(v))/r: \\
N((pu+qv)/r) &=&(p/r)N(u)+(q/r)N(v): \\
N(su+tv) &=&sN(u)+tN(v),\text{ for }s,t\in \mathbb{R}_{+}\text{,}
\end{eqnarray*}%
by continuity of $N.$ So, as $N(u)=-\rho u$ and $N(v)=-\rho v,$ for $s,t\in 
\mathbb{R}_{+}$,%
\begin{eqnarray*}
S(su+tv) &=&1_{\mathbb{A}}+\rho (su+tv)+N(su+tv) \\
&=&1_{\mathbb{A}}+s(\rho u+N(u))+t(\rho v+N(v))=1_{\mathbb{A}}.
\end{eqnarray*}%
That is, $su+tv\in \mathcal{N}$. But $-\mathcal{N}=\mathcal{N}$, so $%
\mathcal{N}$ is a vector space.

Conversely, suppose $\mathcal{N}$ is a vector space. For $u\in \mathcal{N}$,
and $t\in \mathbb{R},$ as $tu\in \mathcal{N}$, $N(tu)=-\rho tu=t(-\rho
u)=tN(u)$ and so $u\in \mathcal{H}$. That is $\mathcal{N\subseteq H}$.

If $N$ acts homogeneously on $x\in $ $\mathcal{N}$, i.e. $\mathcal{N}%
\subseteq \mathcal{H}$, then for $x\in \mathcal{N}$ and $t\in \mathbb{R}$, $%
N(tx)=tN(x)=-t\rho x$ and so $tx\in \mathcal{N}$; so $\mathcal{N}$, being an
additive subgroup, is a vector subspace. Conversely, if $\mathcal{N}$ is a
vector subspace, then, for $t\in \mathbb{R}$ and $x\in \mathcal{N}$, $%
tN(x)=-t\rho x=N(tx),$ as $tx\in \mathcal{N}$; so $N$ acts homogeneously on $%
\mathcal{N}$.

As for the final condition, $N,$ being additive on $\mathcal{N}$ by $%
(GGE_{SS}),$ is linear on $\mathcal{N}$ iff it acts homogeneously on $%
\mathcal{N}$.

(iii) If there is $u$ with $\langle u\rangle \subseteq \mathcal{N}$ (e.g. in
a finite-dimensional setting where this holds by Theorem B of \S 4), then,
for $t\in \mathbb{R}$, $N(tu)=-t\rho u=tN(u),$ since $u$ and $tu\in $\textit{%
\ }$\mathcal{N}$. That is, $\langle u\rangle \subseteq \mathcal{H}$ and so $%
\mathcal{H}$ is non-empty.

(iv) Take\textit{\ }$u,v\in \mathcal{H}\cap \mathcal{N}$ and $t\in \mathbb{R}
$; then, as $N(tu)=tN(u)$ and by (i), 
\begin{eqnarray*}
S(tu) &=&1_{\mathbb{A}}+\rho tu+N(tu)=1_{\mathbb{A}}+t(\rho u+N(u))=1_{%
\mathbb{A}}, \\
N(t(u+v)) &=&t\rho (u+v)=t(N(u)+N(v))=t(N(u+v)),
\end{eqnarray*}%
so $\langle u\rangle \subseteq \mathcal{H}\cap \mathcal{N}$ and $u+v\in $ $%
\mathcal{H}\cap \mathcal{N}$, which is thus a vector subspace.\hfill $%
\square $

\bigskip

It is routine to show, by reference to additivity and closure of $\mathcal{N}
$, using rational scalars as above, that $\langle \mathcal{N}\rangle
=\{tx:x\in \mathcal{N},t\in \mathbb{R}\}$ (cf. [Brz1, Proof of Th. 3]),
hence the interplay of $\mathcal{H}$ and $\mathcal{N}$. The next two results
again require the hypothesis of Th. 5.2 on non-encirclement of the origin by
the spectrum. Here we are about to see that $\mathcal{H}\subseteq \mathcal{N}
$, i.e. the reverse inclusion to that of Lemma 5.4. It emerges below that $%
\mathcal{H}=\mathcal{N}$ if $\mathcal{N}$ is a vector subspace upon which $%
N\ $is linear. Since $N(\mathcal{N})\mathcal{\subseteq N}$ (Th. 5.1), this
is in line with a similar result in the First Popa Homomorphism Theorem, Th.
4A of [BinO3], where the homomorphism is linear on $\mathcal{N=N}(\rho ).$
The result below should be compared to the identity $\mathcal{N}%
=\{x:N(x)=-\rho x\}$ of Lemma 5.4. That $\mathcal{H}=\mathcal{N}$ when $%
\mathbb{A}=\mathbb{R}^{d}$ will eventually follow. For $\mathcal{N}_{\gamma
} $ below, see Cor 2.1.

\bigskip

\noindent \textbf{Theorem 5.3 (Third -- differential -- characterization of }%
$\mathcal{N}$\textbf{).} \textit{For }$S$\textit{\ Fr\'{e}chet
differentiable at }$0$ \textit{and }$\gamma =S^{\prime }(0)$\textit{, if
each element }$\gamma u$\textit{\ has spectrum not separating }$0$\textit{\
from }$\infty ,$ \textit{then} 
\[
\mathcal{H}=\mathcal{N}_{\gamma }=\{u:S^{\prime }(0)u=0\}=\{u:N^{\prime
}(0)u=-\rho u\}\mathcal{\subseteq N}, 
\]%
\textit{so that }$\mathcal{H}$\textit{\ is a vector subspace and }$N$\textit{%
\ is linear on }$\mathcal{H}$\textit{; furthermore, }$\mathcal{H}$\textit{\
is the maximal vector subspace of }$\mathcal{N}$\textit{\ such that if }$%
\mathcal{N}$ \textit{is a vector subspace, then }$\mathcal{N=H}$\textit{\
and }$N$ \textit{is linear on }$\mathcal{N}$\textit{.}

\textit{For }$\mathcal{M}$ \textit{any subspace complementary to }$\mathcal{N%
}$:%
\begin{equation}
N(x)=N(\pi _{\mathcal{N}}(x))+N(\pi _{\mathcal{M}}(x)),  \tag{$+$}
\end{equation}%
\textit{where }$\pi $ \textit{refers to the corresponding projections.}

\bigskip

\noindent \textbf{Proof.} By Lemma 5.4 and Theorem 5.2 and since $\gamma
(u)=S^{\prime }(0)u$ (Fr\'{e}chet differentiability at\textit{\ }$0):$ 
\begin{eqnarray*}
\mathcal{H} &=&\{u:N(tu)=tN(u)(\forall t)\}=\{u:\lambda _{u}(t)=t\mathbf{1}%
(\forall t)\}=\{u:\gamma (u)=0\} \\
&=&\{u:S^{\prime }(0)u=0\}.
\end{eqnarray*}%
So for $u\in \mathcal{H}$, as $DS(0)u=0$,%
\[
0=\lim_{t\rightarrow 0}\frac{S(tu)-1_{\mathbb{A}}}{t}=\lim_{t\rightarrow 0}%
\frac{\rho tu+N(tu)}{t}=\rho u+N(u), 
\]%
and so $u\in \mathcal{N}$ by Lemma 5.4. Thus $\mathcal{H\subseteq N}$ .
Since $\gamma $ is linear, $\mathcal{H}$ as its kernel is a vector subspace.
(Also by Lemma 5.4(iv), as $\mathcal{H=H}\cap \mathcal{N}$.)

We may now deduce from its additivity on $\mathcal{N}$ that $N$ is linear on 
$\mathcal{H}$: for $s,t\in \mathbb{R}$ and $u,v$ $\in \mathcal{H}$, as $%
su,tv\in \mathcal{H}\subseteq \mathcal{N}$ by additivity on $\mathcal{N}$%
\[
N(su+tv)=N(su)+N(tv)=sN(u)+tN(v), 
\]%
the last since $\mathcal{H}$ is homogeneous for $N.$

For $\mathcal{N}$ a vector subspace, $\mathcal{N}\subseteq \mathcal{H}$
(Lemma 5.4); but $\mathcal{H}\subseteq \mathcal{N}$, so $\mathcal{N}=%
\mathcal{H}$.

Given $\mathcal{M}$ complementary to $\mathcal{N}$ in $\mathbb{A}$, for $%
x\in $ $\mathcal{N}$ and $y\in \mathcal{M}$, ($+$) follows from $(GGE_{SS}),$
as $S(x)=1$ and $x+y=x\circ _{S}y.$\hfill $\square $

\bigskip

In the Euclidean case, the spectrum of $\gamma (u)$ is finite, hence does
not separate 0 from infinity; however, the proof that $\mathcal{H}=\mathcal{N%
}$ when $\mathbb{A}=\mathbb{R}^{d}$, must wait until we have studied the
tilting map of Th. 5.3 and the equation%
\begin{equation}
v=T(u):=u(e^{\gamma (u)}-1_{\mathbb{A}})/\gamma (u).  \tag{$T$}
\end{equation}

If $\gamma $ is $\mathbb{A}$-homogeneous (for which see Section 2 above),
the map $T$ is invertible on a part of its range and with the \textit{%
explicit inverse} of Prop. 5.1.

As for the conclusion below about the range of $T$ covering a nhd of the
origin, the simplest case of $\mathbb{A=R}$ and $\gamma (u)=\gamma u,$ with $%
\gamma >0$ say, is illuminating. For given $v,$ $v=T(u)$ is soluble only for 
$-1<\gamma v,$ as%
\[
(e^{\gamma u}-1)/\gamma =v:\qquad e^{\gamma u}=1+\gamma v. 
\]%
So the range of $T$ is bounded in one direction.

This uni-directional scenario above emerges in a more general setting,
enabled by a uniqueness/identity result, kindly established for us by Amol
Sasane. For the proof see \S 7 (Appendix). We apply it to $f(\zeta )=\log
(1+\zeta )$ with domain $\mathbb{C}\backslash (-\infty ,0]$ to validate a
condition placed on $f(\gamma (u)).$ Recall that $\mathbb{A}$ is termed 
\textit{semisimple} if its Gelfand transform has trivial kernel, i.e. is
injective [Rud, 11.9], [Con, Ch. VII \S 8], an instance being $C[0,1]$ (loc.
cit.).

\bigskip

\noindent \textbf{Theorem S (Uniqueness/Identity).} \textit{For} $\mathbb{A}$
\textit{semisimple, }$f(\zeta )$\textit{\ holomorphic on its domain,} $u\in 
\mathbb{A}$, \textit{and}%
\[
\mathcal{D}_{u}:=\{\zeta \in \mathbb{C}:\mathrm{spec}(\zeta \gamma
(u))\subseteq \mathrm{dom}f\}
\]%
\textit{non-empty, open and connected: if, for some set }$\Sigma \subseteq 
\mathcal{D}_{u}$\textit{\ with a limit point in }$\mathcal{D}_{u}$\textit{,} 
\textit{the identity }%
\[
f(\zeta \gamma (u))=\gamma (f(\zeta \gamma (u)\cdot u/\gamma (u))
\]%
\textit{holds for all }$\zeta \in \Sigma $\textit{, then it holds also for
all }$\zeta \in \mathcal{D}_{u}$\textit{.}

\bigskip

\noindent \textbf{Proposition 5.1}. \textit{For a derivative }$\gamma
=S^{\prime }(0)$\textit{\ that is }$\omega $-\textit{homogeneous : if }$%
v=T(u)$\textit{\ for }$u$\textit{\ and }$v$ \textit{with }$\gamma (u),\gamma
(v)\in \mathbb{A}^{-1},$\textit{\ and }$t>0,$ \textit{then}: (i)%
\begin{equation}
u=u(v):=v\log (1_{\mathbb{A}}+\gamma (v))/\gamma (v),  \tag{$T$-inv}
\end{equation}%
\textit{so that }$u$\textit{\ is uniquely determined, and further}: (ii)%
\textit{\ }$1+\gamma (v)\in \exp (\mathbb{A})=\mathbb{A}_{1}$ \textit{and}%
\[
\gamma (u)=\log (1_{\mathbb{A}}+\gamma (v))\in \text{\textrm{Ran}}(\gamma ). 
\]%
\textit{Conversely, for a general derivative }$\gamma $\textit{, if }(ii)%
\textit{\ holds for some }$v$\textit{\ with }$u=u(v)$\textit{\ from }(i)%
\textit{, then }$v=T(v(u)).$ \textit{This condition holds with }$\mathbb{A}$ 
\textit{semi-simple for an }$\omega $-\textit{homogeneous }$\gamma ,$\textit{%
\ and, furthermore, the range of }$T$\textit{\ contains the neighbourhood }$%
\{v:||\gamma (v)||<1\}$ \textit{of} $0.$

\bigskip

\noindent \textbf{Proof. }By Corollary 2.2, if $v=T(u),$ then%
\begin{eqnarray*}
\gamma (v) &=&\gamma (u(e^{\gamma (u)}-1_{\mathbb{A}})/\gamma (u))=e^{\gamma
(u)}-1_{\mathbb{A}}: \\
\gamma (u) &=&\log (1_{\mathbb{A}}+\gamma (v)): \\
u &=&v\gamma (u)/(e^{\gamma (u)}-1_{\mathbb{A}}):u=v\log (1_{\mathbb{A}%
}+\gamma (v))/\gamma (v),
\end{eqnarray*}%
giving (i); the first line above giving (ii).

Conversely, if $\gamma (u(v))=\log (1_{\mathbb{A}}+\gamma (v)),$ then $%
e^{\gamma (u(v))}-1_{\mathbb{A}}=\gamma (v)$ and so 
\[
T(u(v))=u(v)\frac{(e^{\gamma (u(v))}-1_{\mathbb{A}})}{\gamma (u(v))}=v\frac{%
\log (1_{\mathbb{A}}+\gamma (v))}{\gamma (v)}\cdot \frac{(e^{\gamma
(u(v))}-1_{\mathbb{A}})}{\gamma (u(v))}=v. 
\]

In particular this holds for $1+\gamma (v)\in \mathbb{A}^{-1}$ in the $%
\omega $-homogeneous case, as Cor. 2.2 gives $\gamma (u(tv))=\log (1_{%
\mathbb{A}}+t\gamma (v))$ for $0\leq t||\gamma (v)||<1$ and Theorem S
extends the domain of validity. The final claim follows, since $\gamma $ is
continuous and the set in question is an open neighbourhood of $0.$\hfill $%
\square $

\bigskip

\noindent \textbf{Remark. }In the particular case of $\gamma $ real, i.e.
with values in $\langle 1_{\mathbb{A}}\rangle $, if $1_{\mathbb{A}}+\gamma
(v)>0,$ then $\log (1_{\mathbb{A}}+\gamma (v))$ is real, as may be verified
from either the familiar series or that for $\log (1_{\mathbb{A}}-\gamma
(v)/(1_{\mathbb{A}}+\gamma (v))),$ when $\gamma (v)\geq 1;$ then 
\[
\gamma (u(v))=\gamma (\log (1_{\mathbb{A}}+\gamma (v))\cdot v/\gamma
(v))=\log (1_{\mathbb{A}}+\gamma (v))\cdot \gamma (v)/\gamma (v). 
\]

\noindent \textbf{Remark (Standardized tilting).} Guided by the alignment of 
$u$ and $v$ in the formula $(T$-inv) above, for $\mathbb{A}$ with \textit{%
dense invertibles} (\S 1.2) it is natural to measure the `tilt' of $v=T(u)$
relative to $u$ (scaling included) when $u,v\in \mathbb{A}^{-1}$ by $\theta
:=v^{-1}u\in \mathbb{A}^{-1}$. Then $u=\theta v$ solves $(T)$ for a given $%
v, $ provided $w=w(\theta ):=\gamma (\theta v)\theta ^{-1}$ satisfies the
apparently simpler equation:%
\begin{equation}
e^{\theta w}=1+w.  \tag{$ST_{\QTR{Bbb}{A}}$}
\end{equation}%
(ST\ for `Standardized Tilting'.) Sufficiency of this condition (when
holding for some $\theta $) is proved exactly as in Prop. 5.1 (Converse
part).

It is thematic that $(ST_{\mathbb{A}})$ equates a canonical (affine) Popa
function with a degenerate (exponential) one, and its solubility may perhaps
depend on the geometry of $\mathbb{A}.$

For $\gamma $ real-valued and $\theta \in \mathbb{R}$, both sides of $(ST_{%
\mathbb{A}})$ are real; the resulting formula for $u$ coincides with ($T$%
-inv), since $w=\gamma (v)$ here, so that $\theta =\log (1+\gamma
(v))/\gamma (v),$ as before (assuming $1+\gamma (v)$ has a logarithm). In
general, one wants $\theta $ to induce a tilt aligning $\log (1+w(\theta ))$
with $\gamma (\theta v),$ with $\theta $ pointing in the direction of
norm-increase of the exponential function (equivalently of $\gamma (\theta
v) $ -- see Lemma 5.5 below), to yield a solution for $(ST_{\mathbb{A}}),$
giving%
\[
u=\theta v\log (1+\gamma (\theta v)\theta ^{-1})/\gamma (\theta v). 
\]

\noindent \textbf{Example 5.1. }For $\mathbb{A=C}$, it emerges that simple
collinearity can be effected. With $\theta =1$ $(ST_{\mathbb{A}})$ takes on
its simplest form:%
\begin{equation}
e^{\omega }=1+\omega .  \tag{$ST_{\QTR{Bbb}{C}}$}
\end{equation}%
The example is instructive: while this has no solutions in $\func{Re}(\omega
)<0,$ there exists a sequence (necessarily unbounded, on account of the
Identity Theorem [Gam2,V.7], [Rem, Ch. 8], as in the Great Picard Theorem
[Gam2, XII.2]) of solutions in $\func{Re}(\omega )>0.$ (In the half-space of
unbounded growth of $e^{\func{Re}(\omega )},$ the factor $e^{\func{Im}%
(\omega )}$ can dampen exponential growth to linear in $\omega $ on the
sequence, by taking $\func{Im}(\omega )$ arbitrarily close to $\pi /2$ mod $%
2\pi \mathbb{Z}$, by taking $\func{Re}(\omega )$ suitably large: for the
details see \S 7 (Appendix).) This resonates with a similar behaviour in a
complex Banach algebra when it contains elements $h$ with $||e^{ith}||=1$
for all $t\in \mathbb{R}$, the `hermitian' elements -- see [Pal] for their
particular relevance to $C^{\ast }$-algebras.

The next result, on uni-directional unboundedness, has two exceptions, just
as in $\mathbb{C}$: the unit sphere and the vector subspace $\{a\in \mathbb{A%
}:\lim_{t\rightarrow \infty }e^{\pm at}=0\}$. See [BohK], where the bound $%
||e^{ta}||\leq e^{t\varphi (a)}$ for $t>0,$ arises from the function:%
\[
\varphi (a):=\sup_{t>0}(t^{-1}\log ||e^{ta}||). 
\]%
\bigskip

\noindent \textbf{Lemma 5.5 (One-sided unboundedness). }\textit{Assume }$%
||1_{\mathbb{A}}||=1.$\textit{\ Unless }$||e^{a}||=1,$\textit{\ then as }$%
s\rightarrow \infty $ \textit{only one of the two sets }$E_{\pm }:=$\textit{%
\ }$\{e^{\pm sa}:s>0\}$\textit{\ is unbounded, the other being bounded by }$%
0 $\textit{. In particular, unless }$||e^{\gamma (u)}||=1,$\textit{\ the
range of the map }%
\[
s\mapsto T(su)=u(e^{s\gamma (u)}-1)/\gamma (u), 
\]%
\textit{for }$u$\textit{\ fixed, is unbounded as }$s\rightarrow \infty $ 
\textit{in only one of the directions }$\pm u,$\textit{\ with limit }$%
-u/\gamma (u)$ \textit{in the other.}

\bigskip

\noindent \textbf{Proof. }This comes from a routine induction, based on
noting that if $||e^{-a}||<1,$ then $||e^{a}||>1:$ see \S 7
(Appendix).\hfill\ $\square $

\bigskip

\noindent \textbf{Corollary 5.4.} \textit{For} $\mathbb{A}=\mathbb{R}^{d}$:
(i) $\mathcal{H}=\mathcal{N}$;\newline
(ii)\textit{\ in the context of Theorem 3.1, }$\gamma $\textit{\ is }$%
\mathbb{A}$\textit{-homogeneouse and equation }$(T)$\textit{\ is soluble for}
$v+1_{\mathbb{A}\hfill }\in \exp (\mathbb{A})$\newline
(iii) \textit{Likewise, with} $\mathbb{A}=\mathbb{C}$, \textit{equation }$%
(T) $\textit{\ is soluble for} $v+1_{\hfill }\in \exp (\mathbb{C}).$

\bigskip

\noindent \textbf{Proof. }First we prove (ii) (whence (i) will follow). For $%
I\subseteq \{1,...,d\}$ we write $v|_{I}$ for the projection $e_{I}(v)$ of $%
v $ on the span $\langle \{e_{i}:i\in I\}\rangle $ of the corresponding base
vectors, as in Theorem 3.1. As there, we partition $d$ into parts $I_{i}$ of
cardinality $d_{i}$ (with $d_{i}$ summing to $d)$ for $i=1,...k,$ obtaining,
as equivalent to $(T)$ on $\mathbb{R}^{d},$ the de-coupled simultaneous
system below of $k$ equations $(T)$ on $\mathbb{R}^{d_{i}}$ in terms of $%
\gamma _{i}(u)=\gamma (u)\mathbf{1|}_{I_{i}}$ with corresponding solutions. 
\[
\left. 
\begin{array}{c}
v|_{I_{i}}=u|_{I_{i}}(e^{\gamma _{i}(u|_{I_{i}})}-1_{\mathbb{A}})/\gamma
_{i}(u|_{I_{i}})\mathbf{1|}_{I_{i}} \\ 
u|_{I_{i}}=(v|_{I_{i}}/\gamma _{i}(v))\log (1_{\mathbb{A}}+v|_{I_{i}})%
\end{array}%
\right\} \qquad (i=1,...,k). 
\]%
Re-combining yields a solution in $\mathbb{R}^{d}$ (with $\gamma $ here $%
\mathbb{A}$-homogeneous).

(i) This now follows by partitioning $\{1,...,d\}$ into the part where the
`degenerate' components of $S$ take value $1$ (Th. 3.2) and applying (ii) to
the complementary part. Compare the final assertion of Th. 3.1.

(iii) This follows from the solubility of $(ST_{\mathbb{C}})$ above. \hfill $%
\square $

\bigskip

Even if $\gamma $ is not $\mathbb{A}$-homogeneous, the iteration%
\[
u_{n+1}=u_{n}+v-T(u_{n}) 
\]%
establishes more generally the final conclusion of Prop. 5.1:

\bigskip

\noindent \textbf{Proposition 5.2}. \textit{The range of the tilting map }$T$%
\textit{\ contains a nhd of }$0.$

\bigskip

We begin with a technical

\bigskip

\noindent \textbf{Lemma 5.6. }\textit{For }$S$\textit{\ Fr\'{e}chet
differentiable, take }%
\[
H(u):=(e^{\gamma (u)}-1-\gamma (u))/\gamma (u)=\gamma (u)(\frac{1}{2}+\frac{1%
}{3!}\gamma (u)+...). 
\]%
\textit{Then, for some }$\delta >0$\textit{,}%
\begin{eqnarray*}
||H(u)|| &\leq &||\gamma ||.||u||e^{||\gamma ||.||u||}, \\
||aH(a)-bH(b)|| &\leq &\frac{1}{2}||a-b||\qquad (a,b\in B_{\delta }(0)).
\end{eqnarray*}

\noindent \textbf{Proof. }Below we write $1$ for $1_{\mathbb{A}}.$ By
hypothesis, $S$ is Fr\'{e}chet differentiable, so $||\gamma (u)||\leq
||\gamma ||.||u||.$ By the triangle inequality applied termwise to the
series defining $H$, 
\[
||H(u)||\leq ||\gamma ||.||u||e^{||\gamma ||.||u||}. 
\]%
A similar approach gives for any $u,w$%
\[
||wH^{\prime }(u)||\leq ||w||.||\gamma ||e^{||\gamma ||.||u||}. 
\]%
Indeed 
\[
H^{\prime }(u)h=(\frac{1}{2}+\frac{2}{3!}\gamma (u)+...)\gamma (h). 
\]%
Now%
\begin{eqnarray*}
aH(a)-bH(b) &=&a(H(a)-H(b))+(a-b)H(b) \\
&=&aH^{\prime }(a)(a-b)+o(a-b)+(a-b)H(b).
\end{eqnarray*}%
Thus for $a,b$ small enough, say for $a,b\in B_{\delta }(0),$%
\[
||aH(a)-bH(b)||\leq \frac{1}{2}||a-b||. 
\]%
Indeed, $||H(x)||\leq ||x||.||\gamma ||.e^{||\gamma ||.||x||}<1/3$ and $%
||xH^{\prime }(x)||\leq ||x||.||\gamma ||e^{||\gamma ||}<\frac{1}{3}$
provided $||x||<\min \{1,1/(3||\gamma ||e^{||\gamma ||})\}.$\hfill $\square $

\bigskip

\noindent \textbf{Proof of Proposition 5.2.} We assume $||\gamma ||>0,$
otherwise $T(u)=u$ and the result is immediate. With $H$ as in Lemma 5.6%
\[
T(u)-u=u(e^{\gamma (u)}-1-\gamma (u))/\gamma (u)=uH(u). 
\]%
With $\delta $ as in Lemma 5.6, put%
\[
\eta =\min \{1,\delta /2,\delta /(2||\gamma ||e^{||\gamma ||})\}<\delta /2. 
\]%
Take $v\in B_{\eta }(0)$ and $u_{1}=v,$ so $||u_{1}||<\delta /2.$ Define a
recurrence by%
\[
u_{n+1}=v-u_{n}H(u_{n})=v-(\frac{1}{2}u_{n}\gamma (u_{n})+...). 
\]%
Then $||u_{2}||<\delta $, since by Lemma 5.6, 
\[
||u_{2}-u_{1}||=||vH(v)||\leq ||\gamma ||.||v||^{2}e^{||\gamma ||}\leq
||v||.||\gamma ||e^{||\gamma ||}<\delta /4. 
\]%
Apply Lemma 5.6 inductively, with the inductive hypothesis%
\[
u_{n},u_{n-1}\in B_{\delta }(0)\qquad \text{and}\qquad
||u_{n}-u_{n-1}||<\delta /2^{n}, 
\]%
which holds for $n=2.$ Since $u_{n},u_{n-1}\in B_{\delta }(0),$ by Lemma 5.6%
\begin{eqnarray*}
||u_{n+1}-u_{n}|| &=&||u_{n-1}H(u_{n-1})-u_{n}H(u_{n})||\leq \frac{1}{2}%
||u_{n}-u_{n-1}||<\delta /2^{n}: \\
||u_{n+1}-u_{1}|| &<&(\delta /4)+(\delta /8)+...<\delta /2;
\end{eqnarray*}%
so $u_{n+1},u_{n}\in B_{\delta }(0),$ as $||u_{1}||<\delta /2,$ completing
the induction.

So the sequence $\{u_{n}\}$ is Cauchy. Say $u=\lim u_{n};$ then%
\[
u=v-uH(u)\text{, i.e. }v=T(u). 
\]%
As $v$ was arbitrary, $B_{\eta }(0)$ is in the range of $T.$\hfill $\square $

\bigskip

The next result connects the solubility of ($T$) to whether $\mathcal{N}$ is
a vector subspace and is motivated by Example 5.1.

\bigskip

\noindent \textbf{Proposition 5.3. }\textit{With the assumptions and
notation of Theorem 4.2: }

\noindent (i) \textit{the following hold for }$t\geq 0$\textit{:}%
\begin{eqnarray}
N^{\prime }(0)u &=&(\gamma (u)/(e^{\gamma (u)}-1))N(u(e^{\gamma
(u)}-1)/\gamma (u))\text{ for any }u,  \label{$1$} \\
N(T(tu)) &=& \lambda _{u}(t)N^{\prime }(0)T(u)\text{ for }u\in \mathcal{N}:
\label{$2$} \\
N(tu) &=& tN^{\prime }(0)u\text{ for }u\in \mathcal{H}.  \label{$3$}
\end{eqnarray}%
\textit{So }$N$\textit{\ is linear on }$\mathcal{H}$\textit{\ and }$%
N^{\prime }(0)u=-\rho u$\textit{\ for} $u\in \mathcal{H}.$

\noindent (ii) \textit{Furthermore, for} $u\in \mathcal{N}$,%
\[
\{\lambda _{u}(t)N^{\prime }(0)T(u):t\geq 0\}\subseteq \mathcal{N}. 
\]%
\textit{In particular, for }$u\in \mathcal{N}$ \textit{and any }$t\geq 0,$ 
\textit{provided }$\lambda _{u}(t)$\textit{\ is invertible},\textit{\ }$%
T(tu)\in \mathcal{N}$ \textit{iff }$T(u)\in \mathcal{H}$. \textit{%
Furthermore, }$T(u)\in \mathcal{N}$ \textit{iff} $u\in \mathcal{H}$. \textit{%
So for }$u\in \mathcal{N}$\textit{:}

\[
T(u)\in \mathcal{H}\text{ iff }T(u)\in \mathcal{N}\text{ iff }u\in \mathcal{H%
}. 
\]%
\noindent (iii) $\mathcal{N}$\textit{\ is a vector subspace if, for all
large }$v$\textit{, one of }$\pm v=T(u)$\textit{\ is soluble.}

\bigskip

\noindent \textbf{Proof. }(i)(1) For any $u$ and subject to the L'Hospital
convention, since $N$ is Fr\'{e}chet differentiable, recalling $(Rad)$ of
Theorem 5.2:%
\begin{eqnarray*}
N(u(e^{t\gamma (u)}-1)/\gamma (u)) &=&\lambda _{u}(t)N(u(e^{\gamma
(u)}-1]/\gamma (u))\qquad \text{(i.e. }(Rad)\text{):} \\
N(tu) &=&t\gamma (u)/[e^{\gamma (u)}-1]N(u(e^{\gamma (u)}-1)/\gamma
(u))+o(t), \\
N(tu)/t &=&\gamma (u)/[e^{\gamma (u)}-1]N(u(e^{\gamma (u)}-1)/\gamma
(u))+o(t)/t.
\end{eqnarray*}%
Now passage to the limit $t\rightarrow 0$ yields the claimed formula (as $%
N(0)=0)$.

(2) Differentiating the radiality formula (left to right) with respect to $t$
and using commutativity and post-multiplication yields%
\[
N(u(e^{\gamma (u)}-1)/\gamma (u))\gamma (u)e^{t\gamma (u)}/(e^{\gamma
(u)}-1)=N^{\prime }(u(e^{t\gamma (u)}-1)/\gamma (u))ue^{t\gamma (u)}, 
\]%
\[
N(u(e^{\gamma (u)}-1)/\gamma (u))=N^{\prime }(u(e^{t\gamma (u)}-1)/\gamma
(u))u(e^{\gamma (u)}-1)/\gamma (u). 
\]%
Setting $t=0$ gives%
\[
N(T(u))=N^{\prime }(0)T(u)\quad (=N^{\prime }(u)T(u)\text{ for }u\in 
\mathcal{N}), 
\]%
the last since by $(\ddag ),$ $N^{\prime }(u)=N^{\prime }(0)$ for $u\in 
\mathcal{N}$. Formula (2) now follows from $N(T(tu))=\lambda _{u}(t)N(T(u)).$

(3) Since $N(0)=0,$ for $u\in \mathcal{H}$, 
\[
N^{\prime }(0)u=\lim (N(tu)-N(0))/t=N(u). 
\]%
(ii) Since $N$ maps into $\mathcal{N},$ the first assertion follows from
formula (2), since $\lambda _{u}(t)N^{\prime }(0)T(u)=N(T(tu))\in \mathcal{N}
$, for $u\in \mathcal{N}$. Here $\lambda _{u}(t)\sim t\gamma (u)/[e^{\gamma
(u)}-1],$ so $\mathcal{N}$ is tangentially dense along $\gamma
(u)/[e^{\gamma (u)}-1]N^{\prime }(0)T(u).$

Suppose $\lambda _{u}(t)$ is invertible and $u\in \mathcal{N}$, then $%
T(tu)\in \mathcal{N}$ iff%
\begin{eqnarray*}
-\rho T(tu) &=&N(T(tu)):\text{(Lemma 4.4; use Th. 4.2)} \\
-\rho \lambda _{u}(t)T(u) &=&\lambda _{u}(t)N^{\prime }(0)T(u):\text{(using
(2) for }u\in \mathcal{N}\text{)} \\
-\rho T(u) &=&N^{\prime }(0)T(u)\text{ (cancelling)\quad iff }T(u)\in 
\mathcal{H}.
\end{eqnarray*}%
As $\mu (\gamma (u))$ is invertible, $T(u)=u(e^{\gamma (u)}-1)/\gamma (u)\in 
\mathcal{N}$ iff%
\begin{eqnarray*}
N(u(e^{\gamma (u)}-1)/\gamma (u)) &=&-\rho u(e^{\gamma (u)}-1)/\gamma (u)%
\text{ (Lemma 4.4)} \\
(\gamma (u)/(e^{\gamma (u)}-1))N(u(e^{\gamma (u)}-1)/\gamma (u)) &=&-\rho u%
\text{ (cross multiply)} \\
N^{\prime }(0)u &=&-\rho u\text{ (using (1))\quad iff }u\in \mathcal{H}.
\end{eqnarray*}%
The final claim comes, since $\lambda _{u}(1)=1_{\mathbb{A}},$ and so we may
combine $u\in \mathcal{H}$ iff $T(u)\in \mathcal{N}$ with $T(u)\in \mathcal{N%
}$ iff $T(u)\in \mathcal{H}$.

(iii) As to the last claim, take $v\in \mathcal{N}$. Hence $\pm nv\in 
\mathcal{N}$ for $n\in \mathbb{N},$ as $\mathcal{N}$ is an additive
subgroup. Say $+nv=T(u),$ for some $n\in \mathbb{N}$ and $u.$ Then $%
nv=T(u)\in \mathcal{H}$ since $nv=T(u)\in \mathcal{N}$. Then $v\in \mathcal{H%
},$ as $\mathcal{H}$ is a vector subspace. Thus $\mathcal{N}\subseteq 
\mathcal{H}$ and so $\mathcal{N}=\mathcal{H};$ that is, $\mathcal{N}$ is a
vector subspace.\hfill $\square $

\bigskip

To clarify a first application of Th. 5.1 below (to the case of $\mathbb{C}$%
) we offer

\bigskip

\noindent \textbf{Example 5.2.} If $S(\zeta )=1+a\func{Re}(\zeta )+b\func{Im}%
(\zeta )$ with $a,b$ real and $b\neq 0,$ then $\mathcal{N}:=\{\zeta \in 
\mathbb{C}:S(\zeta )=1\}=\langle b-ai\rangle .$ For, writing $\zeta =x+iy,$ 
\[
ax+by=0\text{ iff }z=x-i(ax)/b=(x/b)[b-ai]. 
\]%
Here $S(1)=1+a,$ so $\rho =a$ and the adjustor for $1+\rho \zeta $ is 
\[
N(\zeta )=-a\zeta +a\func{Re}(\zeta )+b\func{Im}(\zeta )=(b-ai)\func{Im}%
(\zeta ), 
\]%
which is $\mathbb{R}$-linear on $\mathbb{C}$ but not $\mathbb{C}$%
-differentiable.

\bigskip

We now use the adjustor of Theorem 5.1 to characterize the continuous
solutions of $(GS_{\mathbb{C}}).$ This analysis shows that the two kinds of
solution correspond to the distinction between \textit{analyticity} and 
\textit{real analyticity}, for which see [KraP]. Like Baron [Bar], we rely
on Theorem B of \S 4 above (albeit through a corollary), which is key also
to the Gebert proof [Geb]. (The latter source also uses the Wo\l od\'{z}%
ko-Javor characterization, for which see \S 6.1, to infer very elegantly
from Th. B that, unless $S$ is real-valued, $\mathcal{N}$ is discrete, and
so $\mathcal{N}$ is the trivial vector subspace $\{0\}$, since $S(\mathbb{C}%
) $ is connected and $S(\mathbb{C})\mathcal{N}=\mathcal{N}$; so Corollary
5.2 applies here.)

\bigskip

\noindent \textbf{Corollary 5.6 }(cf. [Bar])\textbf{. }\textit{For} $\mathbb{%
A=C}$, \textit{if }$S$\textit{\ solves }$(GS_{\mathbb{C}})$\textit{\ and is }%
$\mathbb{C}$\textit{-differentiable, then for some} $\rho \in \mathbb{C}$%
\begin{equation}
\mathit{\ }S(z)=1+\rho z.  \tag{$Can_{\QTR{Bbb}{C}}$}
\end{equation}%
\textit{If }$S$ \textit{is only continuous, so that }$N:\mathbb{C}%
\rightarrow \mathcal{N}\ $\textit{is Fr\'{e}chet differentiable, then for
some} $a,b\in \mathbb{R}$ \textit{\ }%
\begin{equation}
S(z)=1+a\func{Re}(z)+b\func{Im}(z). 
\tag{$\QTR{rm}{\func{Re}}$-$\QTR{rm}{\func{Im}}$}
\end{equation}

\noindent \textbf{Proof. }The special case of a $\mathbb{C}$-differentiable $%
S$ is covered by Th. 5.1.

We now assume only that $S$ is continuous. As noted in \S 3, the Popa group $%
\mathbb{G}_{S}^{\ast }(\mathbb{R}^{2})$ is a Lie group, so interpreting $%
z=x+iy$ as $(x,y)$ and $S(z)$ as $S(x,y),$ this, and so also its adjustor $N$%
, is differentiable in the usual Euclidean sense. By Prop. 4.2 the closed
subgroup $\mathcal{N}:=\{z:S(z)=1\}$ is a vector subspace of $\mathbb{R}%
^{2}. $ If the subspace $\mathcal{N}$ is two-dimensional, then $\mathcal{N}=%
\mathbb{C}$ and $S\equiv 1.$ Otherwise $\mathcal{N=}\alpha \mathbb{R}$ for
some $\alpha \in \mathbb{C}$. W.l.o.g. we take $\alpha =1$ (otherwise apply
to $\mathbb{C}$ the transformation $z\mapsto \alpha ^{-1}z).$

We now use the fact that $N:\mathbb{C}\rightarrow \mathbb{R}\ $is Fr\'{e}%
chet (= ordinarily) differentiable. Note that if $L:\mathbb{C}\rightarrow 
\mathbb{R}$ is linear, then for some $a_{L},b_{L}\in \mathbb{R}$ 
\[
L(z)=a_{L}\func{Re}(z)+b_{L}\func{Im}(z) 
\]%
(since $L(x+iy)=L(x)+L(iy)$). Fix $z$ and put $L:=N^{\prime }(z)$ and $%
L_{0}:=N^{\prime }(0).$ Then, with $a:=a_{L_{0}}$ and $b:=b_{L_{0}},$ by
Theorem 5.1,%
\[
S(z)^{-1}[a_{L}\func{Re}(S(z)h)+b_{L}\func{Im}(S(z)h)]=a\func{Re}(h)+b\func{%
Im}(h)\qquad (h\in \mathbb{C}). 
\]%
First suppose that $ab\neq 0.$ Then $S(z)^{-1}$ is real (take e.g. $h=1).$
We may cancel through the real and imaginary parts, giving 
\[
\lbrack a_{L}\func{Re}(h)+b_{L}\func{Im}(h)]=a\func{Re}(h)+b\func{Im}%
(h)\qquad (h\in \mathbb{C}). 
\]%
So $a_{L}=a$ and $b_{L}=b.$ Now suppose $a=b=0;$ then%
\[
\lbrack a_{L}\func{Re}(S(z)h)+b_{L}\func{Im}(S(z)h)]=0\qquad (h\in \mathbb{C}%
). 
\]%
Take $h=\overline{S(z)}$ so that $S(z)h>0$ (as $S$ takes invertible values),
and similarly take $h=i\overline{S(z)},$ to deduce that $a_{L}=b_{L}=0$.

So, in either case $a_{L}=a$ and $b_{L}=b,$ and so $N^{\prime }(z)=N^{\prime
}(0)$ for all $z\in \mathbb{C}.$ Then, since $N(0)=0,$ again by Lemma 5.1,%
\begin{eqnarray*}
N(z) &=&L_{0}(z)=a\func{Re}(z)+b\func{Im}(z): \\
S(z) &=&1+\rho z+a\func{Re}(z)+b\func{Im}(z).
\end{eqnarray*}%
But $\rho :=S(1)-1=\rho +a;$ so $a=0,$ and so%
\[
S(z)=1+\rho z+b\func{Im}(z). 
\]%
For $b=0$ we again obtain $(Can_{\mathbb{C}}).$ Substitution into $(GS)$
shows that $\rho =0,$ if $b\neq 0.$ Now the transformation $z\mapsto (u+iv)z$
of $\mathbb{C}$ with $u,v$ real yields%
\begin{equation}
\func{Im}((u+iv)z)=v\func{Re}(z)+u\func{Im}(z),  \tag{$\func{Im}$}
\end{equation}%
since with $x,y$ real%
\[
\func{Im}((u+iv)(x+iy))=\func{Im}((ux-vy)+i(vx+uy)). 
\]%
As $b$ above is real, this yields ($\mathrm{\func{Re}}$-$\mathrm{\func{Im}}$%
) with $bv$ and $ub$ for $a$ and $b$. Thus, we have obtained both ($Can_{%
\mathbb{C}}$) and ($\mathrm{\func{Re}}$-$\mathrm{\func{Im}}$) assuming
continuity of $S$. \hfill $\square $

\bigskip

The next result corresponds to $S(1_{\mathbb{A}})=1_{\mathbb{A}}$ with $N$
arbitrary except for the condition $N(1_{\mathbb{A}})=0.$ Take a linear $%
\sigma :\mathbb{A}\rightarrow \mathbb{R}$ with $\sigma (1_{\mathbb{A}})=0$.
Then for $\nu :\mathbb{A}\rightarrow \mathcal{N}:=\langle e_{i}:i\in
I\rangle $ as below, $\nu (1_{\mathbb{A}})=0:$ indeed, adjoin $e_{0}=1_{%
\mathbb{A}},$ as a further orthogonal idempotent, and then $\sigma (e_{i}1_{%
\mathbb{A}})=0$ all $i.$

Mutatis mutandis, the idempotents $e_{i}$ below may also be interpreted as
orthogonal projections onto one-dimensional vector subspaces: compare $\func{%
Re}$ and $\func{Im}$ in Example 5.2.

\noindent \textbf{Proposition 5.4.} \textit{For }$e_{i}\in \mathbb{A}$%
\textit{\ mutually orthogonal idempotents (i.e. with }$e_{i}e_{j}=\delta
_{ij}e_{i}$\textit{) and} $\sigma :\mathbb{A}\rightarrow \mathbb{R}$ \textit{%
linear, take}%
\[
\nu (x):=\Sigma _{i}\sigma (e_{i}x)e_{i}, 
\]%
\textit{assumed convergent, then} $\nu :\mathbb{A}\rightarrow \langle
e_{i}:i\in I\rangle $ \textit{is linear and }$(GS)$ \textit{has solution}%
\[
S(x):=1_{\mathbb{A}}+\nu (x). 
\]

\noindent \textbf{Proof. }See \S 7 (Appendix). \hfill $\square $

\bigskip

\noindent \textbf{Remark.} Assume that the $e_{i}$ are projections with
one-dimensional ranges, spanning $\mathbb{A},$ with $1_{\mathbb{A}}=\Sigma
_{i}e_{i}.$Then the Proposition above also includes the case $S(x)=1+\rho x$
for some $\rho \in \mathbb{A}$. For, take $\rho _{i}:=e_{i}\rho ,$ $%
x_{i}=e_{i}x,$ and, interpreting these as scalar multipliers (as $e_{i}x=\xi
_{i}e_{i}$ for some scalar $\xi _{i}),$ put 
\[
f_{\rho }(x):=\Sigma _{i}\rho _{i}x_{i}, 
\]%
so that $f_{\rho }:\mathbb{A}\rightarrow \mathbb{R}$ is linear. Then%
\[
\nu (x):=\Sigma _{i}f_{\rho }(e_{i}x)e_{i}=\Sigma _{i}\rho _{i}x_{i}e_{i}. 
\]%
But $\rho =\rho 1_{\mathbb{A}}=\Sigma _{i}\rho e_{i}e_{i}=\Sigma _{i}\rho
_{i}e_{i}$ and so%
\[
\rho x=\Sigma _{i}\rho _{i}e_{i}\Sigma _{j}x_{j}e_{j}=\Sigma _{ij}\rho
_{i}e_{i}x_{i}e_{j}=\Sigma _{i}\rho _{i}x_{i}e_{i}=\nu (x). 
\]

This may be viewed as the totally independent case in that%
\[
\nu _{i}(x)=\nu (x)e_{i}=\rho _{i}x_{i}:\qquad S(x)e_{i}=1+\rho _{i}x_{i}. 
\]

\section{Complements}

\subsection{Wo\l od\'{z}ko-Javor theory}

The following result of Wo\l od\'{z}ko [Wol] was presented in 1968 as a
`construction' yielding, for $F$ a commutative field, all $F$\textit{-valued}
solutions $S$ of $(GS)$ over a vector space $X$ and was cited as a theorem
first in [Jav], also in 1968, and again later in the textbook [AczD, Ch. 19
Th. 5] (attributed there to [Jav], but see the comment in [BriD, Prop. 4]).
The idea, however, may be traced back to [GolS, Th. 4]. We check that this
characterization of solutions of $S$ continues to hold also over $\mathbb{A}$
by interpreting the invertible elements of $F$ there by $\mathbb{A}^{-1}$
here. This yields, mutatis mutandis, a characterization of the restriction $%
S|G_{S}^{\ast }(\mathbb{A)}:\mathbb{A}\rightarrow $ $\mathbb{A}^{-1}$ and a
connection with the Goldie equation, noted in the Remark below.

\bigskip

\noindent \textbf{Theorem 6.1 (Wo\l od\'{z}ko-Javor Theorem).} $S:\mathbb{A}%
\rightarrow $ $\mathbb{A}$ \textit{solves }$(GS)$\textit{\ iff there exist:
an additive subgroup }$\mathcal{N}$\textit{\ of }$\mathbb{A}$\textit{, a
multiplicative subgroup }$\Lambda $\textit{\ of }$\mathbb{A}^{-1}$\textit{,
and a function }$W:\Lambda \rightarrow \mathbb{A}$ \textit{such that:}%
\newline
i) $\Lambda \mathcal{N}=\mathcal{N}$,\newline
ii) $W(\lambda )\in \mathcal{N}$ \textit{iff} $\lambda =1,$\newline
iii) \textit{for all }$\lambda _{1},\lambda _{2}\in \Lambda ,$ \textit{the Wo%
\l od\'{z}ko equation holds:} 
\begin{equation}
W(\lambda _{1}\lambda _{2})=W(\lambda _{1})+\lambda _{1}W(\lambda
_{2})\qquad \func{mod}\mathcal{N},  \tag{$W$}
\end{equation}%
iv) 
\[
S(x)=\left\{ 
\begin{array}{cc}
\lambda , & \text{if }x=W(\lambda )\quad \func{mod}\mathcal{N},\text{ for
some }\lambda \in \Lambda , \\ 
0, & \text{otherwise.}%
\end{array}%
\right. 
\]%
\noindent \textbf{Proof. }See \S 7 (Appendix).\hfill $\square $

\bigskip

\noindent \textbf{Remark. }$K(u):=W(\exp (u))$ ($u\in \mathbb{A}^{-1}$)
converts $(W)$ to a $(GFE)$ variant:%
\[
K(u_{1}+u_{2})=K(u_{1})+e^{u_{1}}K(u_{2})\qquad \func{mod}\mathcal{N}. 
\]

\subsection{Functional equations and probability theory}

We mention briefly some of the probability background to some of the
functional equations we use here. For a monograph treatment of the
intersection of these two areas, see [BalL] and [KagLR].

\noindent \textit{1. The Cauchy functional equation,} $(CFE).$ The Cauchy
functional equation of \S 1.1 is ubiquitous, for example in regular
variation [BinGT, \S 1.1], and in the lack-of-memory property of the
exponential law, and so in the Markov property in continuous time (see e.g.
[GriS, \S 6.9]). In the setting of e.g. renewal theory, survival times X are
in $\mathbb{R}_{+}:=[0,\infty )$, survival probabilities $S(x):=P(X>x)(x\geq
0)$ are monotone, so solutions $S$ to the $(CFE)$ here are exponential (see
e.g. [GriS, 4.14.5]), $e^{-\lambda x}$ ($\lambda \geq 0,$ and $\lambda >0$
in the non-trivial case).

\noindent \textit{2.The Go\l \k{a}b-Schinzel equation} $(GS).$ Replacing$(%
\mathbb{R}_{+},+)$ here by $(G_{\rho }(\mathbb{R}),\circ _{\rho })$ takes $%
(CFE)$ to $(GS)$, with solutions the generalised Pareto laws of EVT
([BinO5], and for higher dimensions, [RooT]). A relative of $(GS)$,%
\[
S(x+\theta (x)y)=S(x)S(y)\quad (\theta (x)=1+cx), 
\]%
appears in the probability literature in [OakD] (cf. [AsaRS, Th. 3.4]) in
connection with characterisation of the Hall-Wellner laws [HalW, Prop. 6].

\noindent \textit{3. (CFE) in higher dimensions. }This is sensitive to
quantifier weakening. See Marshall and Olkin [MarO1]; for further
developments, see [MarO2, \S 4], [MarO3, Ex. 5.1].

\bigskip

\noindent \textbf{Acknowledgement.} It is a great pleasure to thank our
colleague Amol Sasane for the benefits of extensive and encouraging
discussions of Prop. 5.1, enabling its significant improvement and a greater
understanding of the context.

\bigskip

\noindent \textbf{References}

\noindent \lbrack AczD] J. Acz\'{e}l and J. Dhombres, \textsl{Functional
equations in several variables. With applications to mathematics,
information theory and to the natural and social sciences.} Encycl. Math.
App.\textbf{\ 31}, Cambridge Univ. Press, 1989.\newline
\noindent \lbrack AsaRS] M. Asadi, C. R. Rao, D. N. Shanbhag, Some unified
characterization results on generalized Pareto distributions. \textsl{J.
Statist. Plann. Inference} \textbf{93} (2001), no. 1-2, 29--50.\newline
\noindent \lbrack BajK] B. Baj\v{s}anski and J. Karamata, Regularly varying
functions and the principle of equicontinuity, \textsl{Publ. Ramanujan Inst.}
\textbf{1} (1969), 235-246.\newline
\noindent \lbrack BalK] R. Balasubrahmanian and K.-S. Lau, \textsl{%
Functional equations in probability}, Acad. Press, 1991\newline
\noindent \lbrack Bal] A. A. Balkema, \textsl{Monotone transformatio\textsl{n%
} and limit laws.} Mathematical Centre Tracts, No. 45. Math. Centrum,
Amsterdam, 1973.\newline
\noindent \lbrack Bar] K. Baron, On the continuous solutions of the Go\l 
\k{a}b-Schinzel equation. \textsl{Aequationes Math.} \textbf{38} (1989), no.
2-3, 155--162.\newline
\noindent \lbrack Ber] M. S. Berger, \textsl{Nonlinearity and functional
analysis}, Acad. Press, 1977.\newline
\noindent \lbrack BinGT] N. H. Bingham, C. M. Goldie and J. L. Teugels, 
\textsl{Regular variation}, 2nd ed., Cambridge University Press, 1989 (1st
ed. 1987).\newline
\noindent \lbrack BinO1] {N. H. Bingham and A. J. Ostaszewski, }General
regular variation, Popa groups and quantifier weakening. \textsl{J. Math.
Anal. Appl.} \textbf{483} (2020) 123610 (arXiv1901.05996). \newline
\noindent \lbrack BinO2] N. H. Bingham and A. J. Ostaszewski, The
Steinhaus-Weil property: I. subcontinuity and amenability, \textsl{Sarajevo
J. Math}, \textbf{16} (2020), 13-32 (arXiv:1607.00049v3).\newline
\noindent \lbrack BinO3] N. H. Bingham and A. J. Ostaszewski, Multivariate
general regular variation: Popa groups on vector spaces (arXiv: 1910.05816).%
\newline
\noindent \lbrack BinO4] {N. H. Bingham and A. J. Ostaszewski, }Multivariate
Popa groups and the Goldie Equation (arXiv:1910.05817).\newline
\noindent \lbrack BinO5] {N. H. Bingham and A. J. Ostaszewski, }Extremes and
regular variation, to appear in the R. A. Doney Festschrift
(arXiv:2001.05420).\newline
\noindent \lbrack BohL] H. F. Bohnenblust and S. Karlin. Geometrical
properties of the unit sphere of Banach algebras. \textsl{Ann. Math.} (2) 
\textbf{62} (1955), 217--229.\newline
\noindent \lbrack Bou] N. Bourbaki, \textsl{General Topology Part 2},
Addison-Wesley, 1966.\newline
\noindent \lbrack BriD] N. Brillou\"{e}t and J. Dhombres, \'{E}quations
fonctionnelles et recherche de sous-groupes. \textsl{Aequationes Math.} 
\textbf{31} (1986), 253--293.\newline
\noindent \lbrack Brz] J. Brzd\k{e}k, Subgroups of the group Z$_{n}$ and a
generalization of the Go\l \k{a}b-Schinzel functional equation. \textsl{%
Aequationes Math. }\textbf{43} (1992), 59--71.\newline
\noindent \lbrack CabC] F. Cabello S\'{a}nchez and J. M. F. Castillo, Banach
space techniques underpinning a theory for nearly additive mappings, \textsl{%
Dissertationes Math}. \textbf{404} (2002), 73 pp\newline
\noindent \lbrack Chu1] J. Chudziak, Semigroup-valued solutions of the Go\l 
\k{a}b-Schinzel type functional equation. \textsl{Abh. Math. Sem. Univ.
Hamburg,} \textbf{76} (2006), 91-98.\newline
\noindent \lbrack Chu2] J. Chudziak, Semigroup-valued solutions of some
composite equations. \textsl{Aequationes Math.} \textbf{88} (2014), 183--198.%
\newline
\noindent \lbrack Con] J. B. Conway, \textsl{A course in functional analysis}%
, Grad. Texts Math. 96, Springer, 2$^{\text{nd.}}$ ed. 1996 (1$^{\text{st }}$%
ed. 1985).\newline
\noindent \lbrack Dal] H. G. Dales, \textsl{Banach algebras and automatic
continuity.} London Math. Soc. Monographs \textbf{24}, Oxford Univ. Press,
2000.\newline
\noindent \lbrack DalF] H. G. Dales and J. F. Feinstein, Banach function
algebras with dense invertible group. \textsl{Proc. Amer. Math. Soc.} 
\textbf{136} (2008), 1295--1304.\newline
\noindent \lbrack Gam1] T. W. Gamelin, \textsl{Uniform algebras}.
Prentice-Hall, 1969.\newline
\noindent \lbrack Gam2] T. W. Gamelin, \textsl{Complex Analysis}. Springer,
2001.\newline
\noindent \lbrack Geb] H. Gebert, A direct proof of a theorem of K. Baron,\ 
\textsl{Ann. Math. Silesianae} \textbf{9} (1995), 101-103.\newline
\noindent \lbrack GerS] H.-U. Gerber and E. S. W. Shiu, Option pricing by
Esscher transforms (with discussion). \textsl{Trans. Soc. Actuaries} \textbf{%
46} (1994). 99-191.\newline
\noindent \lbrack GolM] C. M. Goldie and A. Mijatovi\'{c} (ed.), \textsl{%
Probability, analysis and number theory} (N. H. Bingham Festschrift), 
\textsl{Adv. Appl. Probab.} \textbf{48A }(2016).\newline
\noindent \lbrack GolS] S. Go\l \k{a}b and A. Schinzel, Sur l'\'{e}quation
fonctionelle $f(x+yf(x))=f(x)f(y)$. \textsl{Publ. Math. Debrecen.} \textbf{6}
(1960), 113-125.\newline
\noindent \lbrack GriS] G. R. Grimmett and D. R. Stirzaker, \textsl{%
Probability and random processes.} 3$^{\text{rd}}$ ed. Oxford University
Press, 2001.\newline
\noindent \lbrack HalW] W. J. Hall and J. A. Wellner, (1981). Mean residual
life, \textsl{Proceedings of the International Symposium on Statistics and
Related Topics} (eds. M. Cs\"{o}rg\H{o} et al.), 169--184, North Holland.%
\newline
\noindent \lbrack Jab1] E. Jab\l o\'{n}ska, Continuous on rays solutions of
an equation of the Go\l \c{a}b-Schinzel type. \textsl{J. Math. Anal. Appl.} 
\textbf{375} (2011), 223--229.\newline
\noindent \lbrack Jab2] E. Jab\l o\'{n}ska, The pexiderized Go\l \c{a}%
b-Schinzel functional equation. \textsl{J. Math. Anal. Appl.} \textbf{381}
(2011), 565-572.\newline
\noindent \lbrack Jab3] E. Jab\l o\'{n}ska, Remarks concerning the
pexiderized Go\l \c{a}b-Schinzel functional equation. \textsl{J. Math. Appl.}
\textbf{35} (2012), 33-38.\newline
\noindent \lbrack Jav] P. Javor, On the general solution of the functional
equation $f(x+yf(x))=f(x)f(y)$. \textsl{Aequationes Math.} \textbf{1}
(1968), 235-238.\newline
\noindent \lbrack KagLR] A. M. Kagan, Yu. V. Linnik and C. R. Rao, \textsl{%
Characterization problems in mathematical statistics.} Wiley, 1973.\newline
\noindent \lbrack Kec] A. S. Kechris, \textsl{Classical Descriptive Set
Theory.} Grad. Texts in Math. \textbf{156}, 1995.\newline
\noindent \lbrack KirRSW] A. Kiriliouk, H. Rootz\'{e}n, J. Segers and J. L.
Wadsworth, Peaks over thresholds modeling with multivariate generalized
Pareto distributions. \textsl{Technometrics} \textbf{61} (2019), no. 1,
123--135.\newline
\noindent \lbrack KraP] S. G. Krantz and H. R. Parks, \textsl{A primer of
real-analytic functions}, 2nd ed., Birkh\"{a}user, 2002 (1st ed. 1992).%
\newline
\noindent \lbrack Lev] T. Levi-Civita, Sulle funzioni che ammettono una
formula d'addizione del tipo $f(x+y)=\sum_{1}^{n}X_{j}(x)Y_{j}(y)$. \textsl{%
Atta Accad. Naz. Lincei Rend.} \textbf{22(5)} (1913), 181-183.\newline
\noindent \lbrack Lyu] Yu. I. Lyubich, \textsl{Functional analysis I: Linear
functional analysis,} Encyclopaedia Math. Sci. \textbf{19}, Springer, 1992.%
\newline
\noindent \lbrack MarO1] A. W. Marshall and I. Olkin. A multivariate
exponential distribution. \textsl{J. Amer. Statist. Assoc. }\textbf{62}
(1967), 30--44.\newline
\noindent \lbrack MarO2] A. W. Marshall and I. Olkin. Functional equations
for multivariate exponential distributions. \textsl{J. Multivariate Anal.} 
\textbf{39} (1991), 209--215.\newline
\noindent \lbrack MarO3] A. W. Marshall and I. Olkin. Multivariate
exponential and geometric distributions with limited memory. \textsl{J.
Mult. Anal.} \textbf{53} (1995), 110--125\newline
\noindent \lbrack MonZ] D. Montgomery and L. Zipppin, \textsl{Topological
transformation groups.} Interscience, 1955.\newline
\noindent \lbrack Mur] A. Mure\'{n}ko,On solutions of the Go\l \c{a}%
b-Schinzel equation. \textsl{Int. J. Math. Math. Sci.} \textbf{27} (2001),
541--546.\newline
\noindent \lbrack OakD] D. Oakes and T. Dasu, A note on residual life. 
\textsl{Biometrika} \textbf{77} (1990), 409--410.\newline
\noindent \lbrack Ost1] A. J. Ostaszewski, Beyond Lebesgue and Baire III:
Steinhaus' Theorem and its descendants, \textsl{Top. Appl.} \textbf{160}
(2013), 1144-1155.\newline
\noindent \lbrack Ost2] A. J. Ostaszewski, Beurling regular variation, Bloom
dichotomy, and the Go\l \k{a}b-Schinzel functional equation, \textsl{Aequat.
Math.} \textbf{89} (2015), 725-744. \newline
\noindent \lbrack Ost3] A. J. Ostaszewski, Homomorphisms from Functional
Equations: The Goldie Equation. \textsl{Aequat. Math. }\textbf{90} (2016),
427-448 (arXiv: 1407.4089).\newline
\noindent \lbrack Ost4] A. J. Ostaszewski, Stable laws and Beurling kernels,
[GolM], 239-248.\newline
\noindent \lbrack Ost5] A. J. Ostaszewski, Homomorphisms from Functional
Equations in Probability, in: \textsl{Developments in Functional Equations
and Related Topics}, ed. J. Brzd\k{e}k et al., Springer Optim. Appl. 124
(2017), 171-213.\newline
\noindent \lbrack Oxt] J. C. Oxtoby: \textsl{Measure and category}, 2nd ed.
Graduate Texts in Math. \textbf{2}, Springer, 1980 (1$^{\text{st}}$ ed.
1971).\newline
\noindent \lbrack Pal] T. W. Palmer, Characterizations of $C^{\ast }$%
-algebras. II. \textsl{Trans. Amer. Math. Soc.} \textbf{148} (1970),
577--588.\newline
\noindent \lbrack PitP] E. J. G. Pitman and J. W. Pitman, A direct approach
to the stable distributions, [GolM], 261-282.\newline
\noindent \lbrack Pop] C. G. Popa, Sur l'\'{e}quation fonctionelle $%
f[x+yf(x)]=f(x)f(y).$ \textsl{Ann. Polon. Math.} \textbf{17} (1965), 193-198.%
\newline
\noindent \lbrack Rem] R. Remmert, \textsl{Theory of complex functions.}
Springer, 1991.\newline
\noindent \lbrack Ric] C. E. Rickart, \textsl{General theory of Banach
algebras.} Van Nostrand, 1960.\newline
\noindent \lbrack RooSW] H. Rootz\'{e}n, J. Segers and J. L. Wadsworth.
Multivariate peaks over thresholds models. \textsl{Extremes} \textbf{21}
(2018), 115--145.\newline
\noindent \lbrack RooT] H. Rootz\'{e}n and N. Tajvidi, Multivariate
generalized Pareto distributions. \textsl{Bernoulli} \textbf{12} (2006),
917--930.\newline
\noindent \lbrack Rud] W. Rudin, \textsl{Functional analysis}, McGraw-Hill,
1991 (1$^{\text{st}}$ ed. 1973).\newline
\noindent \lbrack Sil] B. W. Silverman, \textsl{Density estimation for
statistics and data analysis.} Monogr. Stat. Appl. Prob. \textbf{26},
Chapman \& Hall, 1986.\newline
\noindent \lbrack Stet] H. Stetk\ae r, \textsl{Functional equations on
groups.} World Scientific, 2013.\newline
\noindent \lbrack Tao] T. Tao, \textsl{Hilbert's fifth problem and related
topics}, Grad. Studies Math. \textbf{153}, Amer. Math. Soc., 2014, 338 pp.%
\newline
\noindent \lbrack Wol] S. Wo\l od\'{z}ko, Solution g\'{e}n\'{e}rale de l'%
\'{e}quation fonctionelle. $f(x+yf(x))=f(x)f(y)$ \textsl{Aequationes Math.} 
\textbf{2} (1968-69), 12-29.\newline
\bigskip

\noindent Mathematics Department, Imperial College, London SW7 2AZ;
n.bingham@ic.ac.uk \newline
Mathematics Department, London School of Economics, Houghton Street, London
WC2A 2AE; A.J.Ostaszewski@lse.ac.uk\newpage

\section{Appendix (arXiv only):}

To lighten the burden of the main text above, the four sections here serve
as a repository: first for some further examples, second for an independent
derivation of the radiality formula of Theorem 5.2 in the Euclidean case
that illuminates the role of spectral conditions, third for the fulfilment
of earlier promised referrals to some lengthier, straightforward, albeit
necessary, arguments, and lastly to reproduce Sasane's proof of his Theorem
S.

\subsection{Further Examples}

\noindent \textbf{Example 7.1. }Take $\mathbb{A=R}^{2}$; then the $\mathbb{A}
$-differentiable solutions are of the form%
\[
S(x)=\mathbf{1}+\rho x=(1+\rho _{1}x_{1},1+\rho _{2}x_{2}). 
\]%
Here provided $\rho $ is invertible, $t1_{\mathbb{A}}=S((t-1)/\rho
_{1},(t-1)/\rho _{2}).$

We consider the alternative Fr\'{e}chet differentiable (but not $\mathbb{A}$%
-differentiable) solution given by Th.3.1(ii) (with $d=2$ and trivial
partition)%
\begin{eqnarray*}
S(x) &:&=(1+\sigma _{1}x_{1}+\sigma _{2}x_{2},1+\sigma _{1}x_{1}+\sigma
_{2}x_{2})\in \langle (1,1)\rangle =\langle 1_{\mathbb{A}}\rangle : \\
\gamma (u) &=&DS(0)u=\left[ 
\begin{array}{cc}
\sigma _{1} & \sigma _{2} \\ 
\sigma _{1} & \sigma _{2}%
\end{array}%
\right] u=\left[ 
\begin{array}{c}
\tau (u) \\ 
\tau (u)%
\end{array}%
\right] \text{ for }\tau (u):=\sigma _{1}u_{1}+\sigma _{2}u_{2}, \\
\mathcal{N} &:&=\langle (\sigma _{2},-\sigma _{1})\rangle =\{u:\tau (u)=0\}.
\end{eqnarray*}%
So $e^{\gamma (u)}-1=(e^{\tau (u)}-1,e^{\tau (u)}-1)$ is non-invertible for $%
\tau (u)=0,$ implying $N$ is homogeneous on $\mathcal{N}$ and so linear, as
we shall see. 
\[
\rho :=S(\mathbf{1})-\mathbf{1}=(\sigma _{1}+\sigma _{2},\sigma _{1}+\sigma
_{2})=(\sigma _{1}+\sigma _{2})(1,1). 
\]%
Adjusting $\mathbf{1}+\rho x$ into agreement with $S(x)$ gives%
\[
N(x):=(\sigma _{2}(x_{2}-x_{1}),\sigma
_{1}(x_{1}-x_{2}))=(x_{2}-x_{1})(\sigma _{2},-\sigma _{1})\in \mathcal{N}, 
\]%
which is indeed linear. Then \textrm{ran}$N=\langle (\sigma _{2},-\sigma
_{1})\rangle =\mathcal{N}$, $\ker N=\langle (1,1)\rangle $ and%
\[
\mathcal{N}=\ker (S^{\prime }(0))=\{(x_{1},x_{2}):\sigma _{1}x_{1}+\sigma
_{2}x_{2}=0\}=\langle (\sigma _{2},-\sigma _{1})\rangle . 
\]%
Tilting here is :%
\[
T(u)=u(e^{\gamma (u)}-1)/\gamma (u)=u\left( \frac{e^{(\sigma
_{1}u_{1}+\sigma _{2}u_{2})}-1}{\sigma _{1}u_{1}+\sigma _{2}u_{2}},\frac{%
e^{(\sigma _{1}u_{1}+\sigma _{2}u_{2})}-1}{\sigma _{1}u_{1}+\sigma _{2}u_{2}}%
\right) \text{ or }u\text{ if }\tau (u)=0. 
\]

\noindent \textbf{Example 7.2. }As in Example 7.1, take $\mathbb{A=R}^{3}$;
then the $\mathbb{A}$-differentiable solutions are of the form%
\[
S(x)=1+\rho x=(1+\rho _{1}x_{1},1+\rho _{2}x_{2},1+\rho _{3}x_{3}). 
\]%
We consider the alternative (Fr\'{e}chet differentiable) solution given by
Th. 3.3(iv):

Here, provided $(\sigma _{1},\sigma _{2})\neq 0$ and $\sigma _{3}\neq 0$, if 
$(x_{1},x_{2})$ solves $\sigma _{1}x_{1}+\sigma _{2}x_{2}=(t-1)$, then $t1_{%
\mathbb{A}}=S(x_{1},x_{2},(t-1)/\sigma _{3}),$ in which case%
\[
(\mathrm{ran}S)u\supseteq \langle 1_{\mathbb{A}}\rangle u=\langle u\rangle . 
\]%
\begin{eqnarray*}
S(x) &:&=(1+\sigma _{1}x_{1}+\sigma _{2}x_{2},1+\sigma _{1}x_{1}+\sigma
_{2}x_{2},1+\sigma _{3}x_{3}): \\
\gamma (u) &=&DS(0)u=\left[ 
\begin{array}{ccc}
\sigma _{1} & \sigma _{2} & 0 \\ 
\sigma _{1} & \sigma _{2} & 0 \\ 
0 & 0 & \sigma _{3}%
\end{array}%
\right] u=\left[ 
\begin{array}{c}
\sigma _{1}u_{1}+\sigma _{2}u_{2} \\ 
\sigma _{1}u_{1}+\sigma _{2}u_{2} \\ 
\sigma _{3}u_{3}%
\end{array}%
\right] , \\
(e^{\gamma (u)}-1)/\gamma (u) &=&\left( (e^{\tau (u)}-1)/\tau (u),(e^{\tau
(u)}-1)/\tau (u),\frac{e^{\sigma _{3}u_{3}}-1}{\sigma _{3}u_{3}}\right) , \\
\mathcal{N} &:&=\{(x_{1},x_{2},0):\sigma _{1}x_{1}+\sigma
_{2}x_{2}=0\}=\langle (\sigma _{2},-\sigma _{1},0)\rangle , \\
\rho &:&=S(\mathbf{1})-\mathbf{1}=(\sigma _{1}+\sigma _{2},\sigma
_{1}+\sigma _{2},\sigma _{3}).
\end{eqnarray*}%
Adjusting $\mathbf{1}+\rho x$ into agreement with $S(x)$ gives%
\[
N(x):=(\sigma _{2}(x_{2}-x_{1}),\sigma
_{1}(x_{1}-x_{2}),0)=(x_{2}-x_{1})(\sigma _{2},-\sigma _{1},0)\in \mathcal{N}%
, 
\]%
which is linear; then \textrm{ran}$N=\mathcal{N}$, $\ker N=\langle
(1,1,0),(0,0,1)\rangle .$

\bigskip

\noindent \textbf{Example 7.3. }(Non-linear adjustor for a degenerate $S$).
In the previous examples the solution function $S$ had an underlying linear
form. This is not so for $\mathbb{A=R}^{2}=\mathbb{G}_{S}^{\ast }(\mathbb{R}%
^{2}\mathbb{)}$, when $S(x_{1},x_{2}):=(1,e^{x_{1}}).$ Note that $%
e(x)=(0,e^{x_{1}}-1-x_{1})=(0,\frac{1}{2}x_{1}^{2}+...)$ and $e(x\circ
x)=(0,2x_{1}^{2}+...).$ Here%
\[
e(x)/||x||^{2}=\frac{1}{2}\frac{x_{1}^{2}}{x_{1}^{2}+x_{2}^{2}}\rightarrow 
\text{depends on }x_{2}/x_{1} 
\]%
with variable limit as $x\rightarrow 0.$%
\[
S(x)(0,u)=(0,e^{x_{1}}u):\qquad (\mathrm{ran}S)(0,u)=\langle (0,u)\rangle . 
\]%
\begin{eqnarray*}
S(x_{1},x_{2}) &:&=(1,e^{x_{1}})=\mathbf{1+(}0,e^{x_{1}}-1): \\
\mathcal{N} &:&=\{x:S(x)=\mathbf{1}\}=\{x:x_{1}=0\}=\langle e_{2}\rangle , \\
\gamma (u) &=&DS(0)u=\left[ 
\begin{array}{cc}
0 & 0 \\ 
1 & 0%
\end{array}%
\right] u=\left[ 
\begin{array}{c}
0 \\ 
u_{1}%
\end{array}%
\right] =0\text{ iff }u_{1}=0. \\
\rho &=&S(\mathbf{1})-\mathbf{1}=(0,e-1):\qquad \rho x=(0,(e-1)x_{2}): \\
N(x_{1},x_{2}) &:&=S(x_{1},x_{2})-\mathbf{1}-\rho
x=(0,(e^{x_{1}}-1)-(e-1)x_{2})\in \mathcal{N}.
\end{eqnarray*}%
Applying the L'Hospital convention:%
\begin{eqnarray*}
\mathrm{spec}(\gamma (u)) &=&\{0,u_{1}\}:\qquad \lambda _{u}(t)=(t,\frac{%
e^{tu_{1}}-1}{e^{u_{1}}-1}),\text{ or }(t,t)\text{ if }u_{1}=0: \\
N(te_{2}) &=&t\mathbf{1}N(e_{2})=tN(e_{2})\qquad \text{ (for }u_{1}=0);
\end{eqnarray*}%
but for $u_{1}\neq 0:$%
\begin{eqnarray*}
N((t,e^{t}-1)e_{1}) &=&(t,\frac{e^{t}-1}{e-1})N((1,e^{1}-1)e_{1}): \\
N(te_{1}) &=&(t,\frac{e^{t}-1}{e-1})N(e_{1}): \\
N(t(1,0)) &=&(0,e^{t}-1)=\frac{e^{t}-1}{e-1}(0,e-1)=\lambda (t)N(1,0).
\end{eqnarray*}%
Here 
\[
S(x)(0,u)=(0,e^{x_{1}}u):\qquad (\mathrm{ran}S)(0,u)=\langle (0,u)\rangle . 
\]%
and tilting is given by%
\[
T(u):=u(e^{\gamma (u)}-1)/\gamma (u)=u(1,\frac{e^{u_{1}}-1}{u_{1}})\text{ or 
}u\text{ if }u_{1}=0. 
\]

\noindent \textbf{Example 7.4.} $S(x)=(1,1,e^{x_{1}+x_{2}})$. Here $\mathcal{%
N}=\{x:x_{1}+x_{2}=0\},$ which is 2-dimensional. As in Cor. 4.4:%
\begin{eqnarray*}
S(x)(u_{1},u_{2},u_{3}) &=&(u_{1},u_{2},e^{x_{1}+x_{2}}u_{3}):\qquad u\in (%
\mathrm{ran}S)u: \\
\qquad (-u)+(\mathrm{ran}S)u &\supseteq &\{(0,0,(e^{t}-1)u_{3}):t\in \mathbb{%
R}\};
\end{eqnarray*}%
\begin{eqnarray*}
(\mathrm{ran}S)\mathcal{N} &=&\mathcal{N}:\qquad \mathcal{N}\supseteq
\{(0,0,\lambda u_{3}):\lambda >-1\}\text{ for }u\in \mathcal{N}\text{ and so 
}\mathcal{H}\neq \emptyset . \\
\{u &:&(\mathrm{ran}S)u=\{u\}\}=\{u:u_{3}=0\}\text{ orthogonal to }(0,0,1).
\end{eqnarray*}%
\begin{eqnarray*}
DS(0)u &=&\left[ 
\begin{array}{ccc}
0 & 0 & 0 \\ 
0 & 0 & 0 \\ 
1 & 1 & 0%
\end{array}%
\right] u=\left[ 
\begin{array}{c}
0 \\ 
0 \\ 
u_{1}+u_{2}%
\end{array}%
\right] =0\text{ iff }u_{1}+u_{2}=0, \\
\lambda _{u}(t) &=&\left( t,t,\frac{e^{t(u_{1}+u_{2})}-1}{e^{(u_{1}+u_{2})}-1%
}\right) \text{ or }(t,t,t)\text{ if }u_{1}+u_{2}=0.
\end{eqnarray*}%
So $\gamma (e_{1}-e_{2})=\gamma (e_{1})-\gamma (e_{2})=0$ and $\gamma
(e_{3})=0,$ confirming that $\mathcal{N=\langle }e_{1}-e_{2},e_{3}\mathcal{%
\rangle }$.

So tilting is given by%
\[
T(u)=u(e^{\gamma (u)}-1)/\gamma (u)=u(1,1,\frac{e^{(u_{1}+u_{2})}-1}{%
u_{1}+u_{2}})\text{ or }u\text{ if }u_{1}+u_{2}=0. 
\]

We may check this directly: 
\begin{eqnarray*}
\rho &=&S(\mathbf{1})-\mathbf{1}=(0,0,e^{2}-1):\qquad \rho x=(e^{2}-1)x_{3}:
\\
N(x) &=&(0,0,(e^{x_{1}+x_{2}}-1)-(e^{2}-1)x_{3}).
\end{eqnarray*}%
Here%
\begin{eqnarray*}
N(tx) &=&(0,0,(e^{tx_{1}+tx_{2}}-1)-(e^{2}-1)tx_{3}), \\
N(tx_{1},-tx_{1},tx_{3}) &=&(0,0,-(e^{2}-1)tx_{3})=tN(x_{1},-x_{1},x_{3}): \\
N(tx_{1},-tx_{1},0) &=&tN(x_{1},-x_{1},0)\text{ and }N(te_{3})=tN(e_{3}),
\end{eqnarray*}%
yielding the two directions along which $N$ is homogenous:%
\[
N(t(1,-1,0)=tN(1,-1,0),\text{ and }N(te_{3})=tN(e_{3}). 
\]

\noindent \textbf{Remark. }Along the natural base directions (with $\lambda $
subscript $i$ for $e_{i}$):%
\[
\lambda _{1}(t)=\lambda _{2}(t)=(t,t,(e^{t}-1)/(e-1))\qquad \lambda
_{3}(t)=(t,t,t). 
\]%
So%
\begin{eqnarray*}
N(u(e^{t\gamma (u)}-1)/\gamma (u)) &=&\lambda _{u}(t)N(u(e^{\gamma
(u)}-1]/\gamma (u)): \\
N(e_{i}(t,t,(e^{t}-1)) &=&(t,t,(e^{t}-1)/(e-1))N(e_{i}(1,1,(e-1)):\qquad
(i=1,2) \\
N(te_{i}) &=&(t,t,(e^{t}-1)/(e-1))N(e_{i})\qquad (i=1,2), \\
N(te_{3}) &=&t(1,1,1)N(e_{3})=tN(e_{3}).
\end{eqnarray*}

\subsection{Radiality formula in $\mathbb{R}^{d}$}

We give a direct proof of Theorem 4.2 when $\mathbb{A}=\mathbb{R}^{d}$. We
begin with two Lemmas. In $\mathbb{R}^{d}$ we will write $x_{i}$ or $(x)_{i}$
for the $i$-th component of $x;$ of course $x_{i}=xe_{i}$ where the $e_{i}$
are the natural base vectors. We write $\gamma (u)=S_{u}^{\prime }(0).$
Recall that the \textit{spectrum} of $x$ [Rud, 10.10] is defined by%
\[
\mathrm{spec}(x):=\{\lambda :(\lambda 1_{\mathbb{A}}-x)\notin \mathbb{A}%
^{-1}\}. 
\]%
So in $\mathbb{A=R}^{d},$ $\mathrm{spec}(\gamma (u))$ comprises the
components $\gamma (u)_{i}.$

\bigskip

\noindent \textbf{Lemma 7.1. }\textit{In} $\mathbb{R}^{d}$, $e^{\gamma
(u)}-1_{\mathbb{A}}$ \textit{is invertible iff} $0\notin \mathrm{spec}%
(\gamma (u)).$

\bigskip

\noindent \textbf{Proof.} $\gamma (u)$ has real components, so $\mathrm{spec}%
(\gamma (d))\subseteq $ $\mathbb{R}$. So by [Rud, 10.28], $e^{\gamma
(u)}-1_{A}$ is invertible iff $e^{\lambda }-1\neq 0$ for all $\lambda \in 
\mathrm{spec}(\gamma (d))\subseteq $ $\mathbb{R}$ iff $0\notin \mathrm{spec}%
(\gamma (d)).$\hfill $\square $

\bigskip\ 

\noindent \textbf{Lemma 7.2.} \textit{In} $\mathbb{R}^{d}$, \textit{if }$%
S_{u}^{\prime }(0)_{i}=0$\textit{\ for some }$i$ \textit{and }$u=e_{j}$%
\textit{, then }$S(te_{j})_{i}=1$\textit{\ for all }$t\in \mathbb{R}$, 
\textit{i.e.} $\langle e_{j}\rangle \subseteq \{x:S(x)_{i}=1\}.$ \textit{In
particular, if }$S_{u}^{\prime }(0)_{i}=0,$\textit{\ for }$u=e_{j}$\textit{\
and for all }$i,$ \textit{then }$S(te_{j})=1$\textit{\ for all }$t\in 
\mathbb{R}$, \textit{i.e.} $\langle e_{j}\rangle \subseteq \mathcal{N}%
=\{x:S(x)=1\}.$\textit{\ }

\bigskip

\noindent \textbf{Proof.} Fix $i.$ Working in direction $e_{j}$, $(GS)$
implies that for $s,t\in \mathbb{R}$%
\[
S(se_{j}+tS(se_{j})_{j}e_{j})_{i}=(S(se_{j}))_{i}(S(te_{j}))_{i}. 
\]%
This is a pexiderized form of the real-valued univariate version of $(GS)$
(see e.g. [Jab]). Since $t\mapsto S(te_{j})_{j}$ is continuous for all $j,$
it follows that, for some $\sigma _{ij}\in \mathbb{R}$, 
\[
S(te_{j})_{j}=S(te_{j})_{i}=1+t\sigma _{ij}\qquad (t\in \mathbb{R}). 
\]%
Then for $u=e_{j}$%
\[
S(tu)_{i}-S(0)_{i}=t\sigma _{ij}:\qquad \sigma _{ij}=S_{u}^{\prime
}(0)_{i}=0. 
\]%
So $S(te_{j})_{j}=S(te_{j})_{i}=1+t\sigma _{ij}=1.$ \hfill $\square $

\bigskip

\noindent \textbf{Theorem 7.1}$_{Euclidean}$. \textit{For }$S$\textit{\ Fr%
\'{e}chet differentiable, the solution }$(N,S)$\textit{\ to the Goldie
equation in} $\mathbb{R}^{d},$%
\[
N(x+S(x)y)=N(x)+S(x)N(y), 
\]%
\textit{satisfies, for }$u$\ \textit{any natural base vector,} \textit{the
radiality formulae below (with the L'Hospital convention); that is, with }$%
\gamma (u)=S_{u}^{\prime }(0)$\textit{\ and}%
\[
\lambda _{u}(t):=[e^{t\gamma (u)}-1]/[e^{\gamma (u)}-1]:\qquad \lambda
_{u}(t)_{i}=[e^{t\gamma (u)_{i}}-1]/[e^{\gamma (u)_{i}}-1]\text{,} 
\]%
\[
N(u(e^{t\gamma (u)}-1)/\gamma (u))=\lambda _{u}(t)N(u(e^{\gamma
(u)}-1)/\gamma (u)). 
\]%
\textit{In particular, if }$\gamma (u)_{i}=0$\textit{\ for all }$i,$\textit{%
\ then}%
\[
N(tu)=t\mathbf{1}N(u)=tN(u). 
\]

\noindent \textbf{Proof.} Referring again to the polynomials $\wp
_{n}(x)=1+x+...+x^{n-1}$ and rational polynomials $[\wp _{m}/\wp _{n}](x),$
[BinO4] Theorem 2 gives with $N$ here for $K$ there that for any $u$%
\[
N(\wp _{m}(S(u/n))u/n)=[\wp _{m}/\wp _{n}](S(u/n))N(\wp _{n}(S(u/n))u/n). 
\]%
By Gateaux differentiability, write 
\[
S(u/n)-S(0)=S_{u}^{\prime }(0)(1/n)+\varepsilon _{n}(u), 
\]%
with $n\varepsilon _{n}(u)\rightarrow 0.$ Fix $u=e_{j}.$ We work with odd
integers $n$ below and take $\gamma (u):=S_{u}^{\prime }(0)$ as above and $%
\xi _{n}:=S(u/n).$

Case 1. $\gamma (u)_{i}\neq 0$ for all $i.$ Provided $(\xi _{n})_{i}\neq 1$
for infinitely many (odd) $n:$%
\begin{eqnarray*}
\lbrack \wp _{m}/\wp _{n}](\xi _{n})_{i}
&=&[(S(u/n)_{i})^{m}-1]/[(S(u/n)_{i})^{n}-1] \\
&=&[(1+S_{u}^{\prime }(0)_{i}/n+\varepsilon
_{n}(u))^{nt(n)}-1]/[(1+S_{u}^{\prime }(0)_{i}/n+\varepsilon _{n}(u))^{n}-1]
\\
&\rightarrow &[e^{t\gamma (u)_{i}}-1]/[e^{\gamma (u)_{i}}-1],
\end{eqnarray*}%
as $n\rightarrow \infty $ with $m(n)/n\rightarrow t\in \mathbb{R}.$ If for
some $i,$ $(\xi _{n})_{i}=1$ for all large (odd) $n$, then%
\begin{eqnarray*}
0 &=&S(u/n)_{i}-S(0)_{i}=S_{u}^{\prime }(0)_{i}(1/n)+\varepsilon _{n}(u): \\
S_{u}^{\prime }(0)_{i} &=&-n\varepsilon _{n}(u)\rightarrow 0:\qquad \gamma
(u)_{i}=0,
\end{eqnarray*}%
a contradiction to this case.

Case 2. For some $i,S_{u}^{\prime }(0)_{i}=0.$ Here $u=e_{j},$ so by Lemma 2 
$(\xi _{n})_{i}=S(e_{j}/n)_{i}=1$ for all $n;$ then, again as $%
m/n\rightarrow t$ and as $\wp _{n}(1)=n,$ 
\[
\lbrack \wp _{m}/\wp _{n}](\xi _{n})_{i}=\frac{m(n)}{n}\rightarrow t. 
\]%
So we may interpret the earlier displayed formula using the L'Hospital
convention. \hfill $\square $

\subsection{Expanded arguments}

Below we make good on the promise to set out certain routine arguments.

\subsubsection{Proof of Prop 1.1}

(i) $\mathbb{G}_{S}^{\ast }(\mathbb{A})$ is closed under the operation $%
\circ _{S}$ since 
\[
S(a\circ _{S}b)=S(a)S(b), 
\]%
so that $S(a)S(b)$ is invertible when $S(a)$ and $S(b)$ are invertible. Note
that for $a\in \mathbb{G}_{S}(\mathbb{A})$%
\[
a=a\circ _{S}b=a+S(a)b\text{ iff }b=0. 
\]%
Also for $a\in \mathbb{G}_{S},$ as $S(a)$ is invertible,%
\[
S(a)=S(a\circ _{S}0)=S(a)S(0):\quad S(0)=1_{\mathbb{A}}. 
\]%
Thus the neutral element for $\circ _{S},$ i.e. $0$ is in $\mathbb{G}%
_{S}^{\ast }$. The $\circ _{S}$-inverse of $a$ is $b=-aS(a)^{-1}$. As in
[Jav], the operation is associative: 
\begin{eqnarray*}
(a\circ _{S}b)\circ _{S}c &=&(a+S(a)b)+S(a+S(a)b)c=a+S(a)b+S(a)S(b)c, \\
a\circ _{S}(b\circ _{S}c) &=&a+S(a)[b+S(b)c]=a+S(a)b+S(a)S(b)c.
\end{eqnarray*}

As the elements of $\mathbb{A}_{1}$ are invertible, $\mathbb{G}_{S}:=S^{-1}(%
\mathbb{A}_{1})\subseteq \mathbb{G}_{S}^{\ast }$. As $S(0)=1_{\mathbb{A}},$ $%
0\in \mathbb{G}_{S},$ and furthermore for $a,b\in \mathbb{G}_{S}:$ as $%
S(a),S(b)\in \mathbb{A}_{1},$%
\[
S(a\circ _{S}b)=S(a)S(b)\in \mathbb{A}_{1} 
\]%
(as $\mathbb{A}_{1}$ is a multiplicative group), so $a\circ _{S}b\in \mathbb{%
G}_{S},$ and also $a_{S}^{-1}\in \mathbb{G}_{S}$ because 
\[
S(a_{S}^{-1})=S(a)^{-1}\in \mathbb{A}_{1}, 
\]%
since%
\[
1_{\mathbb{A}}=S(0)=S(a\circ _{S}a_{S}^{-1})=S(a)S(a_{S}^{-1}). 
\]%
So $\mathbb{G}_{S}$ is a subgroup of $\mathbb{G}_{S}^{\ast }.$

(ii) Evidently, for $a,b\in \mathcal{N}$ 
\[
a+b=a+S(a)b=a\circ _{S}b. 
\]%
So if $a,b\in \mathcal{N},$ then $a+b\in \mathcal{N}$ and $-a\in \mathcal{N}%
, $ since%
\begin{eqnarray*}
S(a+b) &=&S(a+S(a)b)=S(a)S(b)=1_{\mathbb{A}}, \\
1_{\mathbb{A}} &=&S(0)=S(a-a)=S(a-S(a)a)=S(a)S(-a)=S(-a).
\end{eqnarray*}%
Also, for $a\in \mathcal{N},$ $-a=-aS(-a)=-aS(a)^{-1}$ is the $\circ _{S}$%
-inverse of $a,$ so $\mathcal{N}$ is both an additive subgroup of $\mathbb{A}
$ and a subgroup of $\mathbb{G}_{S}^{\ast }.$

\subsubsection{\textbf{Proof of Th. 2.2}}

Below we use the notation $a\circ b:=a+bS(a)$ on $\mathbb{G}_{S}:=(a\in 
\mathbb{A}:S(a)\in \mathbb{A}^{-1}\}.$ Differentiation of $(GS)$ w.r.t. to $%
b $ gives%
\[
S^{\prime }(a+bS(a))S(a)=S(a)S^{\prime }(b); 
\]%
now take $a=b_{S}^{-1}$ so $S(a)^{-1}=S(b)$ to obtain the \textit{similarity
relations}:%
\begin{eqnarray*}
S(a)^{-1}S^{\prime }(0)S(a) &=&S^{\prime }(b): \\
S(b)\gamma S(b)^{-1} &=&S^{\prime }(b).
\end{eqnarray*}%
Below we repeatedly write%
\[
S(a+h)=S(a)+S^{\prime }(a)h+o(h). 
\]%
We consider $(GS)$ after $a$ and $b$ are equally incremented by $h:$%
\begin{eqnarray*}
&&S((a+h)+(b+h)(S(a)+S^{\prime }(a)h+o(h)) \\
&=&S((a+bS(a))+h+bS^{\prime }(a)h+hS(a)+o(h)) \\
&=&S(a\circ b)+S^{\prime }(a\circ b)[h+bS^{\prime }(a)h+S(a)h+o(h)]+o(h) \\
&=&(S(a)+S^{\prime }(a)h+o(h))(S(b)+S^{\prime }(b)h+o(h)) \\
&=&S(a)S(b)+S(a)S^{\prime }(b)h+S(b)S^{\prime }(a)h+o(h).
\end{eqnarray*}%
Comparison of the two sided gives to within $o(h)$ 
\[
S^{\prime }(a\circ b)[1+bS^{\prime }(a)+S(a)]=S(a)S^{\prime
}(b)+S(b)S^{\prime }(a). 
\]%
Applying the similarity relations gives%
\begin{eqnarray*}
&&S(a\circ b)\gamma S(a\circ b)^{-1}[1+bS^{\prime }(a)+S(a)] \\
&=&S(a)S(b)\gamma S(b)^{-1}+S(b)S(a)\gamma S(a)^{-1}.
\end{eqnarray*}%
Cancelling $S(a)S(b)$ on the left-hand side$:$%
\[
\gamma S(a)^{-1}S(b)^{-1}[1_{\mathbb{A}}+bS(a)\gamma S(a)^{-1}+S(a)]=\gamma
S(b)^{-1}+\gamma S(a)^{-1}. 
\]%
Absorbing $S(a)^{-1}$ on the left-hand side:%
\[
\gamma S(b)^{-1}[S(a)^{-1}+b\gamma S(a)^{-1}+1_{\mathbb{A}}]=\gamma
S(b)^{-1}+\gamma S(a)^{-1}. 
\]%
Cancelling the term $\gamma S(b)^{-1}$ appearing on each side:%
\[
\gamma \lbrack S(b)^{-1}S(a)^{-1}+S(b)^{-1}b\gamma S(a)^{-1}]=\gamma
S(a)^{-1}. 
\]%
Cancelling $S(a)^{-1}$ on the right on both sides gives:%
\[
\gamma \lbrack S(b_{s}^{-1})-b_{s}^{-1}\gamma ]=\gamma , 
\]%
with $b_{s}^{-1}$ the $\circ _{S}$-inverse of $b.$ Put $c=b_{s}^{-1};$ then,
on rearranging,%
\begin{equation}
\gamma (S(c)h)=\gamma (h)+\gamma (c\gamma h).  \tag{$\ast ^{\prime }$}
\end{equation}%
All the steps above are reversible. We now improve on the last equation.

Differentiating $(\ast ^{\prime })$ with respect to $c$ in direction $k:$%
\[
\gamma (S^{\prime }(c)(k)h)=\gamma (k\gamma h). 
\]%
Using the similarity relations gives%
\[
\gamma (S(c)\gamma S(c)^{-1}(k)h)=\gamma (k\gamma h), 
\]%
which for $a=S(c)$ yields the claim $(\ast \ast \ast )$ and for $c=0$ yields
the claim $(\ast \ast ).$

Writing $S(c)k$ for $k$ in the last equation gives 
\begin{equation}
\gamma (S(c)\gamma (k)h)=\gamma (S(c)k\gamma h).  \tag{$\sharp $}
\end{equation}%
It now follows from $(\ast ^{\prime })$ and $(\ast \ast )\ $that%
\[
\gamma (S(c)h)=\gamma (h)+\gamma (c\gamma h)=\gamma (h)+\gamma (h\gamma
(c))=\gamma (h(1+\gamma (c)). 
\]%
That is, $(\ast )\ $holds. Evidently $(\ast )$ and $(\ast \ast )$ yield $%
(\ast ^{\prime }),$ and so the conjuction of $(\ast )$ and $(\ast \ast )$
yields $(GS).$

It is immediate from $(\ast )$ that%
\[
\gamma (S(c))=\gamma (1_{\mathbb{A}}+\gamma (c)). 
\]%
Furthermore, $(\ast \ast )$ with $k=1_{\mathbb{A}}$ yields%
\[
\gamma (\gamma (h))=\gamma (\gamma (1_{\mathbb{A}})h). 
\]

Differentiating $(\natural )$ with respect to $c$ in direction $u$ and
setting $c=0$ gives%
\[
\gamma (\gamma (u)\gamma (k)h)=\gamma (\gamma (u)k\gamma (h)). 
\]%
Taking $k=1_{\mathbb{A}}$ yields%
\begin{eqnarray*}
\gamma (\gamma (1_{\mathbb{A}})\gamma (u)h) &=&\gamma (\gamma (u)\gamma (h)):
\\
\gamma (\gamma (\gamma (u)h) &=&\gamma (\gamma (u)\gamma (h)).\qquad \qquad
\qquad \square
\end{eqnarray*}

\subsubsection{Proof of Prop. 2.2.}

Assume the weaker homogeneity property and suppose for some $n$ that%
\[
\gamma (x\gamma (x)^{n})=\gamma (x)^{n+1}, 
\]%
which is valid for $n=0;$ the inductive step is provided by%
\[
\gamma (x\gamma (x)^{n+1})=\gamma ((x\gamma (x)^{n}).\gamma (x))=\gamma
(x\gamma (x)^{n}).\gamma (x)=\gamma (x)^{n+1}\gamma (x). 
\]%
Conversely, the case of $(\times )$ for $k$ implies for $y=x\gamma (x)^{k}$
that%
\begin{eqnarray*}
\gamma (y\gamma (x) &=&\gamma (x\gamma (x)^{k}\gamma (x))=\gamma (x\gamma
(x)^{k+1})=\gamma (x)^{k+1}\gamma (x) \\
&=&\gamma (x\gamma (x)^{k})\gamma (x)=\gamma (y)\gamma (x).\qquad \qquad
\qquad \qquad \square
\end{eqnarray*}

\subsubsection{Proof of Proposition 2.3.}

We first prove the case $k=1$ of $(\times ),$ which case also establishes
the converse. By linearity of $\gamma ,$ it is enough to prove this identity
for any $u$ of norm $1.$ Fix $u$ of norm 1, and take $x=tu$ with $%
t=||x||\rightarrow 0$. Substitution into $GS$ leads to%
\begin{equation}
\gamma (x\gamma (x))-\gamma (x)^{2}+\gamma (xe(x))-2e(x)\gamma
(x)=2e(x)-e(x\circ x).  \tag{\dag }
\end{equation}%
Since $e(x)=o(x),$%
\[
||\gamma (xe(x))||/||x||^{2}=||\gamma (u\cdot e(x)/||x||)||\leq ||\gamma
||.||u||.(||e(x)/||x||)\rightarrow 0. 
\]%
Furthermore, 
\[
||x\circ x||/||x||=||2u+tu\gamma (u)+ue(tu)||\rightarrow 2, 
\]%
and so RHS\ of (\dag ) gives%
\[
\lbrack 2e(x)-e(x\circ x)]/||x||^{2}=2e(x)/||x||^{2}-\frac{||x\circ x||^{2}}{%
||x||^{2}}\cdot e(x\circ x)/||x\circ x||^{2}\rightarrow 0. 
\]%
Hence 
\[
\gamma (u\gamma (u))=\gamma (u)^{2}, 
\]%
and holds for all $u.$ This last identity implies the converse: from $(\dag
),$%
\[
\gamma (u\gamma (u))-\gamma (u)^{2}+\gamma (ue(tu)/t)-2\gamma (u)e(tu)/t=%
\frac{2e(x)}{t^{2}}(1_{\mathbb{A}}-\frac{e(x\circ x)}{2e(x)}). 
\]%
Passage to the limit yields $e(x)/||x||^{2}\rightarrow 0.$

To deduce the cases $k\geq 2$ of $(\times ),$ it again suffices to establish
them for $u$ of norm $1$. Comparing the two sides of $S(x\circ y)=S(x)S(y)$
yields:%
\begin{eqnarray*}
&&\gamma (y\gamma (x))-\gamma (x)\gamma (y)+\gamma (ye(x))-e(x)\gamma
(y)-e(y)\gamma (x) \\
&=&e(x)+e(y)-e(x\circ y).
\end{eqnarray*}%
With $x=tu$ and $y=x\gamma (x)^{k}=t^{k+1}u\gamma (u)^{k},$ the analogous
estimates of ratios against $t^{2}$ are similar and easier because the power
of $t$ in $y$ (strictly) exceeds $2$ and $e(x)=o(x).$ Note that $||x\circ
y||/t\rightarrow 1.$\hfill $\square $

\subsubsection{Proof of Th. 3.2 (Exhaustivity)}

Put $s(x,y)=(s_{1}(x,y),s_{2}(x,y))=(\sigma (x,y),\tau (x,y));$ evaluating
components on right- and left-hand sides, the two of the right are%
\[
RHS_{1}=\sigma (a_{1},a_{2})\sigma (b_{1},b_{2}):\qquad RHS_{2}=\tau
(a_{1},a_{2})\tau (b_{1},b_{2}). 
\]%
Likewise the first component on the left is 
\[
LHS_{1}=\sigma (a+s(a)b)=\sigma (a_{1}+\sigma (a)b_{1},a_{2}+\tau (a)b_{2}). 
\]%
So%
\begin{equation}
\sigma (a_{1}+\sigma (a)b_{1},a_{2}+\tau (a)b_{2})=\sigma
(a_{1},a_{2})\sigma (b_{1},b_{2}).  \tag{$S1$}
\end{equation}%
Similarly,%
\begin{equation}
\tau (a_{1}+\sigma (a)b_{1},a_{2}+\tau (a)b_{2})=\tau (a_{1},a_{2})\tau
(b_{1},b_{2}).  \tag{$S2$}
\end{equation}

In $(S1),$ taking $a_{2}=b_{2}=0$ gives%
\[
\sigma (a_{1}+\sigma (a_{1},0)b_{1},0)=\sigma (a_{1},0)\sigma (b_{1},0). 
\]%
So $\varphi (x):=\sigma (x,0)$ solves the standard $(GS)$ equation%
\[
\varphi (x+y\varphi (x))=\varphi (x)\varphi (y), 
\]%
and so for some $\sigma _{1}$%
\[
\sigma (x,0)=1+\sigma _{1}x. 
\]%
Likewise, working with $(S2)$ with $a_{1}=b_{1}=0$ gives%
\begin{eqnarray*}
\tau (0,a_{2}+\tau (0,a_{2})b_{2}) &=&\tau (0,a_{2})\tau (0,b_{2}): \\
\tau (0,y) &=&1+\tau _{2}y.
\end{eqnarray*}%
Now taking $a_{2}=b_{2}=0$ in $(S2)$ gives%
\[
\tau (a_{1}+\sigma (a_{1},0)b_{1},0)=\tau (a_{1},0)\tau (b_{1},0). 
\]%
So $\varphi (x):=\tau (x,0)$ and $h(x):=\sigma (x,0)$ solve the pexiderized
version of the (scalar) $(GS):$%
\[
\varphi (x+h(x)y)=\varphi (x)\varphi (y), 
\]%
for which see [Jab, Cor. 1(iv)]. Here there are two possibilities%
\begin{eqnarray*}
h &=&1\text{ and }\varphi =e^{\gamma x}, \\
h &=&1+cx\text{ and }\varphi =(1+cx)^{\gamma }.
\end{eqnarray*}%
That is,%
\begin{eqnarray*}
\sigma (x,0) &=&1+\sigma _{1}x=1\text{ and }\tau (x,0)=e^{\gamma x},\text{
i.e. }\sigma _{1}=0. \\
\sigma (x,0) &=&1+\sigma _{1}x=1+cx\text{ and }\tau (x,0)=(1+cx)^{\gamma }.
\end{eqnarray*}%
The latter case should be interpreted as either $\gamma \neq 0$ and $\tau
_{1}=c=\sigma _{1},$ or $\gamma =0$ and so $\tau _{1}=0.$

Also taking $a_{1}=b_{1}=0$ in $(S1)$ gives%
\[
\sigma (0,a_{2}+\tau (0,a_{2})b_{2})=\sigma (0,a_{2})\sigma (0,b_{2}). 
\]%
So $\varphi (x):=\sigma (0,y)$ and $h(x):=\tau (0,y)$ solve the pexiderized
version%
\[
\varphi (x+h(x)y)=\varphi (x)\varphi (y). 
\]

Here again there are two possibilities:%
\begin{eqnarray*}
\tau (0,y) &=&1+\tau _{2}y=1\text{ and }\sigma (0,y)=e^{\gamma y},\text{i.e. 
}\tau _{2}=0; \\
\tau (0,y) &=&1+\tau _{2}y=1+cy\text{ and }\sigma (0,y)=(1+cy)^{\gamma },%
\text{ i.e. }\sigma (0,y)=(1+\tau _{2}y)^{\gamma }.
\end{eqnarray*}

In summary,%
\begin{equation}
\sigma (x,0)=1+\sigma _{1}x,\qquad \tau (0,y)=1+\tau _{2}y, 
\tag{\QTR{rm}{summary}}
\end{equation}%
and the following `side' conditions hold.

\noindent If $\tau _{2}\neq 0,$ then $\sigma _{2}:=\tau _{2}$ and $\sigma
(0,y)=(1+\sigma _{2}y)^{\gamma };$ otherwise $\sigma (0,y)=e^{\gamma y}.$

\noindent If $\sigma _{1}\neq 0,$ then $\tau _{1}:=\sigma _{1}$ and $\tau
(x,0)=(1+\tau _{1}x)^{\delta };$ otherwise $\tau (x,0)=e^{\delta x}.$

\bigskip

Finally in $(S1)$ and in $(S2)$ take $a_{2}=b_{1}=0.$ Then%
\[
\sigma (a_{1},\tau (a_{1},0)b_{2})=\sigma (a_{1},0)\sigma (0,b_{2})\qquad
\tau (a_{1},\tau (a_{1},0)b_{2})=\tau (a_{1},0)\tau (0,b_{2}). 
\]%
Thus%
\begin{eqnarray*}
\tau (x,y) &=&\tau (x,0)(1+\tau _{2}y/\tau (x,0))\text{\quad (taking }y=\tau
(x,0)b_{2}\text{)} \\
&=&\tau (x,0)+\tau _{2}y=\left\{ 
\begin{array}{cc}
(1+\tau _{1}x)^{\delta }+\tau _{2}y & \text{if }\tau _{1}=\sigma _{1}\neq 0
\\ 
e^{\delta x}+\tau _{2}y & \text{if }\sigma _{1}=0.%
\end{array}%
\right.
\end{eqnarray*}

Similarly%
\begin{eqnarray*}
\sigma (a_{1},\tau (a_{1},0)b_{2}) &=&\sigma (a_{1},0)\sigma (0,b_{2}) \\
\sigma (x,\tau (x,0)b_{2}) &=&(1+\sigma _{1}x)\sigma (0,y/\tau (x,0))\text{
(take }y=\tau (x,0)b_{2}\text{)} \\
&=&\left\{ 
\begin{array}{cc}
(1+\sigma _{1}x)(1+\sigma _{2}y/\tau (x,0))^{\gamma } & \text{if }\sigma
_{2}=\tau _{2}\neq 0 \\ 
(1+\sigma _{1}x)e^{\gamma y/\tau (x,0)} & \text{if }\sigma _{2}=0%
\end{array}%
\right. ,
\end{eqnarray*}

\[
=\left\{ 
\begin{array}{cc}
((1+\sigma _{1}x)^{1/\gamma }+\sigma _{2}y(1+\sigma _{1}x)^{[1/\gamma
-\delta ]})^{\gamma } & \text{if }\tau _{1}=\sigma _{1}\neq 0\text{ \& }%
\sigma _{2}=\tau _{2}\neq 0 \\ 
(1+\sigma _{2}ye^{-\delta x})^{\gamma } & \text{if }\sigma _{1}=0\text{ \& }%
\sigma _{2}=\tau _{2}\neq 0 \\ 
e^{\gamma ye^{-\delta x}} & \text{if }\sigma _{1}=\sigma _{2}=0 \\ 
(1+\sigma _{1}x)e^{\gamma y(1+\sigma _{1}x)^{-\delta }} & \text{if }\tau
_{1}=\sigma _{1}\neq 0\text{ \& }\sigma _{2}=0.%
\end{array}%
\right. . 
\]%
We examine all the possible cases. The first and second \textit{special}
cases below emerge as the only viable ones, i.e. leading to a solution. The
remaining \textit{general} cases turn out to be impossible on asymptotic
growth grounds.

\textit{First Special Case}: Take $\sigma _{1}\neq 0\neq \sigma _{2}$ and $%
\gamma =\delta =1;$ then $\sigma (x,0)=\tau (x,0)=1+\sigma _{1}x$ and $%
\sigma (0,y)=\tau (0,y)=1+\sigma _{2}y$ with $:$ 
\begin{eqnarray*}
\sigma (a_{1},(1+\sigma _{1}a_{1})b_{2}) &=&(1+\sigma _{1}a_{1})(1+\sigma
_{2}b_{2})\qquad \text{(}z-b_{2}:=\sigma _{1}a_{1}b_{2}\text{),} \\
\sigma (a_{1},z) &=&1+\sigma _{1}a_{1}+\sigma _{2}b_{2}+\sigma _{1}\sigma
_{2}a_{1}b_{2} \\
&=&1+\sigma _{1}a_{1}+\sigma _{2}b_{2}+\sigma _{2}[z-b_{2}] \\
&=&1+\sigma _{1}a_{1}+\sigma _{2}z.
\end{eqnarray*}

We verify that this formula with $\sigma =\tau $ satisfies $(GS).$ It will
suffice to check the first component:%
\begin{eqnarray*}
LHS_{1} &=&\sigma (a_{1}+\sigma (a)b_{1},a_{2}+\tau (a)b_{2}) \\
&=&1+\rho _{1}(a_{1}+(1+\rho _{1}a_{1}+\rho _{2}a_{2})b_{1})+\rho
_{2}(a_{2}+(1+\rho _{1}a_{1}+\rho _{2}a_{2})b_{2}) \\
&=&1+\rho _{1}a_{1}+\rho _{1}b_{1}+\rho _{1}^{2}a_{1}b_{1}+\rho _{1}\rho
_{2}a_{2}b_{1} \\
&&+\rho _{2}a_{2}+\rho _{2}b_{2}+\rho _{2}\rho _{1}a_{1}b_{2}+\rho
_{2}^{2}a_{2}b_{2} \\
&=&1+\rho _{1}a_{1}+\rho _{2}a_{2}+\rho _{1}b_{1}+\rho _{2}b_{2}+(\rho
_{1}a_{1}+\rho _{2}a_{2})(\rho _{1}b_{1}+\rho _{2}b_{2}).
\end{eqnarray*}%
This agrees with%
\begin{eqnarray*}
RHS_{1} &=&(1+\rho _{1}a_{1}+\rho _{2}a_{2})(1+\rho _{1}b_{1}+\rho _{2}b_{2})
\\
&=&1+\rho _{1}a_{1}+\rho _{2}a_{2}+\rho _{1}b_{1}+\rho _{2}b_{2} \\
&&+(\rho _{1}a_{1}+\rho _{2}a_{2})(\rho _{1}b_{1}+\rho _{2}b_{2}).
\end{eqnarray*}%
The calculation of the second component is exactly the same in these
circumstances.

\textit{Second Special Case.} Take $\tau _{2}\neq 0,$ and $\gamma =0,$
equivalently $\sigma _{2}=0.$ Then $\sigma (x,y)=1+\sigma _{1}x.$ Take $%
\delta =0,$ then $\tau (x,0)=1,$ equivalently $\tau _{1}=0.$ Then $\tau
(x,y)=1+\tau _{2}y.$ We check the proposed solution:%
\[
s(x_{1},x_{2})=(1+\sigma _{1}x_{1},1+\tau _{2}x_{2}). 
\]%
The first component on the left is%
\begin{eqnarray*}
LHS_{1} &=&\sigma (a_{1}+\sigma (a)b_{1},a_{2}+\tau (a)b_{2}) \\
&=&1+\sigma _{1}(a_{1}+b_{1}+\sigma _{1}a_{1}b_{1})=1+\sigma
_{1}a_{1}+\sigma _{1}b_{1}+\sigma _{1}^{2}a_{1}b_{1},
\end{eqnarray*}%
and this matches%
\begin{eqnarray*}
RHS_{1} &=&\sigma (a_{1},a_{2})\sigma (b_{1},b_{2}) \\
&=&(1+\sigma _{1}a_{1})(1+\sigma _{1}b_{1}) \\
&=&1+\sigma _{1}a_{1}+\sigma _{1}b_{1}+\sigma _{1}^{2}a_{1}b_{1}.
\end{eqnarray*}%
The calculation of the second component (with $\tau $) is similar.

It turns out that these special cases are the only possible ones. So it now
remains to eliminate the remaining (ostensibly `general') cases.

\textit{General Cases }(non-viable).

\noindent We first compute $\sigma ,$ according to the side condition below $%
($\textrm{summary}$)$ above, treating it `disjunctively' as Cases A and B.

\noindent Case A (If $\tau _{2}\neq 0$ then $\sigma _{2}=\tau _{2}$ and $%
\sigma (0,y)=(1+\sigma _{2}y)^{\gamma }.).$ Here%
\[
\sigma (\sigma (0,a_{2})b_{1},a_{2})=\sigma (0,a_{2})\sigma (b_{1},0). 
\]%
Substitution under this case yields%
\begin{eqnarray*}
\sigma ((1+\sigma _{2}y)^{\gamma }b_{1},y) &=&(1+\sigma _{2}y)^{\gamma
}(1+\sigma _{1}b_{1}) \\
\sigma (x,y) &=&\sigma _{1}x+(1+\sigma _{2}y)^{\gamma }\text{ with }\sigma
_{2}=\tau _{2}\neq 0.
\end{eqnarray*}%
\noindent Case B ($\sigma (0,y)=e^{\gamma y}).$ Substituting as in Case A%
\begin{eqnarray*}
\sigma (e^{\gamma y}b_{1},y) &=&e^{\gamma y}(1+\sigma _{1}b_{1})=e^{\gamma
y}+\sigma _{1}e^{\gamma y}b_{1} \\
\sigma (x,y) &=&e^{\gamma y}(1+\sigma _{1}b_{1})=\sigma _{1}x+e^{\gamma y}%
\text{ with }\tau _{2}=0.
\end{eqnarray*}

\noindent Next we compute $\tau $ analogously to $\sigma $, again
`disjunctively'.

\noindent Case A (If $\sigma _{1}\neq 0$ then $\sigma _{1}=\tau _{1}$ and $%
\tau (x,0)=(1+\sigma _{1}x)^{\delta }).$ Substituting under this case into 
\[
\tau (a_{1},\tau (a_{1},0)b_{2})=\tau (a_{1},0)\tau (0,b_{2}) 
\]%
yields%
\begin{eqnarray*}
\tau (x,(1+\sigma _{1}x)^{\delta }b_{2}) &=&(1+\sigma _{1}x)^{\delta
}(1+\tau _{2}b_{2}) \\
\tau (x,y) &=&(1+\sigma _{1}x)^{\delta }+\tau _{2}b_{2}(1+\sigma
_{1}x)^{\delta } \\
&=&(1+\tau _{1}x)^{\delta }+\tau _{2}y\text{ with }\tau _{1}:=\sigma
_{1}\neq 0.
\end{eqnarray*}

\noindent Case B ($\tau (x,0)=e^{\delta x}).$ Substituting as in the
preceeding Case A yields:%
\begin{eqnarray*}
\tau (x,e^{\delta x}b_{2}) &=&\tau (x,0)\tau (0,b_{2})=e^{\delta x}(1+\tau
_{2}b_{2}) \\
\tau (x,y) &=&e^{\delta x}+\tau _{2}b_{2}e^{\delta x} \\
&=&e^{\delta x}+\tau _{2}y\text{ with }\tau _{1}=\sigma _{1}=0.
\end{eqnarray*}

In summary, we have the following possibilities for $\sigma $ and $\tau .$ 
\textit{\ }%
\[
\sigma (x_{1},x_{2})=\left\{ 
\begin{array}{cc}
\sigma _{1}x_{1}+(1+\sigma _{2}x_{2})^{\gamma } & \text{with }\sigma
_{2}=\tau _{2}\neq 0, \\ 
\sigma _{1}x_{1}+e^{\gamma x_{2}} & \text{with }\sigma _{2}=\tau _{2}=0,%
\end{array}%
\right. 
\]%
\[
\tau (x_{1},x_{2})=\left\{ 
\begin{array}{cc}
(1+\tau _{1}x_{1})^{\delta }+\tau _{2}x_{2} & \text{with }\tau _{1}:=\sigma
_{1}\neq 0, \\ 
e^{\delta x_{1}}+\tau _{2}x_{2} & \text{ with }\tau _{1}=\sigma _{1}=0.%
\end{array}%
\right. 
\]%
We now rule out all four possible pairings of $\sigma $ and $\tau $ by their
asymptotic behaviour for large values of the arguments on both sides of the
first component equation $(S1).$ This last asserts that%
\[
\sigma (a_{1}+\sigma (a)b_{1},a_{2}+\tau (a)b_{2})=\sigma
(a_{1},a_{2})\sigma (b_{1},b_{2}). 
\]%
Case 1. Consider pairing the first choices available to $\sigma $ and $\tau
. $ 
\[
\sigma (x)=\sigma _{1}x_{1}+(1+\sigma _{2}x_{2})^{\gamma },\tau (x)=(1+\tau
_{1}x_{1})^{\delta }+\tau _{2}x_{2}. 
\]%
On the LHS of $(S1)$ the two arguments for $\sigma $ are 
\[
x_{1}:=a_{1}+b_{1}[\sigma _{1}a_{1}+(1+\sigma _{2}a_{2})^{\gamma }],\quad
x_{2}:=a_{2}+b_{2}(1+\tau _{1}a_{1})^{\delta }+\tau _{2}a_{2}b_{2}. 
\]%
So%
\begin{eqnarray*}
LHS &=&\sigma _{1}x_{1}+(1+\sigma _{2}x_{2})^{\gamma } \\
&=&\sigma _{1}a_{1}+\sigma _{1}\sigma _{1}a_{1}b_{1}+\sigma
_{1}b_{1}(1+\sigma _{2}a_{2})^{\gamma } \\
&&+(1+\sigma _{2}a_{2}+\sigma _{2}\tau _{2}a_{2}b_{2}+\sigma
_{2}b_{2}(1+\tau _{1}a_{1})^{\delta })^{\gamma }.
\end{eqnarray*}%
But%
\begin{eqnarray*}
RHS &=&[\sigma _{1}a_{1}+(1+\sigma _{2}a_{2})^{\gamma }][\sigma
_{1}b_{1}+(1+\sigma _{2}b_{2})^{\gamma }] \\
&=&\sigma _{1}a_{1}\sigma _{1}b_{1}+\sigma _{1}a_{1}(1+\sigma
_{2}b_{2})^{\gamma }+\sigma _{1}b_{1}(1+\sigma _{2}a_{2})^{\gamma } \\
&&+(1+\sigma _{2}a_{2})^{\gamma }(1+\sigma _{2}b_{2})^{\gamma }.
\end{eqnarray*}%
Equating sides and setting $b_{1}=a_{2}=0$ gives%
\[
\sigma _{1}a_{1}+(1+\sigma _{2}b_{2}(1+\tau _{1}a_{1})^{\delta })^{\gamma
}=\sigma _{1}a_{1}(1+\sigma _{2}b_{2})^{\gamma }+(1+\sigma
_{2}b_{2})^{\gamma }. 
\]%
Letting $a_{1},b_{2}\rightarrow \infty $ yields%
\[
\sigma _{1}a_{1}+(\sigma _{2}b_{2})^{\gamma }(1+\tau _{1}a_{1})^{\delta
+\gamma }\sim \sigma _{1}a_{1}(\sigma _{2}b_{2})^{\gamma }+(\sigma
_{2}b_{2})^{\gamma }, 
\]%
a contradiction unless $\sigma _{1}=\tau _{1}=0$. This yields $\sigma
=(1+\sigma _{2}x_{2})^{\gamma }$ and $\tau =1+\tau _{2}x_{2}$ in the
univariate format. Clearly $(S2)$ is satisfied, whereas $(S1)$ requires that%
\[
\lbrack 1+\sigma _{2}(x_{2}+y_{2}(1+\tau _{2}x_{2}))]^{\gamma }=[(1+\sigma
_{2}x_{2})(1+\sigma _{2}y_{2})]^{\gamma }. 
\]%
The case $\gamma =0$ gives the `independent' format with $\sigma =1$, as
does $\sigma _{2}=0$. Otherwise one has $\sigma _{2}=\tau _{2},$ a \textit{%
univariate type}.

Case 2. Consider now pairing second choices, so that%
\[
\sigma =\sigma _{1}x_{1}+e^{\gamma x_{2}},\quad \tau =e^{\delta x_{1}}+\tau
_{2}x_{2}. 
\]%
Here the arguments on the left of $(S1)$ are%
\[
x_{1}=a_{1}+b_{1}[\sigma _{1}a_{1}+e^{\gamma a_{2}}],\quad
x_{2}=a_{2}+b_{2}[e^{\delta a_{1}}+\tau _{2}a_{2}]. 
\]%
So%
\[
LHS=\sigma _{1}x_{1}+e^{\gamma x_{2}}=\sigma _{1}a_{1}+\sigma
_{1}b_{1}\sigma _{1}a_{1}+\sigma _{1}b_{1}e^{\gamma a_{2}}+\exp [\gamma
a_{2}+\gamma b_{2}[e^{\delta a_{1}}+\tau _{2}a_{2}]]. 
\]%
But%
\[
RHS=[\sigma _{1}a_{1}+e^{\gamma a_{2}}][\sigma _{1}b_{1}+e^{\gamma
b_{2}}]=\sigma _{1}a_{1}\sigma _{1}b_{1}+\sigma _{1}b_{1}e^{\gamma
a_{2}}+\sigma _{1}a_{1}e^{\gamma b_{2}}+e^{\gamma (a_{2}+b_{2})}. 
\]%
Again the asymptotic behaviours of both sides do not match unless $\sigma
_{1}=\gamma =0$. when both sides reduce to 1. So here $\sigma =1$ and so by $%
(S2)$%
\begin{eqnarray*}
e^{\delta (x_{1}+y_{1})}+\tau _{2}(x_{2}+y_{2}(e^{\delta x_{1}}+\tau
_{2}x_{2})) &=&(e^{\delta x_{1}}+\tau _{2}x_{2})(e^{\delta y_{1}}+\tau
_{2}y_{2}) \\
e^{\delta (x_{1}+y_{1})}+\tau _{2}x_{2}+\tau _{2}y_{2}e^{\delta x_{1}}+\tau
_{2}x_{2}\tau _{2}y_{2} &=&e^{\delta (x_{1}+y_{1})}+\tau _{2}x_{2}e^{\delta
y_{1}}+\tau _{2}y_{2}e^{\delta x_{1}}+\tau _{2}x_{2}\tau _{2}y_{2} \\
\tau _{2}x_{2} &=&\tau _{2}x_{2}e^{\delta y_{1}}
\end{eqnarray*}%
leading to the solutions $\sigma =1,\tau =1+\tau _{2}x_{2}$ covered by the
`independent' format and finally the `univariate' type%
\[
\sigma =1,\tau =e^{\delta x_{1}}. 
\]

The remaining Cases 3 and 4 (cross-choices) are similar. \hfill $\square $

\subsubsection{Th. 3.3 -- Stone-Weierstrass argument}

Using the Structure Theorem of the Euclidean case of $\mathbb{G}_{S}^{\ast }(%
\mathbb{R}^{d})$ established in \S 3, we may now describe $S$ by reference
to $C(T)$ as follows.

Denote by $\delta $ and $\delta _{T}$ the Dirac mass function respectively
for $[0,1]^{2}$ and for $T^{2}$ so that $\delta (t)(s)=1$ iff $t=s;$ then by
Th. 3.1, there are numbers $\sigma _{ij}$ with%
\[
\sigma _{T}(x_{T})(t_{i})=\Sigma _{j}\sigma _{ij}x_{T}(t_{j})\text{ so that }%
\sigma _{i}{}_{j}=\sigma _{T}(t_{i})_{j}=S(\delta _{T}(t_{j}))(t_{i})-1. 
\]%
The corresponding `generator' will be represented here by $\rho
_{T}(x_{T})\in \mathbb{R}^{n+1}=C(T)$ with 
\[
\rho _{T}(x_{T})(t_{i}):=\Sigma _{j}\sigma _{ij}x_{T}(t_{j})=\Sigma
_{i=0}^{n}\sigma _{T}(x_{T})(t_{i})\rightarrow S(x^{T})(t_{i}), 
\]%
with the limit here again under refinement of subdivisions.

Define the partition $\mathcal{P}_{T}$ of $T$ to comprise all the distinct
sets $I_{T}(t)$ for $t\in T$ with 
\[
I_{T}(t):=\{t_{i}\in T:(\forall x)S(x)(t_{i})=S(x)(t)\}. 
\]%
Then for $s,t\in I_{T}(t)$ 
\[
S(\delta (t_{j}))(s)=S(\delta (t_{j}))(t)\qquad (j=0,1,...,n). 
\]%
Here taking limits under refinement of subdivisions $T$ yields 
\[
I_{T}(t)\rightarrow K_{t}. 
\]%
For $s,t\in I\in \mathcal{P}_{T}$, as $S(\delta (t_{j}))(t)=0$ iff $S(\delta
(t_{j}))(s)=0,$ we may define 
\[
J_{T}(I):=\{t_{j}\in T:S(\delta (t_{j}))(t)\neq 0,t\in I\}. 
\]%
Then for $x=x_{T}$ 
\[
\sigma _{I}(x):=\rho _{T}(e_{I}\cdot x)e_{I}\text{ for }e_{I}:=\Sigma _{t\in
J(I)}\delta _{T}(t)\text{ and }I\in \mathcal{P}_{T}. 
\]%
Here $x\mapsto e_{I}\cdot x$ is the projection from $\mathbb{R}^{n+1}$ onto
the span of $\{\delta _{T}(t):t\in J(I)\}.$

Taking limits over subdivisions $T$ under refinement with inclusion of $t$
yields%
\[
\rho _{T}(x_{T})(t)\rightarrow \rho (x)(t); 
\]%
this limit generator gives a continuous linear map $\rho :C[0,1]\rightarrow
C[0,1]$. For $s\in \lbrack 0,1],$ put $\rho _{s}(x):=\rho (x)(s),$ a
continuous linear functional. Then 
\[
\{t:\rho _{s}=\rho _{t}\}=K_{s}. 
\]%
Let $e_{K}$ denote the map $x\mapsto e_{K}\cdot x$ projecting from $C[0,1]$
onto $C(K)$. Now define $\sigma _{K}:C(K)\rightarrow \mathbb{R}$ for $%
K_{t}=K\in \mathcal{P}$ by 
\[
\sigma _{K}(x):=\rho _{t}(e_{K}\cdot x), 
\]%
thereby completing the analysis of the action of $S$. \hfill $\square $

\subsubsection{Proof of Lemma 4.1\textbf{\ }}

With restrictions on $a,b,c,z$ as above, using $(GS)$ and invertibility of $%
S,$%
\begin{eqnarray*}
S(c+(a-b)) &=&S(a+(c-b)S(a)S(b)^{-1})=S(a)S((c-b)S(b)^{-1}) \\
&=&S(b)S((c-b)S(b)^{-1})=S(b+S(b)(c-b)S(b)^{-1})=S(c).
\end{eqnarray*}%
Since%
\[
S(z+S(z)a)=S(z)S(a)=S(z)S(b)=S(z+S(z)b)\in \mathbb{A}^{-1}, 
\]%
replacing in (i) $a$ by $z+S(z)a$ and $b$ by $z+S(z)b$ yields (ii). In
particular, for any $z\in \mathbb{G}_{S}^{\ast }$ 
\[
1_{\mathbb{A}}=S(0)=S(S(z)a) 
\]%
(take $b=c=0),$ giving $S(\mathbb{G}_{S}^{\ast })\mathcal{N}\subseteq 
\mathcal{N}$ and so also the last assertion, as $S(0)\mathcal{N}=\mathcal{N}$%
. \hfill $\square $

\subsubsection{Proof of Corollary 4.4}

Note first that since $\mathbb{C}_{1}\mathbb{=C}\backslash \{0\}=\mathbb{C}%
^{-1},$ $\mathbb{G}_{S}^{\ast }(\mathbb{C})=\mathbb{G}_{S}(\mathbb{C}).$ It
is known [Bar] and reproved in Prop. 4.2 that in the case $\mathbb{A=C}$ any
continuous solution $S$ of $(GS)$ is either canonical, i.e. of the form $%
S(z)=1+\rho z$ so that $\mathbb{G}_{S}(\mathbb{C})=\mathbb{G}_{\rho }(%
\mathbb{C}),$ or takes the `non-canonical' form $S(z):=1+a\func{Re}(z)+b%
\func{Im}(z)$ with $a,b$ real constants. Here if $a=b=0,$ this leads to $%
\mathrm{ran}S=\{1\};$ then $1$ is an isolated point, and the preceeding
result cannot be applied, although its conclusion still holds with $c=0.$
The alternative is that $\mathrm{ran}S=\langle 1\rangle .$ This presents two
possibilities:

\noindent (i) $S(z)=1+az$ for $z\in \langle 1\rangle $ with $a\neq 0,$
yielding a Popa subgroup $\mathbb{G}_{\alpha }^{\ast }(\langle 1\rangle )=%
\mathbb{G}_{\alpha }^{\ast }(\mathbb{R})$ with $\alpha =1/a;$

\noindent (ii) $S(z)=1-ibz$ for $z\in \langle i\rangle $ with $b\neq 0,$
yielding a Popa subgroup $\mathbb{G}_{\beta }^{\ast }(\langle i\rangle
)\approx \mathbb{G}_{|\beta |}^{\ast }(\mathbb{R})$ with $\beta =1/(-ib).$

These two separate restrictions of $S,$ both in `canonical' form, correspond
to $S((z-1)/a)=z$ for $z\in \langle 1\rangle $ (with $c=1/a),$ and $%
S((iz-1)/(-ib))=iz$ for $z\in \langle i\rangle $ (with $c=i/b).$

Write $z_{1}=u+iv$ and $z_{2}=x+iy;$ then, since $S(z_{1})=1+au+bv,$%
\[
z_{1}\circ _{S}z_{2}=(u+iv)+(x+iy)+(1+au+bv)(x+iy). 
\]%
This corresponds to a Popa operation $\circ _{\sigma }$on the set $\mathbb{G}%
_{\alpha }^{\ast }(\mathbb{R})\times \mathbb{G}_{\text{\TEXTsymbol{\vert}}%
\beta |}^{\ast }(\mathbb{R})$ with $\sigma (u,v):=1+au+bv$ and 
\[
(u,v)\circ _{\sigma }(x,y)=(u+x+\sigma (u,v)x,v+y+\sigma (u,v)y). 
\]%
We shall identified this in \S 3 as $\mathbb{G}_{\sigma }^{\ast }(\mathbb{R}%
^{2})$.\hfill $\square $

\subsubsection{Example 5.1 (Standardized Tilting in $\mathbb{C}$)}

Writing $\omega =x+iy,$ the real and imaginary parts give:%
\begin{eqnarray*}
e^{x}\cos y &=&1+x,\qquad e^{x}\sin y=y: \\
y^{2} &=&y^{2}(x):=e^{2x}-(1+x)^{2}>0,\text{ for }x>0.
\end{eqnarray*}%
For $x\geq 0,$ $y(x)$ is monotonic and unbounded, whereas%
\[
x\mapsto e^{-x}y(x)=\sqrt{1-(e^{-x}(1+x))^{2}} 
\]%
increases on $[0,\infty )$ strictly from $0$ to $1,$ yielding solutions in $%
x $ to the equation%
\[
\sin y(x)=e^{-x}y(x), 
\]%
one within each consecutive interval in which $\sin y(x)$ traces the
interval $[-1,1]$ with $x$ near the solution point of $y(x)=\pi /2$ mod $%
2\pi \mathbb{Z}$. It follows that, in alternate intervals where $\cos
y(x)>0, $ 
\[
e^{x}\cos y(x)=\sqrt{e^{2x}-y(x)^{2}}=1+x, 
\]%
thus satisfying $(ST)$.

\noindent \textbf{Remark. }The function $y(x)$ is also defined in a maximal
interval $[-\xi ,0],$ with $\xi $ satisfying $e^{-2\xi }=(1-\xi )^{2},$
equivalently $e^{-\xi }=\xi -1.$ From here $\xi >1,$ and in fact 
\[
\xi =1.27846...; 
\]%
it follows that $y(-\xi )=0.$ (Evidently $y(-1)=1.)$ So $\sin y(-\xi
)=e^{\xi }y(-\xi )=0;$ however, 
\[
e^{-\xi }\cos y(-\xi )=e^{-\xi }=\xi -1\neq 1-\xi , 
\]%
since $\xi \neq 1.$

\subsubsection{\textbf{Proof of Prop. 5.4}}

Using orthogonality,%
\begin{eqnarray*}
RHS_{GS} &=&(1_{\mathbb{A}}+\Sigma _{i}\sigma (e_{i}x)e_{i})(1_{\mathbb{A}%
}+\Sigma _{j}\sigma (e_{i}y)e_{j}) \\
&=&1_{\mathbb{A}}+\Sigma _{i}\sigma (e_{i}x)e_{i}+\Sigma _{j}\sigma
(e_{i}y)e_{j} \\
&&+\Sigma _{i}\sigma (e_{i}x)\sigma (e_{i}y)e_{i}.
\end{eqnarray*}%
Noting that%
\[
x+[1_{\mathbb{A}}+\Sigma _{i}\sigma (e_{i}x)e_{i}]y=x+y+\Sigma _{i}\sigma
(e_{i}x)e_{i}y, 
\]%
we compute, using orthogonality, that%
\begin{eqnarray*}
\sigma (e_{j}[x+y+\Sigma _{i}\sigma (e_{i}x)e_{i}y]) &=&\sigma
(e_{j}x+e_{j}y+\Sigma _{i}\sigma (e_{i}x)e_{j}e_{i}y) \\
&=&\sigma (e_{j}x)+\sigma (e_{j}y)+\sigma \lbrack \sigma (e_{j}x)e_{j}y] \\
&=&\sigma (e_{j}x)+\sigma (e_{j}y)+\sigma (e_{j}x)\sigma \lbrack e_{j}y].
\end{eqnarray*}%
So%
\begin{eqnarray*}
LHS_{GS} &=&1_{\mathbb{A}}+\Sigma _{j}\sigma (e_{j}[x+y+\Sigma _{i}\sigma
(e_{i}x)e_{i}y])e_{j} \\
&=&1_{\mathbb{A}}+\Sigma _{j}\sigma (e_{j}x)e_{j}+\sigma
(e_{j}y)e_{j}+\sigma (e_{j}x)\sigma (e_{j}y)e_{j},
\end{eqnarray*}%
and the two sides match. \hfill $\square $

\subsubsection{Proof of Lemma 5.5}

W.l.o.g. $||e^{-a}||<1.$ Then $||e^{a}||>1,$ for otherwise $||e^{a}||\leq 1,$
leading to%
\[
1=||e^{-a}e^{a}||\leq ||e^{-a}||\cdot ||e^{a}||<1, 
\]%
a contradiction. Furthermore, for $n\in \mathbb{Z}$, 
\[
||e^{na}||=||e^{(n+1)a}e^{-a}||\leq ||e^{-a}||\cdot
||e^{(n+1)a}||<||e^{(n+1)a}||. 
\]%
The sequence $\{||e^{na}||\}_{n\in \mathbb{N}}$ is thus monotonically
increasing. If $||e^{na}||\rightarrow c$ for some finite $c>0,$ then, by the
preceding inequality%
\[
c\leq ||e^{-a}||c<c, 
\]%
also a contradiction. So $||e^{na}||$ is unbounded. Similarly, $%
\{||e^{-na}||\}_{n\in \mathbb{N}}$ is monotonically decreasing, this time
with limit $0.$

As for the final statement, provided $||e^{\gamma (u)}||\neq 1$,%
\[
T(\pm us)=u(e^{\pm s\gamma (u)}-1)/\gamma (u) 
\]%
is unbounded as $s\rightarrow \infty $ in one of the directions $\pm u$%
.\hfill $\square $

\subsubsection{\textbf{Proof of Th. 6.1 (Wo\l od\'{z}ko-Javor Theorem)}}

\textbf{\ }Given $S,$ take $\mathcal{N}:=\{u\in \mathbb{G}_{S}^{\ast
}:S(u)=1_{A}\}$, $\Lambda :=S\mathcal{(}\mathbb{G}_{S}^{\ast }),$ and choose
any right inverse $W$ with $S(W(a))\equiv a.$ Then use Lemma 2.1(i) and (ii).

For the reverse direction,\textbf{\ }as $\mathcal{N}$ is a subgroup, we may
work $\func{mod}\mathcal{N},$ indicating this now with $\equiv _{\mathcal{N}%
}.$ First note that $S$ is well-defined. For if $W(\lambda _{1})=W(\lambda
), $ then, taking $\lambda _{2}=\lambda \lambda _{1}^{-1},$ so that $\lambda
=\lambda _{1}\lambda _{2}$ by (iii),%
\[
\lambda _{1}W(\lambda _{2})=\lambda _{1}W(\lambda _{2})+[W(\lambda
_{1})-W(\lambda _{1}\lambda _{2})]\equiv _{\mathcal{N}}0. 
\]%
So $\lambda _{1}W(\lambda _{2})\in \mathcal{N}$, or $W(\lambda _{2})\in
\lambda _{1}^{-1}\mathcal{N}\subseteq \Lambda \mathcal{N=N}$; so by (ii), $%
\lambda _{2}=1,$ i.e. $\lambda =\lambda _{1},$ as required.

We check that (iv) satisfies $(GS)$ with $x=x_{1}$ and $y=x_{2}.$

If $S(x_{1})=0,$ then $(GS)$ holds trivially.

If $S(x_{1})\neq 0$ and $S(x_{1}+S(x_{1})x_{2})\neq 0,$ then pick $\lambda
_{1}$ and $\lambda $ with $W(\lambda _{1})=x_{1}$ and $W(\lambda )\equiv
x_{1}+S(x_{1})x_{2}.$ Take $\lambda _{2}=\lambda \lambda _{1}^{-1};$ then 
\begin{eqnarray*}
W(\lambda _{1}\lambda _{2}) &=&W(\lambda )\equiv _{\mathcal{N}%
}x_{1}+S(x_{1})x_{2}=W(\lambda _{1})+\lambda _{1}x_{2}: \\
\lambda _{1}x_{2} &\equiv &_{\mathcal{N}}W(\lambda _{1}\lambda
_{2})-W(\lambda _{1})\equiv _{\mathcal{N}}\lambda _{1}W(\lambda _{2}):\qquad
x_{2}\equiv _{\mathcal{N}}W(\lambda _{2}).
\end{eqnarray*}%
So here $S(x_{2})=\lambda _{2}\neq 0$, i.e. passing to the contrapositive:
if $S(x_{2})=0,$ then $S(x_{1}+S(x_{1})x_{2})=0$ and $(GS)$ holds.

Now consider $x_{i}$ with both $S(x_{i})\neq 0.$ Write $x_{i}\equiv _{%
\mathcal{N}}W(\lambda _{i}).$ Then%
\begin{eqnarray*}
x_{1}+S(x_{1})x_{2} &\equiv &_{\mathcal{N}}W(\lambda _{1})+\lambda
_{1}W(\lambda _{2})\equiv _{\mathcal{N}}W(\lambda _{1}\lambda _{2}): \\
S(x_{1}+S(x_{1})x_{2}) &=&\lambda _{1}\lambda _{2}=S(x_{1})S(x_{2}).
\end{eqnarray*}%
This completes the check that $(GS)$ holds.\hfill $\square $

\subsection{Proof of Theorem S}

We are grateful to Amol Sasane for the following proof.

For fixed $\theta \in M_{\mathbb{A}}$, the maximal ideal space of $\mathbb{A}
$ (viewed as comprising characters in $\mathcal{L}(\mathbb{A},\mathbb{C})),$
define for $\zeta \in \mathcal{D}_{u}$%
\begin{eqnarray*}
\varphi _{\theta u}(\zeta ) &:&=\theta (f(\zeta \gamma (u))), \\
\psi _{\theta u}(\zeta ) &:&=\theta (\gamma (f(\zeta \gamma (u))\cdot
u/\gamma (u)));
\end{eqnarray*}%
both are holomorphic on $\mathcal{D}_{u}$, since $\theta $ is linear.
Furthermore, by the assumed identity,%
\[
\varphi _{\theta u}(\zeta )=\psi _{\theta u}(\zeta )\qquad \hfill (\zeta \in
\Sigma ).
\]%
So, by Riemann's Uniqueness (Identity) theorem [Gam2,V.7], [Rem, Ch. 8],
also $\varphi _{\theta u}=\psi _{\theta u}$ on $\mathcal{D}_{u}.$ As $\theta 
$ was arbitrary, this may be restated using Gelfand transforms as%
\[
f(\zeta \gamma (u)\;\widehat{}=\gamma (f(\zeta \gamma (u)\cdot u/\gamma
(u))\;\widehat{}\;.
\]%
For $\mathbb{A}$ semisimple, the Gelfand transform is injective and so%
\[
f(\zeta \gamma (u))=\gamma (f(\zeta \gamma (u)\cdot u/\gamma (u))\qquad
\hfill (\zeta \in \mathcal{D}_{u}),
\]%
thus extending an identity from $\Sigma $ to $\mathcal{D}_{u}.$ \hfill $%
\square $

\newpage

\end{document}